\documentclass[a4paper]{article}

\usepackage[utf8]{inputenc}
\usepackage[T1]{fontenc}
\usepackage{textcomp}

\usepackage[top=2.5cm, bottom=2.5cm, outer=2.5cm, inner=2.5cm, heightrounded, marginparwidth=3cm, marginparsep=0.5cm]{geometry}

\usepackage{fancyhdr}
\usepackage{hyperref}
\hypersetup{
    colorlinks,
    linkcolor={black},
    citecolor={blue},
    urlcolor={blue!80!black}
}

\usepackage[nottoc]{tocbibind}

\usepackage{graphicx}
\usepackage{float}
\usepackage[usenames,dvipsnames,svgnames]{xcolor}

\usepackage{mathtools, amsthm, amssymb}
\mathtoolsset{showonlyrefs} 
\usepackage{aligned-overset}
\usepackage{mathrsfs}
\usepackage{enumitem}
\usepackage{ esint }

\newcommand\R{\ensuremath{\mathbb{R}}}
\newcommand\Z{\ensuremath{\mathbb{Z}}}

\newcommand\Q{\ensuremath{\mathbb{Q}}}
\newcommand\C{\ensuremath{\mathbb{C}}}

\let\epsilon\varepsilon

\newcommand\lth{\ensuremath{\mathscr{L}}}
\newcommand\llm{\ensuremath{\ell}}
\newcommand\vek[1]{\ensuremath{\mathbf{#1}}}

\let\originalleft\left
\let\originalright\right
\renewcommand{\left}{\mathopen{}\mathclose\bgroup\originalleft}
\renewcommand{\right}{\aftergroup\egroup\originalright}

\usepackage{tikz}
\usepackage{tikz-cd}

\newtheoremstyle{boldremark}
    {\dimexpr\topsep/2\relax} 
    {\dimexpr\topsep/2\relax} 
    {}          
    {}          
    {\bfseries} 
    {.}         
    {.5em}      
    {}          

\theoremstyle{plain}
\newtheorem{theorem}{Theorem}[section]
\newtheorem{lemma}[theorem]{Lemma}
\newtheorem{prop}[theorem]{Proposition}
\newtheorem{corollary}[theorem]{Corollary}

\theoremstyle{definition}
\newtheorem{definition}[theorem]{Definition}

\theoremstyle{boldremark}
\newtheorem{remark}[theorem]{Remark}
\newtheorem{notation}[theorem]{Notation}

\newenvironment{myproof}[1][\proofname]{%
  \proof[\rm \bf #1]%
}{\endproof}

\usepackage{import}
\usepackage{xifthen}

\usepackage{pdfpages}
\usepackage{transparent}

\usepackage{subcaption}
\usepackage{multirow}
\def\block(#1,#2)#3{\multicolumn{#2}{c}{\multirow{#1}{*}{$ #3 $}}}

\newcommand{\RN}[1]{%
  \textup{\uppercase\expandafter{\romannumeral#1}}%
}

\usepackage[final]{showlabels}

\usepackage{pdfpages}
\usepackage{lipsum}
\usepackage{parskip}
\usepackage{titletoc}

\numberwithin{figure}{section}
\numberwithin{equation}{section}

\usepackage{mathpazo}
\usepackage{bm}

\title{Relaxation of perturbed circles in flat spaces for the Mullins-Sekerka evolution in two dimensions.}
\author{Sa\v{s}a Luki\'{c}\\
        RWTH Aachen University\\
        lukic@eddy.rwth-aachen.de}
\date{}

\begin{document}
    \pagenumbering{arabic}
    \maketitle
    \begin{abstract}
    We analyze the convergence of a perturbed circular interface for the two-phase Mullins-Sekerka evolution in flat two-dimensional space. Our method is based on the gradient flow structure of the evolution and captures two distinct regimes of the dynamics, an initial---and novel---phase of algebraic-in-time decay and a later---and previously explored---phase of exponential-in-time decay. By quantifying the initial phase of relaxation, our method allows for the investigation of systems with large initial dissipation as long as the isoperimetric deficit is small enough. We include quantitative estimates of the solution in terms of its initial data, including the $C^{1}$-distance to the center manifold of circles and the displacement of the barycenter.
\end{abstract}
    \section{Introduction}

The (two-phase) Mullins-Sekerka (MS) model is a free boundary problem arising in various physical contexts such as phase separation and coarsening. It was originally introduced in the seminal work \cite{[Mullins1963]} in the context of crystallization processes. In \cite{[Pego1989]}, the author formally established the MS model as a sharp interface limit of the Cahn-Hilliard phase field model by means of matched asymptotic expansions. Later, this convergence was shown rigorously in \cite{[Alikakos1994]}.

In its strong form, the MS model describes the curvature-driven evolution of interfaces that can be represented by closed, separating but not necessarily connected hypersurfaces. Unlike mean curvature flow or its two-dimensional variant curve shortening flow, which are  local second-order parabolic flows, the MS evolution is a third-order, nonlocal geometric evolution law. While curve shortening flow preserves convexity of the evolving curve, cf. \cite{[Gage1986]}, MS flow does not preserve convexity during the evolution, cf. \cite{[Mayer1997]} for a two-dimensional example configuration. 

Due to the potential presence of finite-time singularities in the flow, a vast collection of weak solution concepts and existence results exists in the literature, such as in \cite{[Luckhaus1995]}, \cite{[Chen1996]} and \cite{[Hensel2024]}. In \cite{[Julin2022]} and \cite{[Arya2025]}, for instance, the relaxation behavior of so-called \emph{flat flows}, limiting time-continuous evolutions based on the minimizing movement approach from \cite{[Luckhaus1995]}, for the MS model in two space dimensions was studied, starting merely from sets of finite perimeter. Short-time existence results of classical solutions for generic, sufficiently regular initially bounded interfaces were proved in \cite{[EscherSimonett1998]} and \cite{[Chen1993]}, while those for unbounded interfaces in two ambient dimensions were proved in \cite{[Escher2024]}. Recently, a weak-strong uniqueness result was established in \cite{[Fischer2024]} relating weak and classical solutions up to the occurrence of singularities. 

For initial perturbations of equilibrium configurations, relaxation rates to these equilibrium configurations were studied in \cite{[Chen1993]} (two dimensions) and \cite{[EscherSimonett1998]} (arbitrary dimensions) for bounded domains. In \cite{[Chugreeva2018]} (two dimensions) and \cite{[Otto2025]} (two and three dimensions) unbounded domains were considered. Our work bridges the gap between bounded and unbounded domains in the two-dimensional case, as it transfers the framework employed in \cite{[Chugreeva2018]} for unbounded domains to bounded ones. While the compactness of the parameter domain provides a helpful energy-dissipation relation compared with the unbounded case, cf. \eqref{eq:energydissspoiler}, the nontrivial curvature of the equilibrium configurations calls for a more refined analysis of the relaxation rates.

In a two-dimensional manifold $\mathcal{M}$, the classical MS model describes the evolution of a sufficiently smooth curve $\Gamma(t) \subset \mathcal{M}$ on some time interval $[0,T]$. In our case, $\mathcal{M}$ is either the Euclidean plane $\R^2$ or the flat torus of edge length $2L$ for some arbitrary $L \in \R_{>0}$. Sufficiently smooth in this context requires at least the existence of a global normal vector field on $\Gamma(t)$ for each $t \in [0,T]$.
Suppose the map
\begin{equation}
    \gamma: \bigsqcup_{t \in [0,T]} \{t\} \times I(t)  \to \mathcal{M}
    \label{eq:curvefamilygeneral}
\end{equation}
represents a family of parametrizations for $\left(\Gamma(t)\right)_{t \in [0,T]}$. Let subsequently $t \in  [0,T]$ be arbitrary but fixed. By assumption, the complement of $\Gamma(t)$ in $\mathcal{M}$ decomposes into two connected components $\Omega_{\text{in}}(t)$ and $\Omega_{\text{out}}(t)$. The labelling of these sets will be made more precise below after the introduction of parametrizations for $\Gamma(t)$. We always choose the normal vector $n(t,s)$ for each $s \in I(t)$ to be outward-pointing for $\Omega_{\text{in}}(t)$, cf. Figure \ref{fig:angles}, in accordance with \cite{[EscherSimonett1998]}.

In the MS model, the curve $\Gamma(t)$ evolves in the normal direction with a normal velocity $V$ determined by
\begin{equation}
    V(t,s) =  \frac{\partial u_{\text{out}}(t,s)}{\partial n(t,s)} - \frac{\partial u_{\text{in}}(t,s)}{\partial n(t,s)}, \quad \forall s \in I(t),
    \label{eq:msnormal}
\end{equation}
where the $u_{\text{in}}(t,\cdot)$ and $u_{\text{out}}(t,\cdot)$ are the respective solutions to the boundary value problems
\begin{equation}
\begin{aligned}
    \Delta u_i(t,\cdot) &= 0 \quad &&\text{ on } \Omega_i(t) \quad &&\text{ for } i \in \{\text{in}, \text{out}\} \\
    u_i(t,\cdot) &= \kappa(t,\cdot) \quad &&\text{ on } \Gamma(t) \quad &&\text{ for } i \in \{\text{in}, \text{out}\}. \label{eq:harmonicproblem}
\end{aligned}
\end{equation}
Here, $\kappa(t,\cdot)$ is the curvature of the interface $\Gamma(t)$ and $\Omega_{\text{in}}(t)$ and $\Omega_{\text{out}}(t)$ denote the two connected components of $\mathcal{M} \setminus \Gamma(t)$. The solutions $u_{\text{in}}(t,\cdot)$ and $u_{\text{out}}(t,\cdot)$ of \eqref{eq:harmonicproblem} are also referred to as the \emph{harmonic extensions of curvature}. Gluing them along $\Gamma(t)$ defines 
\begin{equation}
\begin{aligned}
    u(t,\cdot): \mathcal{M} \to \R, \, x \to \begin{cases}
        u_{i}(x) \quad &\text{ if } x \in \Omega_i(t), \quad \text{ for } i \in \{\text{in}, \text{out}\} \\
        \kappa(x) \quad &\text{ if } x \in \Gamma(t).
    \end{cases}
\end{aligned}
\end{equation}
Physically, $u$ represents the \emph{chemical potential} of the separating phases, which drive the diffusive separation of the phases until a uniform distribution is achieved. The condition
\begin{equation}
     \left| \nabla u_{\text{out}}(t, x) \right| = o( \left| x \right|) \quad \text{ for } x \in \R^2 \setminus \Gamma(t) \to \infty
\end{equation}
is imposed in the case $\mathcal{M} = \R^2$ to enforce the uniqueness of the solution. Notice that any circle of radius $R$ is an equilibrium configuration under the MS evolution, as the system \eqref{eq:harmonicproblem} is satisfied by constant extensions of curvature with value $1/R$ in this case.

In this work, we focus on a special class of interfaces and investigate their convergence to circles including estimates for their relaxation rates. We do not study the existence of solutions to the MS problem, but instead assume the existence of a global in time solution $\left( \Gamma(t) \right)_{t \in [0,\infty)}$. Under the assumption that the initial interface $\Gamma(0)$ admits a polar parametrization, we establish the preservation of the polar graph structure throughout time, i.e., the interface $\Gamma(t)$ admits a polar parametrization of the form
\begin{equation}
    \gamma_{\text{pol}}(t,\cdot): [0, 2 \pi] \to \Gamma(t) \subset \mathcal{M}, \, \phi \mapsto \vek{c}(t) + \rho(t,\phi) \begin{pmatrix} \cos(\phi) \\ \sin(\phi) \end{pmatrix}
    \label{eq:polarparamcurve}
\end{equation}
for each $t \in  [0, \infty)$. The graphing function
\begin{equation}
    \rho(t,\cdot): [0, 2 \pi] \to \R_{>0}
\end{equation}
is shown to remain close to a sphere in the $C^{1}([0, 2 \pi],\R)$-norm, cf. Definition \ref{def:nearlycircular}, if this holds for the initial graphing function $\rho(0, \cdot )$. The center coordinates
\begin{equation}
    \vek{c}: [0,T] \to \mathcal{M}
\end{equation}
in \eqref{eq:polarparamcurve} are precisely the coordinates of the barycenter of $\Omega_{\text{in}}(t)$ at time $t$, cf. \eqref{eq:bulkbarycenter}. Without loss of generality, we can assume $\vek{c}(0) = 0$.

Two major defining characteristics of the flow are the facts that the perimeter $\lth(\Gamma(t))$ is non-increasing and that the area $  \left| \Omega_{\text{in}} \right| $ is preserved throughout the evolution, both of which are proved in Lemma \ref{lemma:geometricvariations}. These are characteristics of the flow shared with the area-preserving mean curvature flow and the surface diffusion flow, which can be realized as $L^{2}$- and $\dot{H}^{1}$-gradient flows of the perimeter functional of curves on $\mathcal{M}$ on some suitable Banach manifold of curves. Thus, it is not surprising that this system can be realized as an $\dot{H}^{-1}$- or $\dot{H}^{-1/2}$-gradient flow of the perimeter functional. The former structure equips the bulk with a Riemannian metric, while the latter restricts the metric to the interface; see, for instance, the survey article \cite{[Garcke2005]} or \cite{[Hensel2024]}. The area-preservation of the flow is reflected by
\begin{equation}
     \left| \Omega_{\text{in}}(t)   \right| = \text{const} = \pi R^2
\end{equation}
for some $R \in  \R_{>0}$. We refer to this $R$ as the \emph{length scale} of the curve and remark that we restrict ourselves to the case
\begin{equation}
    R < 2L,
\end{equation}
where $L$ is the length scale of the flat torus ( with the convention $L = \infty$ in the case of the plane).

To study the relaxation rates under the evolution, we monitor the three quantities \emph{squared distance} $H(t)$,  \emph{energy gap} $E(t)$ and \emph{dissipation} $D(t)$ of the gradient flow, cf. \cite{[Chugreeva2018]}, defined as
\begin{alignat}{2}
    & H(t) = \int_{\mathcal{M}}  \left| \nabla \vartheta(t,\sigma) \right|^2 \,d\sigma &&\quad \left[  R^{4} \right], \label{eq:squareddistancedef} \\
    & E(t) = \lth(\Gamma(t)) - 2 \pi R &&\quad \left[  R \right], \label{eq:energydef} \\
    & D(t) = \int_{\mathcal{M}}  \left| \nabla u(t,\sigma)  \right|^2 \,d\sigma &&\quad \left[  R^{-2} \right], \label{eq:dissipationdef} 
\end{alignat}
where $\vartheta$ solves
\begin{equation}
    - \Delta \vartheta(t,\cdot) = \chi_{\Omega_{\text{in}}(t)}(\cdot) - \chi_{B_{R}(\vek{c}(t))} \quad \text{ on } \mathcal{M}, \label{eq:harmonicdistance}
\end{equation}
with $B_R(\vek{c}(t))$ being the circle of radius $R$ centered at the barycenter $\vek{c}(t)$ of $\Omega_{\text{in}}(t)$. The entries in the square brackets in \eqref{eq:squareddistancedef} through \eqref{eq:dissipationdef} indicate how the quantities scale with $\Gamma$ and $\mathcal{M}$. The dimensionless quantity
\begin{equation}
    \left(E^2D\right)(t)
\end{equation}
will serve as a smallness parameter.

Both $E(t)$ and $H(t)$ measure the distance of $\Gamma(t)$ from the center manifold of compact equilibrium configurations under the flow, i.e., the set of circles in $\mathcal{M}$. While the energy gap measures the excess length of $\Gamma(t)$ compared with any potential limiting circle, the squared distance is the $\dot{H}^{-1}$-norm (of the indicator function)  of the symmetric difference between $\Omega_{\text{in}}(t)$ and a disk of radius $R$ around the barycenter $\vek{c}(t)$. In the field of \emph{geometric stability}, the energy gap is often referred to as the \emph{isoperimetric deficit}, cf. Section \ref{sec:barycenterrole} for a discussion. In a similar spirit, the squared distance may be viewed as an $\dot{H}^{-1}$-analog of the $L^{1}$-based \emph{barycentric asymmetry}, an upper bound for the \emph{Fraenkel asymmetry}, cf. \cite{[Fuglede1993]}. 

In order to state our main result, we make the following definition:
\begin{definition}[Almost circular radial functions]
    For any fixed $R > 0$, and sufficiently small $\delta<1/2$, the set of \emph{(centered) nearly circular radial functions} $\mathcal{R}_{R} := \mathcal{R}_{R,\delta}$ is a subset of $H^{5/2}([0,2 \pi ])$ satisfying the following conditions:
    \begin{align}
        \sup_{\phi \in [0, 2 \pi]} |\rho(\phi) - R| &\leq \delta, \label{eq:admissannulusbound}\\
        \sup_{\phi \in [0, 2 \pi]} \left| \rho_{\phi}(\phi) \right| &\leq \delta, \label{eq:admissslopebound}\\
        \int_{[0, 2 \pi ]} \rho(\phi)^3\cos(\phi) \,d\phi &= 0 = \int_{[0, 2 \pi ]} \rho(\phi)^3\sin(\phi) \,d\phi, \label{eq:admissbarycentric}\\
        \int_{[0,2 \pi]} \rho(\phi)^2 - R^2 \,d\phi &= 0,\label{eq:admissmassconservation}
    \end{align}
    where $\rho_{\phi}$ denotes the derivative of $\rho$ with respect to the polar angle $\phi$.
    Any parametrization of the form \eqref{eq:polarparamcurve} for a $\rho \in \mathcal{R}_{R,\delta}$ is referred to as a \emph{nearly circular (polar) parametrization}.
    \label{def:nearlycircular}
\end{definition}
\begin{remark}
    For sufficiently small $\delta >0$, the first two conditions above render the corresponding curve $\Gamma(t)$ a \emph{nearly spherical surface}, cf.  \cite{[Fuglede1989]}, a class of hypersurfaces enjoying remarkable stability properties in the context of isoperimetric inequalities. Condition \eqref{eq:admissbarycentric} guarantees that the barycenter of the area enclosed by the corresponding curve coincides with the pole of the parametrization, while condition \eqref{eq:admissmassconservation} states that $  \left| \Omega_{\text{in}} \right| = \pi R^2$.
    $\hfill \bigtriangleup$
\end{remark}
We are now ready to state our main result. Below and throughout, we use the $\lesssim$ notation explained in Notation \ref{not:one} below.
\begin{theorem}[Relaxation to a circle]
    Let $\mathcal{M}$ be either $\R^2$ or the flat torus $\mathbb{T}^2_{2L}$ of edge length $L>0$. For a simple closed, separating initial curve $\Gamma(0)$ with 
    \begin{equation}
        \left| \Omega_{\text{in}}(0) \right| = 2 \pi R^2,
    \end{equation}
    assume that $\left( u(t,\cdot), \Gamma(t) \right)_{t\in [0,\infty)}$ is a global solution of the MS system \eqref{eq:msnormal}-\eqref{eq:harmonicproblem}, which is continuously differentiable in time, and let $H, E, D$ represent the corresponding functions of squared distance, energy gap and dissipation. Let further $\left( \rho(t, \cdot) \right)_{t \in [0,\infty)}$ denote the family of radial functions parametrizing the family of interfaces.

    There exist universal constants $C' > 0$ and $\varepsilon > 0$ such that if
    \begin{equation}
        4C' R^2 \leq L^2,
    \end{equation}
    \begin{equation}
        \left(E^2D\right)(0) \leq \varepsilon
    \end{equation}
    and if $\rho(0, \cdot) \in \mathcal{R}_R$, then there holds
    \begin{equation}
        \left( E^2 D \right)(t) \leq \varepsilon 
        \label{eq:stableeed}
    \end{equation}
    and $\rho(t, \cdot) \in \mathcal{R}_R$ for all $t \in [0, \infty)$. The squared distance and the energy gap obey
    \begin{equation}
        H(t) \lesssim \max \left\{ H(0), E(0)R^3 \right\} =: \alpha(0) \label{eq:squareddistancebound}
    \end{equation}
    and
    \begin{equation}
        E(t) \lesssim \min\left\{E(0), \frac{\alpha(0)}{t}, E(0) \exp(-\frac{t}{CR^3}) \right\} \label{eq:decayenergy}
    \end{equation}
    respectively with some universal constant $C \in \R_{>0}$.  
    Further, there exists a crossover time $ T_{1} \lesssim R^3$ separating two relaxation regimes for the energy gap: Initially, i.e., for times $t \in [0, T_{1}]$, the energy gap decays algebraically as
    \begin{equation}
        E(t) \lesssim \min\left\{E(0), \frac{\alpha(0)}{t}\right\}, \label{eq:algebraicenergy}
    \end{equation}
    and for $t \in [T_{1}, \infty)$ it decays exponentially as
    \begin{equation}
        E(t) \lesssim E(0) \exp(-\frac{t}{CR^3}). \label{eq:exponentialenergy}
    \end{equation}
    Furthermore, there exists a time $T_D \lesssim \alpha_0^{3/4}$ such that the rate of decay of the dissipation is algebraic in time for $t \in [T_D, T_{1}]$,
    \begin{equation}
        D(t) \lesssim \frac{\alpha_{0}}{t^2}, \label{eq:algebraicdissipation}
    \end{equation}
    and then exponentially for $t \in [T_{1}, \infty)$:
    \begin{equation}
        D(t) \lesssim \frac{E(0)}{t} \exp(- \frac{t}{CR^3}). \label{eq:exponentialdissipation}
    \end{equation}

    Moreover, the barycenter of the curve $c(t)$ is confined to a ball of radius $\left( E(0) R \right)^{1/2}$ at the origin, i.e.,
    \begin{equation}
          \left| \vek{c}(t) \right| \lesssim \left( E(0) R \right)^{1/2}
          \label{eq:barycenterbound}
    \end{equation}
    for all times $t \in [0,\infty)$.
    \label{thm:main}
\end{theorem}

\subsubsection*{Previous results and main contribution}

Exponential convergence to a sphere for well-prepared initial data was established in \cite{[Chen1993]} and \cite{[EscherSimonett1998]} under the assumption of either a small dissipation or $C^{2+\beta}$-closeness of the initial manifold. Algebraic convergence rates for a non-compact interface with graph structure was proved in \cite{[Chugreeva2018]}; see also \cite{[Otto2025]}. 

The recent article \cite{[Arya2025]}, building on the work \cite{[Julin2022]}, established exponential convergence of flat flows in the flat 2-torus starting from sets of finite perimeter to collections of circles of equal area or disjoint parallel strips. While the focus in [AGK25] is on the existence of global attractors for quite general initial data, our focus is instead on identifying natural conditions on the initial data and, for this class of initial data, quantifying the optimal convergence rates up to universal prefactors.

Our assumption on the initial data is that it has graph structure and that the nondimensional quantity $E(0)^2D(0)$ is sufficiently small. These assumptions allow, for instance, for initial configurations with large dissipation and large $C^2$-distance to the center manifold, provided the energy gap (isoperimetric deficit), and thus the $C^{1}$-distance, is sufficiently small. The perturbation of a circle with low amplitude and high frequency sine oscillations is an example.

The main novelty of Theorem \ref{thm:main} is the resolution of two successive relaxation regimes in the MS evolution: an initial regime, in which relaxation is dominated by an algebraic rate, followed by a regime of exponential decay. This change of regimes takes place at a time $T_{1} \lesssim R^3$. The resolution of the algebraic decay regime is based on the gradient flow structure of the evolution, while the existence of the exponential decay regime follows from a useful energy-energy-dissipation (EED) relationship, cf. \eqref{eq:energydissipationrelation}. The fact that such an EED inequality exists in the compact case reflects the existence of a spectral gap for the linearization of the evolution operator at elements of the center manifold.

In addition to relaxation rates, we provide a bound for the maximal motion of the barycenter of the interface in terms of the initial energy gap. These bounds resemble estimates for the difference between the outer and inner radius of the minimal annulus containing the curve from the Bonnesen inequality, cf. \eqref{eq:bonnesen}. The proof of this bound rests on the explicit resolution of the algebraic and the exponential relaxation regime.

On the technical side, a novelty of this work is the application of periodic single-layer potential theory via elliptic functions for the Mullins-Sekerka flow; see Lemma \ref{lemma:improvedsobolevnormalvelo}.

\subsection{Method}
The two different modes of relaxation reflect two different mechanisms shaping the basin of attraction for the evolution: The first one, the mode of algebraic decay, is induced by the mildly non-convex gradient flow structure of the evolution: In \cite{[Brezis1973]}, it was shown that a gradient flow in a metric space with respect to a convex energy decays algebraically. More precisely, convexity of some energy $\mathcal{E}$ implies
\begin{equation}
    \frac{d \mathcal{E}}{d t} \leq - \mathcal{D}, \quad \frac{d \mathcal{D}}{d t} \leq 0, \quad \frac{d \mathcal{H}}{d t} \leq 0, \quad \mathcal{E} \leq \left( \mathcal{H} \mathcal{D} \right)^{1/2},
    \label{eq:brezisflow}
\end{equation}
which immediately implies
\begin{equation}
    \mathcal{E}(t) \leq \min\left\{ \mathcal{E}(0), \frac{\mathcal{H}(0)}{t} \right\}, \quad \mathcal{H}(t) \leq \mathcal{H}(0), \quad \mathcal{D}(t) \leq \frac{2 \mathcal{H}(0)}{t^2}, \quad \forall t > 0,
\end{equation}
where $\mathcal{D}(t)$ is the dissipation and $\mathcal{H}(t)$ the squared distance of the current state to the center manifold of equilibria (measured with respect to the metric of the space). In \cite{[Otto2014]}, the authors applied this relaxation framework to the one-dimensional Cahn-Hilliard equation -- a gradient flow with respect to a non-convex energy. Thus, the authors showed that the framework can be extended to non-convex settings, by showing that the estimates in \eqref{eq:brezisflow} hold up to (benign) error terms; see the corresponding discussion in \cite[p. 734]{[Otto2014]}. In \cite{[Chugreeva2018]}, this framework was applied to the Mullins-Sekerka evolution of a graph over the entire $x$-axis in the plane.

The second mode of decay, the exponential one, is a consequence of the existence of a spectral gap for the linearization of the evolution operator around the center manifold: If the linearization of the evolution operator of a (continuous) dynamical system around its center manifold possesses a spectral gap, the convergence of the flow along the stable manifold is exponentially fast. In \cite{[EscherSimonett1998]}, for instance, the exponentially fast convergence of suitable initial data to spheres under the Mullins-Sekerka evolution was established by carefully analyzing the spectral gap of the linearized evolution operator. In the article at hand, however, we do not make explicit use of this spectral gap. Instead, we exploit the useful EED-inequality
\begin{equation}
    E \lesssim R^3 D,
\end{equation}
which immediately yields the exponential decay of $E$ upon combination with $dE/dt = -D$. The above inequality is a consequence of the fact that the energy gap $E$ can be controlled in terms of $\left\| \kappa - \overline{\kappa} \right\|^2_{L^2(\Gamma)}$, with $\overline{\kappa}$ being the average curvature, provided it is sufficiently small. The fact that the proof of this geometric stability result rests upon the spectral properties of the Laplacian, is another indication of the importance of the presence of a spectral gap for the exponential relaxation regime.

In order to resolve both relaxation regimes, we thus need to correctly capture both mechanisms at play: To capture the algebraic decay, we derive the necessary algebraic and differential relationships among $H, E$ and $D$, cf. Proposition \ref{prop:interpolation}, Lemma \ref{lemma:geometricvariations} and Proposition \ref{prop:diffrelations}, leading to a structure similar to that of \eqref{eq:brezisflow}. To capture the exponential decay, we derive the EED-inequality, cf. Corollary \ref{cor:energydissipation}. Because many of these estimates rely on the smallness of the nondimensional quantity $\left(E^2D\right)(t)$ and the explicit graph structure of the curve, we establish that this quantity is non-increasing in time. This indeed suffices to maintain both conditions throughout the evolution, provided they are satisfied initially, cf. Corollary \ref{cor:stableeed} and Remark \ref{rmk:graph}.

\subsection{Organization and Notation} \label{sec:organization}

As pointed out above, the proof of Theorem \ref{thm:main} is based on a combination of algebraic and differential estimates relating $H$, $E$ and $D$ with elementary ODE-type arguments. We begin by stating some general geometric facts and observations in Section \ref{sec:geometry}, before proving first the algebraic relations in $H$, $E$ and $D$ in Section \ref{sec:algebraic}, followed by the differential ones in Section \ref{sec:differential}. The proof of Theorem \ref{thm:main} is then given in Section \ref{sec:main}. Some technical facts about fractional Sobolev spaces, elementary differential geometry and elliptic functions are collected in the appendix.

\begin{notation} \label{not:one}
    The notation
    \begin{equation}
        A \lesssim B
    \end{equation}
    expresses that for two quantities $A$ and $B$ there exists a universal constant $C \in \R_{>0}$ such that
    \begin{equation}
        A \leq C B
    \end{equation}
    holds. Symmetry in this condition, i.e.,
    \begin{equation}
        A \lesssim B \quad \text{ and } \quad A \gtrsim B
    \end{equation}
    is expressed via
    \begin{equation}
        A \sim B.
    \end{equation}

    Throughout this article, we omit time dependence whenever it is clear from the context. For the sake of readability, (partial) derivatives may be denoted in the form $f_{x}$.
    To avoid confusion with partial derivatives, subscripts indicating components are put into parenthesis. The $x$-coordinate of the center function $\vek{c}(t)$, for example, is expressed as $c_{(x)}(t)$. Time derivatives may be denoted by a dot, e.g., $\dot{c}_{(x)}(t)$.

    Throughout the text, we identify domains of the kind $[0, 2 \pi]$ or $[0, \lth(\Gamma)]$ with the one-dimensional tori of equal measure. In particular, any functions defined on these domains are periodic.

    Moreover, we write $ \left| A \right|$ to denote the Lebesgue measure of a subset $A \subset \mathcal{M}$. For vectors $\vek{b}, \vek{c} \in \R^{n}$, $ \left| \vek{b}\right|$ is used to denote the Euclidean norm of $\vek{b}$, while $\vek{b} \cdot \vek{c}$ denotes the Euclidean inner product.
    $\hfill \bigtriangleup$
\end{notation}

    \section{Elementary geometric facts}
\label{sec:geometry}

\begin{figure}[ht]
    \centering
    \includegraphics[width=0.5\textwidth]{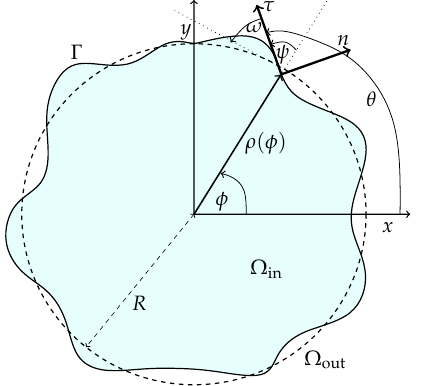}
    \caption{Polar parametrization of $\Gamma$ and definition of the relevant angles}
    \label{fig:angles}
\end{figure}

For the remainder of this section, we fix some arbitrary point in time $t \in [0,\infty)$ with corresponding translation $c(t)$. Figure \ref{fig:angles} depicts the geometric setup, where the curve $\Gamma(t) \equiv \Gamma$ is parametrized via some radial function $\rho$ according to \eqref{eq:polarparamcurve}. 

Before we collect some elementary geometric relations that will be used throughout this article, we first define and motivate the angles $\omega$, $\psi$ and $\theta$ depicted in Figure \ref{fig:angles}:

For every $\phi \in [0, 2 \pi]$ the \emph{inclination angle} $\theta(\phi)$ is defined as the (positively oriented) angle between the $x$-axis and the tangent vector $\tau(\phi)$. The inclination angle is of particular relevance due to the well-known relation
\begin{equation}
    \kappa(s) =  \left. \frac{d \theta}{d s} \right|_{s} \quad \text{ for all } s \in [0, \lth(\Gamma)],
    \label{eq:inclinationcurvature}
\end{equation}
where $s$ denotes the arc-length parameter of the corresponding arc-length parametrization of the curve. Curvature can thus be expressed as an exact form (w.r.t arc-length parametrization) with non-periodic primitive $\theta$.

To obtain a periodic quantity, we subtract the polar angle $\phi$ from $\theta$ to obtain the \emph{polar tangential angle} $\psi(\phi)$
\begin{equation}
    \psi(\phi) = \theta(\phi) - \phi \quad \text{ for all } \phi \in [0, 2 \pi],
\end{equation}
which is the (positively oriented) angle between the polar ray through $(\cos(\phi), \sin(\phi))$ and the tangent vector $\tau(\phi)$. 

A similar behavior can be observed for the positively oriented angle between the tangent vector $\tau(\phi)$ of $\Gamma$ and the corresponding tangent vector of the circle of radius $R$ centered at the barycenter of $\Omega_{\text{in}}$. This angle, subsequently referred to as the \emph{localized angle} $\omega(\phi)$, satisfies
\begin{equation}
    \omega(\phi) = \frac{\pi}{2} - \psi(\phi) = \frac{\pi}{2} + \phi - \theta(\phi).
\end{equation}
Notice, however, that $\omega$ and $\psi$ are not primitives of $\kappa$ with respect to $\phi$. Instead, curvature decomposes into two components: the intrinsic curvature of the reference circle and the curvature induced by the bending of the curve with respect to the circle. This is reflected in \eqref{eq:kappapolar} of the following collection of elementary geometric facts:
\begin{lemma}[Elementary Geometric Relations]
    Let $\rho \in \mathcal{R}_R$, let $\gamma_{\text{pol}}$ denote the corresponding polar parametrization and let $\phi \in [0, 2 \pi]$ be arbitrary. Then the induced curve $\Gamma$ satisfies the following relations:

    The tangent vector $\tau(\phi)$ is computed as
    \begin{equation}
        \tau(\phi) = \begin{pmatrix} \rho_{\phi}(\phi)(\cos(\phi)) - \rho(\phi) \sin(\phi) \\ \rho_{\phi}(\phi)(\sin(\phi)) + \rho(\phi) \cos(\phi)  \end{pmatrix},
        \label{eq:tangentvec}
    \end{equation}
    while the outward pointing unit normal vector $n(\phi)$ is given by
    \begin{equation}
        n(\phi) = \frac{1}{\sqrt{\rho(\phi)^2 + \rho_{\phi}(\phi)^2}} \begin{pmatrix} \rho_{\phi}(\phi)(\sin(\phi)) + \rho(\phi) \cos(\phi)  \\ - \rho_{\phi}(\phi)(\cos(\phi)) + \rho(\phi) \sin(\phi) \end{pmatrix},
        \label{eq:normalvec}
    \end{equation}
    and the corresponding \emph{length element} $\llm$ of the curve, defined as
    \begin{equation}
        \llm: [0, 2 \pi] \to \R_{>0}; \, \phi \mapsto \left\| \frac{\partial \gamma}{\partial \phi}  \right\|(\phi) =  \left. \frac{\partial s}{\partial \phi}  \right|_{\phi},
    \end{equation}
    where $s$ is the arc length parameter, satisfies
    \begin{equation}
        \llm(\phi) = \sqrt{\rho(\phi)^2 + \rho_{\phi}(\phi)^2}.
        \label{eq:lengthelement}
    \end{equation}
    Further, there holds
    \begin{equation}
        \lth(\Gamma) = \int_{\Gamma} ds = \int_{[0, 2 \pi]} \llm(\phi)d\phi \sim R,
        \label{eq:curvelength}
    \end{equation}
    where $\lth(\Gamma)$ denotes the length of the curve.

    The normal velocity $V$ is related to the parametrization via
    \begin{equation}
        \left( V(t,\phi) - \left(\dot{\vek{c}}(t) \cdot n(\phi)\right)  \right) \llm(t,\phi)= \rho_t(t,\phi) \rho(t,\phi),
        \label{eq:normalveloparam}
    \end{equation}
    and the curvature of $\Gamma$ is given by
    \begin{equation}
        \kappa(\phi) = \left( -  \left. \frac{\partial \omega}{\partial \phi} \right|_{\phi}  + 1\right) \frac{1}{\llm(\phi)} = \frac{\rho_{\phi}(\phi)^2 - \rho(\phi) \rho_{\phi \phi}(\phi)}{\llm(\phi)^{3}} + \frac{1}{\llm(\phi)},
        \label{eq:kappapolar}
    \end{equation}
    i.e,
    \begin{equation}
        \left. \frac{\partial \omega}{\partial \phi} \right|_{\phi} = \frac{\rho(\phi) \rho_{\phi \phi}(\phi) - \rho_{\phi}(\phi)^2 }{\llm(\phi)^2}.
        \label{eq:localizedanglederivative}
    \end{equation}

    The localized angle furthermore satisfies
    \begin{equation}
        \tan(\omega)(\phi) = \frac{\rho_{\phi}(\phi)}{\rho(\phi)} = \frac{\partial}{\partial \phi}  \ln\left( \frac{\rho}{R}. \right) \label{eq:tanomega}
    \end{equation}
    and 
    \begin{equation}
        \tan(\omega)(\phi) \sim \omega(\phi).
        \label{eq:tanomegaisomega}
    \end{equation}
    
    \label{lemma:geometricquantitiesrho}
\end{lemma}
\begin{myproof}
    Equations \eqref{eq:tangentvec}, \eqref{eq:normalvec} and \eqref{eq:lengthelement} follow from elementary differential geometric computations, and the bounds on $\rho$ and $\rho_{\phi}$, cf. \eqref{eq:admissslopebound} and \eqref{eq:admissannulusbound}, imply  \eqref{eq:curvelength}. To obtain \eqref{eq:normalveloparam}, one expands
    \begin{equation}
        V(t,\phi) = \left. \frac{\partial \gamma_{\text{pol}}}{\partial t} \right|_{(t,\phi)} \cdot n(t,\phi).
    \end{equation}
    
    Reexpressing \eqref{eq:inclinationcurvature} in polar coordinates yields \eqref{eq:kappapolar} and \eqref{eq:localizedanglederivative}, while \eqref{eq:tanomega} follows from basic trigonometry and leads to  \eqref{eq:tanomegaisomega} through a Taylor approximation. 
\end{myproof}

\subsection{Bulk parametrization}

To account for the bulk diffusion mechanism driving the MS evolution, we introduce two kinds of parametrizations of the phases $\Omega_{\text{in}}$ and $\Omega_{\text{out}}$ and the ambient manifold: The first is a polar parametrization, where the geometry of the phases is reflected in the parametrization domain. In contrast, the second kind, subsequently referred to as the flattening parametrization, encodes the geometry in the parametrization function. While the former will turn out to be useful for explicit integral computations later, the latter will help to approximate some bulk quantities in the physical phases by corresponding quantities on disks and their complements.

Introducing 
\begin{equation}
    P: (0, \infty) \times [0, 2 \pi];\, (r, \phi) \mapsto  r \begin{pmatrix} \cos(\phi) \\ \sin(\phi)  \end{pmatrix},
    \label{eq:polarparam}
\end{equation}
we can parametrize $\Omega_{\text{in}}(t)$ on both ambient manifolds by restricting and translating via
\begin{equation}
    P_{\text{in}}(t, \cdot, \cdot ) := \vek{c}(t) + \left. P \right|_{\mathcal{P}_{\text{in}}(t)} 
\end{equation}
for
\begin{equation}
    \mathcal{P}_{\text{in}}(t) := \bigsqcup_{\phi \in [0,2 \pi]} [0, \rho(t, \phi)] \times \{\phi\}.
\end{equation}
Analogously, 
\begin{equation}
     P_{\text{out},\R^2}(t, \cdot , \cdot ) := \vek{c}(t) +  \left. P \right|_{\mathcal{P}_{\text{out},\R^2}(t)} 
\end{equation}
and
\begin{equation}
    P_{\text{out},\mathbb{T}^2_{2L}}(t, \cdot , \cdot ) := \vek{c}(t) +  \left. P \right|_{\mathcal{P}_{\text{out},\mathbb{T}^2_{2L}}(t)}
\end{equation}
for
\begin{equation}
    \mathcal{P}(t)_{\text{out}, \R^2} := \bigsqcup_{\phi \in [0,2 \pi]} [\rho(t,\phi), \infty) \times \{\phi\},
\end{equation}
and
\begin{equation}
    \mathcal{P}(t)_{\text{out}, \mathbb{T}^2_{2L}} := \bigsqcup_{\phi \in [0,2 \pi]} [\rho(t,\phi), R(\phi)_{\text{max}, 2L}) \times \{\phi\},
\end{equation}
parametrize $\Omega_{\text{out}}(t)$ on $\R^2$ and $\mathbb{T}^2_{2L}$ respectively. Here, $R(\cdot)_{\text{max}, 2L}$ is the continuous function such that the image of
\begin{equation}
    \phi \mapsto R(\phi)_{\text{max}, 2L} \begin{pmatrix} \cos(\phi) \\ \sin(\phi)  \end{pmatrix}
\end{equation}
is the centered square of edge length $2L$.

In contrast, the flattening parametrization for $\Omega_{\text{in}}(t)$ is given as
\begin{equation}
    F_{\text{in}}(t, \cdot, \cdot): \mathcal{F}_{\text{in}} := (0, R) \times [0, 2 \pi] \to \Omega_{\text{in}}; \, (r, \phi) \mapsto \vek{c}(t) + r \frac{\rho(\phi)}{R} \begin{pmatrix} \cos(\phi) \\ \sin(\phi)  \end{pmatrix},
    \label{eq:interiorflat}
\end{equation}
while the corresponding flattening parametrizations for $\Omega_{\text{out},\R^2}(t)$ and $ \Omega_{\text{out},\mathbb{T}^2_{2L}}(t)$ are chosen as
\begin{equation}
    F_{\text{out}, \R^2}(t, \cdot, \cdot ): \mathcal{F}_{\text{out},\R^2} := (R, \infty) \times [0, 2 \pi] \to \Omega_{\text{out}(t),\R^2}; \, (r, \phi) \mapsto \vek{c}(t) + r \frac{\rho(\phi)}{R} \begin{pmatrix} \cos(\phi) \\ \sin(\phi)  \end{pmatrix}
    \label{eq:planeflat}
\end{equation}
and
\begin{align}
    F_{\text{out}, \mathbb{T}^2_{2L}}(t, \cdot, \cdot): &\mathcal{F}_{\text{out}, \mathbb{T}^2_{2L}} \to \Omega_{\text{out}}(t),\\
    &(r,\phi) \mapsto \vek{c}(t) + r\left( 1 + \left( \frac{\rho(\phi)}{R} - 1 \right) \beta(r) \right) \begin{pmatrix} \cos(\phi) \\ \sin(\phi)  \end{pmatrix}
    \label{eq:twopolarparambulkout},
\end{align}
with 
\begin{equation}
    \mathcal{F}_{\text{out}, \mathbb{T}^2_{2L}} := \bigsqcup_{\phi \in [0,2 \pi]} (R, R(\phi)_{\text{max}, 2L}) \times \{\phi\}.
    \label{eq:torusdomain}
\end{equation}

Here $\beta \in \C^{\infty}_c(\R_{\geq 0})$ is a radially non-increasing bump function satisfying
\begin{alignat}{2}
    &\beta(r) = 1  \quad && \text{if }   \left| r \right| \leq R, \\
    0 \leq &\beta(r) < 1 \quad && \text{if }  R <  \left| r \right| < \frac{3R}{2}, \label{eq:betaone} \\
    &\beta(r) = 0  \quad && \text{if }   \frac{3R}{2} < \left| r \right|,  \\
\end{alignat}
and
\begin{equation}
     \left| \partial_r \beta(r) \right| \leq \frac{3}{R}. \label{eq:betatwo}
\end{equation}
Some differential geometric quantities related to these parametrizations are collected in the appendix, Section \ref{app:bulkparams}.
Notice that inserting $\beta \equiv 1$ into \eqref{eq:twopolarparambulkout} yields the formulas for $F_{\text{in}}$ and $F_{\text{out}}$.

\begin{remark}
    Parametrizations of the form \eqref{eq:twopolarparambulkout}, which transform bulk problems on evolving perturbations of some model domain to problems on that model domain, are also referred to as \emph{Hanzawa transformations}, cf. \cite{[Hanzawa1981]}. The Hanzawa transformation, combined with the nearly-circular nature of the considered interfaces, allows for approximations of bulk terms on $\Omega_{\text{in}}$ and $\Omega_{\text{out}}$ in terms of quantities on a circle and its complement.
    $\hfill \bigtriangleup$
\end{remark}

\begin{remark}
    From now on, we will, unless stated otherwise, use $\mathbb{T}^2_{2L}$ as the ambient manifold and thus generally drop the corresponding subscripts, as all the techniques used here equally apply to (or even simplify in) the case of the plane.
    $\hfill \bigtriangleup$
\end{remark}

\subsection{Role of the barycenter}
\label{sec:barycenterrole}
Before we study algebraic relations among $H$, $E$ and $D$ and their relations with (Sobolev norms of) basic geometric quantities in the next section, a brief remark is due on the role of the barycenter condition \eqref{eq:admissbarycentric}. 

Generally, the barycenter $\vek{c} \in \mathcal{M}$ of a simply connected open set $\Omega \subset \mathcal{M}$ whose boundary is a simple closed curve $\partial \Omega$ is computed via
\begin{equation}
    \vek{c} = \begin{pmatrix} c_{(x)} \\ c_{(y)} \end{pmatrix} = \frac{1}{ \left| \Omega \right|} \begin{pmatrix} \displaystyle \int_{\Omega} x \,dV \\[2ex] \displaystyle \int_{\Omega} y \,dV\end{pmatrix}.
    \label{eq:bulkbarycenter}
\end{equation}
If $\partial \Omega$ admits a polar parametrization of the form
\begin{equation}
    \gamma(\phi) = \rho(\phi) \begin{pmatrix} \cos(\phi) \\ \sin(\phi) \end{pmatrix},
\end{equation}
this simplifies to
\begin{equation}
    \vek{c} = \frac{1}{ 3 \left| \Omega  \right|} \begin{pmatrix} \displaystyle \int_{[0, 2 \pi]} \rho(\phi)^3 \cos(\phi) \,d\phi \\[2ex] \displaystyle \int_{[0, 2 \pi]} \rho(\phi)^3 \sin(\phi) \,d\phi \end{pmatrix},
\end{equation}
through explicit integration over the radial component.

In \cite{[Fuglede1989]}, it is shown that nearly spherical $d$-dimensional surfaces, i.e., surfaces $C^{1}$-close to the unit sphere with a centered barycenter, admit strong isoperimetric stability properties, cf. Lemma \ref{lm:nearlyspherical} below. This means that the isoperimetric deficit $\Delta(S)$ of the surface $S$ satisfies
\begin{equation}
    \left\| \rho - R \right\|^2_{L^2(S^{d})} + \left\| \rho_{\phi} \right\|^2_{L^2(S^{d})} \lesssim \Delta(S) \lesssim \left\| \rho_{\phi} \right\|^2_{L^2(S^{d})},
\end{equation}
where $\rho$ is the radial function in a polar parametrization of the surface over the unit sphere. For a surface $S$ bounding an open set $\Omega$ of finite Lebesgue measure $ \left| \Omega  \right|$, the isoperimetric deficit is defined as
\begin{equation}
    \Delta(S) :=  \left| S \right| - d  \left| B_{1}(0) \right|^{1/d} \left| \Omega \right|^{\frac{d-1}{d}},
\end{equation}
where $\ \left| S \right| = \text{Per}(S)$ denotes the perimeter of $S$, $ \left| B_{1}(0) \right|$ is the volume of the $d$-dimensional unit ball and $ \left|  \Omega \right|$ denotes the volume enclosed by $S$. In the case of a two-dimensional disk of radius $R$, i.e., $\Omega = B_{R}(0)$, there holds
\begin{equation}
    d  \left| \Omega \right|^{\frac{d-1}{d}} = 2 \left( \pi \right)^{1/2} \left( \pi R^2 \right)^{\frac{2-1}{2}} = 2 \pi R= \text{Per}(B_{R}(0)).
\end{equation}
Thus, our energy gap, cf. \eqref{eq:energydef}, equals the isoperimetric deficit.
Paired with the isoperimetric inequality
\begin{equation}
    \Delta \geq 0,
\end{equation}
cf. \cite[(1)]{[Fuglede1989]}, this implies that smallness of the norm $\left\| \rho_{\phi} \right\|^2_{L^2(S^{d})}$ is equivalent to closeness to a sphere, under the assumption  of the barycenter being centered.

Effectively, fixing the barycenter at $0$ eliminates any potential perturbations of $\rho$ stemming from pure translations.

\begin{figure}[ht]
    \centering
    \includegraphics[width=0.35\textwidth]{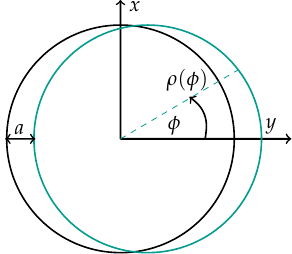}
    \caption{Sketch illustrating the unit circle shifted to the right by $a$.}
    \label{fig:shiftedcircle}
\end{figure}

Thus, it is not surprising that the most basic counterexample to control of $\left\| \rho_{\phi} \right\|_{L^2([0,2pi])}$ in terms of the deficit is constructed by shifting a unit circle by some $a \in (0,1)$, cf. Figure \ref{fig:shiftedcircle}. The radial function in this case is given by
\begin{equation}
    \rho(\phi) = a \cos(\phi) + \sqrt{1 - a^2 \sin^2(\phi)}, 
\end{equation}
for which certainly
\begin{equation}
    \left\| \rho_{\phi} \right\|^2_{L^2([0, 2 \pi])} > 0
\end{equation}
holds, while the translation invariant isoperimetric deficit vanishes. Hence, only radial parametrizations with a centered barycenter allow to mutually infer information about the $\dot{H}^{1}$-norm of the radial function and the energy gap $E$ from one another. After briefly examining an alternative notion of barycenter, we turn to relations of this kind in the next section.

\subsubsection{Different notions of barycenters}
Instead of considering the barycenter of some open set $\Omega \subset \mathcal{M}$ whose boundary is a simple closed curve $\partial \Omega$ as in \eqref{eq:bulkbarycenter}, one can immediately compute the barycenter of the boundary curve $\partial \Omega$ as
\begin{equation}
    \vek{c}_{\partial \Omega} := \frac{1}{\lth(\partial \Omega)} \begin{pmatrix} \displaystyle \int_{\partial \Omega} x \,ds \\[2ex] \displaystyle \int_{\partial \Omega} y \,ds \end{pmatrix}.
\end{equation}
While the former notion is, e.g., being used in \cite{[Fuglede1989]} in the context of isoperimetric stability, in \cite{[Gancedo2023]} in the context of the Muskat flow and in the present article, the latter notion is used in \cite{[KrummelMaggi2016]} in the context of isoperimetric stability. For star-shaped domains, these quantities are equal up to a constant factor, so that $c = 0$ holds if and only if $c_{\partial \Omega} = 0$:

\begin{lemma}[Equivalence of barycenters]
    Let $\Omega \subset \mathcal{M}$ be an open, star-shaped domain with boundary $\partial \Omega$ a simple closed curve. Then the barycenters $\vek{c}$ and $\vek{c}_{\partial \Omega}$ satisfy 
    \begin{equation}
        \vek{c} = \frac{\lth(\partial \Omega)}{ 3 \left| \Omega \right|} \vek{c}_{\partial \Omega}.
    \end{equation}
    \label{lm:barycen}
\end{lemma}
\begin{myproof}
    The main idea is to ``foliate'' $\Omega$ by scaled copies of $\partial \Omega$. To this end, we first remark that scaling the boundary curve scales $\vek{c}_{\partial \Omega}$ accordingly, i.e.,
    \begin{equation}
        \vek{c}_{t \partial \Omega} = t \vek{c}_{\partial \Omega},
    \end{equation}
    for all $t > 0$, where $t \partial \Omega := \{ tx \in \mathcal{M} \mid x \in \partial \Omega \}$. 
    Using the coarea formula, cf. \cite[Chapter 13]{[Maggi2012]}, we can now reexpress $c_{(x)}$ as
    \begin{align}
        c_{(x)} &= \frac{1}{ \left| \Omega  \right|} \int_{\Omega} x\, dV = \frac{1}{ \left| \Omega \right|} \int_{0}^{1} \frac{\lth(t \partial \Omega)}{\lth(t \partial \Omega)} \left(\int_{t \partial \Omega} x\, ds \right) dt 
                = \frac{\lth(\partial \Omega)}{ \left| \Omega \right|} \int_{0}^{1} t^2 \left(c_{\partial \Omega}\right)_{(x)} dt =  \frac{\lth(\partial \Omega)}{ 3 \left| \Omega \right|} \left(c_{\partial \Omega}\right)_{(x)},
    \end{align}
    which establishes the claim.
\end{myproof}

    \section{Algebraic relations} \label{sec:algebraic}
The major goal of this section is to prove the two relations
\begin{align}
    E &\lesssim R^3 D \label{eq:energydissspoiler} \\
    E &\lesssim \sqrt{H D} \label{eq:energyinterpolatonspoiler}
\end{align}
among $H$, $E$ and $D$ as well as the nondimensional relation
\begin{equation}
    \left\| \frac{\rho_{\phi}}{R} \right\|_{L^{\infty}([0, 2 \pi])} \lesssim (E^2 D)^{1/6}.
    \label{eq:slopeboundspoiler}
\end{equation}
In this section, we omit the $t$-dependence, because the estimates are static.

\subsection{An energy-dissipation relation}
The following theorem, cf. \cite[Theorem 1.10]{[KrummelMaggi2016]}, is a main ingredient in the proof of \eqref{eq:energydissspoiler}:
\begin{theorem}[Boundary control in terms of curvature oscillations]
    For any $n \geq 1$ there exists an $\varepsilon(n) > 0$ such that for any open set $\Omega \in \R^{n}$ with $C^{1,1}$-boundary $\partial \Omega := \{ (1 + u(x)) x: x \in S^{n} \}$ satisfying
    \begin{equation}
        \int_{\partial \Omega} x dS = 0
        \label{eq:barycondmaggi}
    \end{equation}
    and for any function $ u \in  C^{1,1}(S^{n})$ satisfying
    \begin{equation}
        \left\| u \right\|_{C^{1}(S^{n})} \leq \varepsilon(n),
        \label{eq:smallgraphcond}
    \end{equation}
    there holds 
    \begin{equation}
        \left\| u \right\|_{W^{1,2}(S^{n})} \leq C(n) \left\| H - n \right\|_{L^2(\partial \Omega)},
        \label{eq:maggistability}
    \end{equation}
    where $H$ denotes the mean curvature of $\partial \Omega$.
    \label{thm:maggiresult}
\end{theorem}
\begin{remark}
    By Lemma \ref{lm:barycen}, the above theorem applies to any nearly spherical curve, since $c = c_{\partial \Omega} = 0$. 
    $\hfill \bigtriangleup$
\end{remark}
The following corollary adapts the above result to our needs:
\begin{corollary}
    There exists $\delta \in (0,\frac{1}{2})$ such that for $\rho \in \mathcal{R}_R$ satisfying
    \begin{equation}
        \left\| \rho - R \right\|_{C^{1}([[0, 2 \pi]])} \leq \delta R,
    \end{equation}
    there holds
    \begin{equation}
        \left\| \rho_{\phi} \right\|^{2}_{L^2([0, 2 \pi])} \lesssim R^{4} \left\| \kappa - \overline{\kappa} \right\|^2_{L^2([0, 2 \pi])} \sim R^3 \left\| \kappa - \overline{\kappa} \right\|^2_{L^2(\Gamma)},
        \label{eq:alexandrov}
    \end{equation}
    where $\kappa$ denotes the curvature of the curve and $\overline{\kappa}$ denotes its integral mean with respect to arc-length parametrization.
\end{corollary}
\begin{myproof}
    We claim that it suffices to prove the result for $R=1$: Suppose $\rho \in \mathcal{R}_{R}$ parametrizes a nearly circular curve $\Gamma$ for some  arbitrary $R > 0$. Then, the rescaled curve $\Gamma^{\ast} := (1/R) \Gamma$ with rescaled radial function $\rho^{\ast} := (1/R) \rho \in \mathcal{R}_{1}$ is a nearly circular curve of length scale $R^{\ast}=1$. Assuming that the claim holds for the case $R=1$, this implies 
    \begin{equation}
        \left\| \rho_{\phi}^{\ast} \right\|^{2}_{L^2([0, 2 \pi])} \lesssim \left\| \kappa^{\ast} - \overline{\kappa^{\ast}} \right\|^2_{L^2([0, 2 \pi])} \sim \left\| \kappa^{\ast} - \overline{\kappa^{\ast}} \right\|^2_{L^2(\Gamma^{\ast})},
        \label{eq:alexandrovone}
    \end{equation}
    where $\kappa^{\ast}$ denotes the curvature of the rescaled curve $\Gamma^{\ast}$. The scaling relations
    \begin{align}
        \left\| \rho_{\phi}^{\ast} \right\|^{2}_{L^2([0, 2 \pi])} &= \frac{1}{R^2} \left\| \rho_{\phi} \right\|^{2}_{L^2([0, 2 \pi])}, \\
        \left\| \kappa^{\ast} - \overline{\kappa^{\ast}} \right\|^2_{L^2([0, 2 \pi])} &= R^2 \left\| \kappa - \overline{\kappa} \right\|^2_{L^2([0, 2 \pi])}, \\
        \left\| \kappa^{\ast} - \overline{\kappa^{\ast}} \right\|^2_{L^2(\Gamma^{\ast})} &= R \left\| \kappa - \overline{\kappa} \right\|^2_{L^2(\Gamma)}
    \end{align}
    combined with \eqref{eq:alexandrovone} then imply the claim for arbitrary $R > 0$.

    Hence, we can subsequently assume $R=1$ without loss of generality. Choosing $\delta$ small enough implies that 
    \begin{equation}
        \left\| \rho_{\phi} \right\|^2_{L^2([0, 2 \pi])} \leq \left\| \rho - R \right\|^2_{L^2([0, 2 \pi])} \lesssim \left\| \kappa - 1 \right\|^2_{L^2(\Gamma)}
    \end{equation}
    holds by Theorem \ref{thm:maggiresult}. Thus, the task is to replace $1$ by $\overline{\kappa}$ in the above estimate. From the Gauss-Bonnet theorem, it is known that
    \begin{equation}
        \overline{\kappa} = \frac{2 \pi }{\lth(\Gamma)}
    \end{equation}
    holds. Define a scaling factor $\lambda := \frac{2 \pi }{\lth(\Gamma)}$ and introduce the rescaled radial function
    \begin{equation}
        \rho' := \lambda \rho,
    \end{equation}
    whose curvature and corresponding curve will be denoted by $\kappa'$ and $\Gamma'$ respectively. As the isoperimetric inequality implies $\lth(\Gamma) \geq 2 \pi$, and homothetic scalings of curves merely scale the barycenter, it follows that $\rho' \in \mathcal{R}_{\lambda}$ and $\left\| \rho'\right\|_{C^{1}(2 \pi )} \leq \delta$ hold, so that Theorem \ref{thm:maggiresult} applies again. Now, the rescaling and the Gauss-Bonnet theorem imply the estimate
    \begin{equation}
        \left\| \rho'_{\phi} \right\|^2_{L^2([0, 2 \pi])} \lesssim \left\| \kappa' - \overline{\kappa'} \right\|^2_{L^2(\Gamma')}.
    \end{equation}
    In light of the smallness condition on $\left\| \rho - 1 \right\|_{C^{1}([0, 2 \pi])}$ , there holds
    \begin{equation}
        \lambda = \frac{2 \pi }{\lth(\Gamma)} \sim 1,
    \end{equation}
    so that a scaling analysis implies
    \begin{equation}
        \left\| \rho_{\phi} \right\|^2_{L^2([0, 2 \pi])} \lesssim \left\| \kappa - \overline{\kappa} \right\|^2_{L^2(\Gamma)},
    \end{equation}
    which establishes the result.
\end{myproof}

Another essential tool in the proof of \eqref{eq:energydissspoiler} and beyond is Lemma \ref{lm:nearlyspherical} below that is an immediate consequence of the following (truncated) theorem from Fuglede's theory of nearly spherical surfaces see \cite[Theorem 1.2]{[Fuglede1989]}:
\begin{theorem}[Isoperimetric stability of nearly spherical surfaces]
    Let $D \in \R^{n+1}$ be a compact set with Lipschitz boundary and volume  $ \left| D \right|$ equal to that of the unit ball $B_1(0) \in \R^{n+1}$. Assume further that $D$ is starshaped with respect to its barycenter
    \begin{equation}
        c = \frac{1}{ \left| D \right|} \int_{D} \vek{x} \,dV,
    \end{equation}
    which is assumed to coincide with the origin, so that the boundary $\partial \Omega$ may be parame\-trized by
    \begin{equation}
        r\left(\xi\right) := 1 + u(\xi)
    \end{equation}
    with $\xi \in S^{n}$ and $u$ Lipschitz. Let furthermore $u$ satisfy
    \begin{equation}
    \begin{aligned}
        \left\| u \right\|_{L^{\infty}(S^{n})} &\leq \frac{3}{20(n+1)} \quad \text{ and } \quad 
        \left\| \nabla u \right\|_{L^{\infty}(S^{n})} \leq 1/2.
    \end{aligned}
    \label{eq:fugledecondition}
    \end{equation}
    Then, the  nondimensional isoperimetric deficit
    \begin{equation}
        \Delta := \frac{ \left| S^{n} \right|}{\left( n+1 \right)  \left| B(1,0) \right|} - 1,
    \end{equation}
    satisfies the stability relation
    \begin{equation}
        \frac{1}{10} \left( \left\| u \right\|^2_{L^2(S^{n})} + \left\| \nabla u \right\|^2_{L^2(S^{n})} \right) \leq \Delta \leq \frac{3}{5} \left\| \nabla u \right\|^2_{L^2(S^{n})}.
        \label{eq:fuglederesult}
    \end{equation}
    \label{thm:fuglederesult}
\end{theorem}
\begin{lemma}[Isoperimetric stability of nearly circular curves] 
    Let $\rho \in \mathcal{R}_R$ parametrize a nearly circular curve. Then $\rho$ satisfies the estimates
    \begin{equation}
        \frac{1}{R} \left\| \rho_{\phi}(\phi) \right\|^2_{L^2([0, 2 \pi])} \sim E,
        \label{eq:nearlysphericalestimateone}
    \end{equation}
    \begin{equation}
        \frac{1}{R} \left\| \rho(\phi) - R \right\|^2_{L^2([0, 2 \pi])} \lesssim E
        \label{eq:nearlysphericalestimatetwo}
    \end{equation}
    and
    \begin{equation}
        \frac{1}{R} \int_{[0, 2 \pi]} \left( \llm(\phi) - R \right)^2 \lesssim E.
        \label{eq:nearlysphericalestimatethree}
    \end{equation}
    \label{lm:nearlyspherical}
\end{lemma}
\begin{myproof}
    The estimates for \eqref{eq:nearlysphericalestimateone} and \eqref{eq:nearlysphericalestimatetwo} follow immediately from first scaling to $R'=1$, applying \eqref{eq:fuglederesult} and accounting for the scaling of the respective norms and the energy gap. The final inequality follows from
    \begin{equation}
    \begin{aligned}
        \frac{1}{R} \int_{L^2[0, 2 \pi]} \left( \llm(\phi) - R \right)^2 \,d\phi 
        &\leq \frac{2}{R} \left( \left\| \rho - R \right\|^2_{L^2([0, 2 \pi])} + \left\| \rho_{\phi} \right\|^2_{L^2([0, 2 \pi])} \right)
        \overset{\eqref{eq:nearlysphericalestimatetwo}, \eqref{eq:nearlysphericalestimateone}}{\lesssim} E.
    \end{aligned}
    \end{equation}
\end{myproof}

The last ingredient in the proof of \eqref{eq:energydissspoiler} is the well-known trace estimate, which couples the Dirichlet energy of a function on $\mathcal{M}$ to the corresponding $\dot{H}^{1/2}$-norm on $\Gamma$:
\begin{theorem}[Trace Estimate]
    Let $\Gamma$ be a nearly circular curve, let $f \in C(\mathcal{M})$ be harmonic on $\mathcal{M} \setminus \Gamma$ and let its trace satisfy $ \left. f \right|_{\Gamma} \in \dot{H}^{1/2}(\Gamma)$. Then the following estimate holds:
    \begin{equation}
        \int_{\mathcal{M} \setminus \Gamma}  \left| \nabla f  \right|^2 \,dV \sim \left\| f \right\|^{2}_{\dot{H}^{1/2}(\Gamma)}
        \label{eq:dissipationtrace}
    \end{equation}
    \label{thm:traceestimate}
\end{theorem}
The proof can be found in the appendix.

Combining the previous results we arrive at:
\begin{corollary}[Algebraic relation between $E$ and $D$]
    For any nearly circular curve para\-metrized by some $\rho \in \mathcal{R}_R$ the corresponding energy gap and dissipation satisfy
    \begin{equation}
        E \lesssim R^3 D.
        \label{eq:energydissipationrelation}
    \end{equation}
    \label{cor:energydissipation}
\end{corollary}

\begin{myproof}
    Viewing $\kappa - \overline{\kappa}$ as a $\lth(\Gamma)$-periodic function with zero mean, we combine the Poincaré-Wirtinger inequality \eqref{eq:poincarewirt} and \eqref{eq:curvelength} to obtain
    \begin{equation}
        \left\| \kappa - \overline{\kappa} \right\|^2_{L^2(\Gamma)} \leq \frac{\lth(\Gamma)}{2 \pi} \left\| \kappa - \overline{\kappa} \right\|^2_{\dot{H}^{1/2}(\Gamma)} \overset{\eqref{eq:curvelength}}{\lesssim} R \left\| \kappa - \overline{\kappa} \right\|^2_{\dot{H}^{1/2}(\Gamma)}.
        \label{eq:ltwokappatoonelhalf}
    \end{equation}
    The above chain of inequalities can be used to infer
    \begin{align}
        E \overset{\eqref{eq:nearlysphericalestimateone}}&{\sim} \frac{1}{R} \left\| \rho_{\phi}(\phi) \right\|^2_{L^2([0, 2 \pi])} \overset{\eqref{eq:alexandrov}}{\lesssim} R^2 \left\| \kappa - \overline{\kappa} \right\|^2_{L^2(\Gamma)} 
        \overset{\eqref{eq:ltwokappatoonelhalf}}{\lesssim} R^3 \left\| \kappa - \overline{\kappa} \right\|^2_{\dot{H}^{1/2}(\Gamma)}
        \overset{\eqref{eq:dissipationtrace}}{\sim} R^3  D,
    \end{align}
    where we have used the isoperimetric stability property of nearly circular curves in the first line, the curvature control of the graph in the second line and the trace estimate in the last line.
\end{myproof}

\subsection{Proof of the interpolation inequality}
The focus of this section will be further auxiliary relations between various Sobolev norms of the radial function and $H$, $E$ and $D$, for the proof of the desired interpolation inequality \eqref{eq:energyinterpolatonspoiler}:
\begin{prop}[Interpolation of $E$ in terms of $H$ and $D$]
    For a nearly circular curve parametrized by a radial function $\rho \in \mathcal{R}_R$, the following estimate holds:
    \begin{equation}
        E \lesssim \left( H D \right)^{1/2}.
        \label{eq:energygapinterpol}
    \end{equation}
    \label{prop:interpolation}
\end{prop}

The following lemma establishes control of the slope of a nearly circular radial function by the nondimensional quantity $E^2 D$, cf. \eqref{eq:slopeboundspoiler}:
\begin{lemma}[Global bound on $\rho$ and  $\rho_{\phi}$]
    Any nearly circular radial function $\rho \in \mathcal{R}_R$ satsifies the estimates
    \begin{equation}
        \left\| \rho - R \right\|_{L^{\infty}([0, 2 \pi])} \lesssim \left( E R \right)^{1/2}
        \label{eq:boundedrho}
    \end{equation}
    and
    \begin{equation}
       \frac{1}{R} \left\| \rho_{\phi} \right\|_{L^{\infty}([0, 2 \pi])}   \lesssim \left( E^2D \right)^{1/6}.
       \label{eq:boundedrhoslope}
    \end{equation}
    \label{lemma:globalbounds}
 \end{lemma}
 \begin{myproof}
    Bonnesen's inequality states that the curve $\Gamma$ is contained in an annulus with radii $R_{\text{in}},R_{\text{out}} \in \R_{>0}$ satisfying
    \begin{equation}
        \pi^2 (R_{\text{out}} - R_{\text{in}})^2 \leq \mathscr{L}(\Gamma)^2 - (2 \pi R)^2,
        \label{eq:bonnesen}
    \end{equation}
    cf.\cite[Theorem 4]{[Osserman1979]}. Together with 
    \begin{equation}
        \left\| \rho - R \right\|_{L^{\infty}([0, 2 \pi])} \leq  R_{\text{out}} - R_{\text{in}}
    \end{equation}
    and
    \begin{equation}
        \mathscr{L}(\Gamma)^2 - (2 \pi R)^2 = \left( \mathscr{L}(\Gamma) - 2 \pi R \right) \left( \mathscr{L}(\Gamma) + 2 \pi R \right) \overset{\mathscr{L}(\Gamma) \sim R}{\sim} E R
    \end{equation}
    this establishes the first claim.

    Notice first that
    \begin{equation}
        \left\| \rho_{\phi} \right\|_{L^{\infty}([0, 2 \pi])} \sim R \left\| \frac{\rho_{\phi}}{\rho} \right\|_{L^{\infty}([0, 2 \pi])}
    \end{equation}
    holds by \eqref{eq:admissannulusbound}. Applying a Gagliardo-Nirenberg (G.N.) type inequality to $\left\| \frac{\rho_{\phi}}{\rho} \right\|_{L^{\infty}([0, 2 \pi])}$ we obtain
    \begin{align}
        \left\| \rho_{\phi} \right\|_{L^{\infty}([0, 2 \pi])} \sim R \left\| \frac{\rho_{\phi}}{\rho} \right\|_{L^{\infty}([0, 2 \pi])} \overset{\text{G.N.}}&{\lesssim} R \left( \left\| \frac{\rho_{\phi}}{\rho} \right\|^2_{L^2([0, 2 \pi])} \left\| \frac{\partial}{\partial \phi} \frac{\rho_{\phi}}{\rho} \right\|^2_{L^2([0, 2 \pi])} \right)^{1/2} \\
        \overset{\eqref{eq:admissannulusbound}}&{\sim} \left( \left\| \rho \right\|^2_{L^2([0, 2 \pi])} \left\| \frac{\partial}{\partial \phi}  \frac{\rho_{\phi}}{\rho} \right\|^2_{L^2([0, 2 \pi])} \right)^{1/2}
        \overset{\eqref{eq:hminusoneomega}, \eqref{eq:ltwoomegatwo}}{\lesssim} R \left( E^2 D \right)^{1/3}.
    \end{align}
 \end{myproof}

In order to prove Proposition \ref{prop:interpolation} we will rely on two results whose proofs we defer to the end of the subsection:
\begin{corollary}[Positive Sobolev Norm Estimates for $\rho^2 - R^2$.]
    For a nearly circular curve parametrized by a radial function $\rho \in \mathcal{R}_R$, the relation  
    \begin{equation}
       \frac{1}{R^3}\left\| \rho^2 - R^2 \right\|_{\dot{H}^{1}([0, 2 \pi])}^2 \sim E
       \label{eq:Honeforconservedrho}
    \end{equation}
    holds.
    There further exists a universal $\varepsilon \in \R_{>0}$ such that
    \begin{equation}
     E^2D \leq \varepsilon
     \label{eq:EDsmallconservedrho}
  \end{equation}
  implies
    \begin{equation}
     \frac{1}{R^{5}} \left\| \rho^2 - R^2\right\|_{\dot{H}^{2}([0, 2 \pi])}^2 \lesssim \left( E D^2 \right)^{1/3}.
     \label{eq:Htwoforconservedrho}
  \end{equation}
    \label{cor:sobolevconservedrho}
\end{corollary}
\begin{corollary}[Negative norm of $\rho^2-R^2$]
    Let $\rho \in \mathcal{R}_R$ parametrize a nearly circular curve. There exists a universal $\varepsilon \in \R_{>0}$ such that if
    \begin{equation}
        E^2 D \leq \varepsilon
        \label{eq:edsmallforhonehalftwo}
    \end{equation}
    holds, the estimate
    \begin{equation}
        \left\| \rho^2- R^2 \right\|^2_{\dot{H}^{-1/2}([0, 2 \pi])} \lesssim H
        \label{eq:conservedonehalf}
    \end{equation}
    holds.
    \label{cor:hminusonehalf}
\end{corollary}
\begin{myproof}[Proof of Proposition \ref{prop:interpolation}]
    The proof is based on interpolation and the previous estimates: More precisely, we estimate 
    \begin{align}
        E \overset{\eqref{eq:Honeforconservedrho}}&{\sim} \frac{1}{R^3} \left\| \rho^2 - R^2 \right\|^2_{\dot{H}^{1}([0, 2 \pi])} 
        \overset{\eqref{eq:interpolation}}{\lesssim} \frac{1}{R^3} \left( \left\| \rho^2 - R^2 \right\|^2_{\dot{H}^{-1/2}([0, 2 \pi])} \right)^{2/5} \left( \left\| \rho^2 - R^2 \right\|^2_{\dot{H}^{2}([0, 2 \pi])} \right)^{3/5} \\
        &= \left( \left\| \rho^2 - R^2 \right\|^2_{\dot{H}^{-1/2}([0, 2 \pi])} \right)^{2/5} \left( \frac{1}{R^{5}} \left\| \rho^2 - R^2 \right\|^2_{\dot{H}^{2}([0, 2 \pi])} \right)^{3/5}
        \overset{\eqref{eq:conservedonehalf}, \eqref{eq:Htwoforconservedrho}}{\lesssim} H ^{2/5} \left( \left( ED^2 \right)^{1/3} \right)^{3/5}.
    \end{align}
    Algebraic manipulations then yield the desired result.
\end{myproof}
\begin{remark}
    The advantage of $\rho^2 - R^2$ over $\rho$ is that it has zero mean due to conservation of mass, allowing for the use of negative-order Sobolev norms.
    $\hfill \bigtriangleup$
\end{remark}

It remains to prove Corollary \ref{cor:sobolevconservedrho} and Corollary \ref{cor:hminusonehalf}. To this end, we begin with we start with estimates that couple norms of the radial function $\rho$ and the localized tangential angle $\omega$ to the energy gap $E$ and the dissipation $D$:

\begin{lemma}
    Let $\rho \in \mathcal{R}_R$ parametrize a nearly circular curve. Then the localized angle $\omega$ and the deviation of curvature from average mean curvature $\kappa - \overline{\kappa}$ of the induced curve $\Gamma$ satisfy
    \begin{alignat}{2}
        &\left\| \frac{\partial \omega}{\partial s} \right\|^2_{\dot{H}^{-1}(\Gamma)} &&\leq \left\| \omega \right\|^2_{L^{2}(\Gamma)} \sim \frac{1}{R} \left\| \rho_{\phi}(\phi) \right\|^2_{L^2([0,2 \pi])} \sim E, \label{eq:hminusoneomega} \\
       &\left\| \frac{\partial \omega}{\partial s}  \right\|^2_{L^2(\Gamma)} &&\sim \frac{1}{R} \left\| \frac{\partial}{\partial \phi} \frac{\rho_{\phi}(\phi)}{\rho(\phi)} \right\|^2_{L^2([0,2 \pi])}\lesssim \left( 1 + \left( E^2D \right)^{1/3} \right) \left( E D^2 \right)^{1/3}, \label{eq:ltwoomegatwo} \\
        &\left\| \kappa - \overline{\kappa} \right\|^2_{\dot{H}^{-1}(\Gamma)}  &&\lesssim E, \label{eq:hminusonekappa} \\
       &\left\| \kappa - \overline{\kappa} \right\|^2_{\dot{H}^{1/2}(\Gamma)} && \sim D,
       \label{eq:honehalfkappa} \\
       &\left\| \kappa - \overline{\kappa}  \right\|^2_{L^2(\Gamma)} &&\lesssim \left( E D^2 \right)^{1/3}.
       \label{eq:ltwokappa}
    \end{alignat}

    \label{lm:angleintegrals}
\end{lemma}
\begin{myproof}

    First, notice that \eqref{eq:lengthelement} and assumptions \eqref{eq:admissslopebound} and \eqref{eq:admissannulusbound} imply
    \begin{equation}
        \llm^2(\phi) \sim \rho^2(\phi)
        \label{eq:lengthrhoequivalence}
    \end{equation}
    for all $\phi \in [0, 2 \pi]$.

    For the proof of \eqref{eq:hminusoneomega} we make use of the dual formulation of the $\dot{H}^{-1}$-norm:
    \begin{equation}
        \left\| \frac{\partial \omega}{\partial s} \right\|_{\dot{H}^{-1}(\Gamma)} = \left(\int_{\Gamma} \left( | \partial_s|^{-1} \frac{\partial \omega}{\partial s} \right)^2 \,ds \right)^{1/2} = \sup_{\xi \in \dot{H}^{1}(\Gamma), \left\| \xi \right\|_{\dot{H}^{1}(\Gamma)} \leq 1} \int_{\Gamma} \frac{\partial \omega}{\partial s} \xi \,ds.
        \label{eq:dualhminusone}
    \end{equation}
    From integration by parts (I.P.) and the Cauchy-Schwarz inequality (C.S.) we obtain
    \begin{align}
        \sup_{\xi \in \dot{H}^{1}(\Gamma), \left\| \xi \right\|_{\dot{H}^{1}(\Gamma)} \leq 1} \int_{\Gamma} \frac{\partial \omega}{\partial s} \xi \,ds \overset{\text{I.P.}}&{=} \sup_{\xi \in \dot{H}^{1}(\Gamma), \left\| \xi \right\|_{\dot{H}^{1}(\Gamma)} \leq 1} \left(- \int_{\Gamma} \omega \xi_s \,ds \right) \\
        \overset{\text{C.S.}}&{\leq} \sup_{\xi \in \dot{H}^{1}(\Gamma), \left\| \xi \right\|_{\dot{H}^{1}(\Gamma)} \leq 1} \left(\int_{\Gamma} \omega^2 \,ds \int_{\Gamma} \xi_s^2 \,ds \right)^{1/2}
        \leq \left( \int_{\Gamma} \omega^2 \,ds \right)^{1/2},
    \end{align}
    which proves the first inequality in \eqref{eq:hminusoneomega}.
    
    From \eqref{eq:tanomegaisomega} and \eqref{eq:tanomega} we infer
    \begin{equation}
        \left\| \omega \right\|^2_{L^2(\Gamma)} \overset{\eqref{eq:tanomegaisomega}}{\sim} \left\| \tan (\omega) \right\|^2_{L^2(\Gamma)} \overset{\eqref{eq:tanomega}}{=} \left\| \frac{\rho_{\phi}(\phi)}{\rho(\phi)} \sqrt{\llm(\phi)} \right\|^2_{L^2([0,2 \pi ])} \sim \frac{1}{R} \left\| \rho_{\phi}(\phi) \right\|^2_{L^2([0,2 \pi ])}.
    \end{equation}

    In light of the isoperimetric stability, cf. \eqref{eq:nearlysphericalestimateone}, this completes the proof of \eqref{eq:hminusoneomega}.

    We now turn to the estimates for $\kappa - \overline{\kappa}$ before proving \eqref{eq:ltwoomegatwo}.
    First, note that both $\kappa - \overline{\kappa}$ and $\tfrac{\partial \phi}{\partial s} - \overline{\kappa}$ have vanishing integral mean. The triangle inequality implies
    \begin{equation}
    \begin{aligned}
        \left\| \kappa - \overline{\kappa} \right\|^2_{\dot{H}^{-1}(\Gamma)} & \lesssim \left\| \kappa - \frac{\partial \phi}{\partial s} \right\|^2_{\dot{H}^{-1}(\Gamma)} + \left\| \frac{\partial \phi}{\partial s} - \overline{\kappa} \right\|^2_{\dot{H}^{-1}(\Gamma)}
        \overset{\eqref{eq:hminusoneomega}}{\lesssim} E + \left\| \frac{\partial \phi}{\partial s} - \overline{\kappa} \right\|^2_{\dot{H}^{-1}(\Gamma)},
    \end{aligned}
    \end{equation}
    so that it suffices for \eqref{eq:hminusonekappa} to estimate the $\dot{H}^{-1}$-norm of $\tfrac{\partial \phi}{\partial s} - \overline{\kappa}$. We express the $\dot{H}^{-1}$-norm using the dual formulation, i.e.,
    \begin{equation*}
        \left(\int_{\Gamma} \left( | \partial_s|^{-1} \left( \frac{\partial \phi}{\partial s} - \overline{\kappa} \right) \right)^2 \,ds \right)^{1/2} = \sup_{\xi \in \dot{H}^{1}(\Gamma), \left\| \xi \right\|_{\dot{H}^{1}(\Gamma)} \leq 1} \int_{\Gamma} \left( \frac{\partial \phi}{\partial s} - \overline{\kappa} \right) \xi \,ds
  \end{equation*}
  and compute
  \begin{align*}
    \sup_{\xi \in \dot{H}^{1}(\Gamma), \left\| \xi \right\|_{\dot{H}^{1}(\Gamma)} \leq 1} \int_{\Gamma} \left( \frac{\partial \phi}{\partial s} - \overline{\kappa} \right) \xi \,ds &= \sup_{\xi \in \dot{H}^{1}(\Gamma), \left\| \xi \right\|_{\dot{H}^{1}(\Gamma)} \leq 1} \int_{\Gamma} \partial_s \left( \phi(s) - \overline{\kappa} s \right) \xi(s) \,ds \\
    \overset{\text{I.P., C.S.}}&{\leq} \left( \int_{\Gamma} \left( \phi(s) - \overline{\kappa} s \right)^2 \,ds \right)^{1/2} \\
    \overset{\text{G.B.}}&{\sim} \left( R \int_{[0,2 \pi ]} \left( \frac{2 \pi }{\lth(\Gamma)} s(\phi) - \phi \right)^2 \,d\phi \right)^{1/2},
\end{align*}
    
    where we have used the Gauss-Bonnet theorem (G.B.).
    This can be combined with a Poincaré-Wirtinger inequality, cf. \eqref{eq:poincarewirt}, with
\begin{equation}
    \int_{[0,2 \pi]} s(\phi) - \frac{\lth(\Gamma)}{2 \pi } \phi d \phi = 0
\end{equation}
and
\begin{equation}
    \partial_{\phi} \left( s(\phi) - \frac{\lth(\Gamma)}{2 \pi } \phi \right) = \llm(\phi) - \frac{\lth(\Gamma)}{2 \pi },
\end{equation}
to obtain
\begin{equation}
    \sup_{\xi \in \dot{H}^{1}(\Gamma), \left\| \xi \right\|_{\dot{H}^{1}(\Gamma)} \leq 1} \int_{\Gamma} \left( \frac{\partial \phi}{\partial s} - \overline{\kappa} \right) \xi \,ds 
    \lesssim  \left( R \left( \frac{2 \pi }{\lth(\Gamma)} \right)^2 \int_{[0,2 \pi ]} \left( \llm(\phi) - \frac{\lth(\Gamma)}{2 \pi } \right)^2 \,d\phi \right)^{1/2}.
    \label{eq:expandederrorHminusone}
\end{equation}

Inserting $\lth(\Gamma) = 2 \pi R + E$ gives
\begin{equation}
    \int_{[0,2 \pi ]} \left( \llm(\phi) - \frac{\lth(\Gamma)}{2 \pi } \right)^2 \,d\phi 
    \leq \int_{[0,2 \pi ]} \left( \llm(\phi) -  R \right)^2 \,d\phi, \label{eq:sphericaldeviationestimate}
\end{equation}
where we have used
\begin{equation}
    \int_{[0,2 \pi ]} \left( \llm(\phi) -  R \right) \frac{E}{2 \pi} \,d\phi = \frac{E}{2 \pi} \int_{[0,2 \pi ]} \left( \llm(\phi) -  R \right)  \,d\phi = \frac{E}{2 \pi} E.
\end{equation}
Thus, the estimate 
\begin{equation}
    \sup_{\xi \in \dot{H}^{1}(\Gamma), \left\| \xi \right\|_{\dot{H}^{1}(\Gamma)} \leq 1} \int_{\Gamma} \left( \frac{\partial \phi}{\partial s} - \overline{\kappa} \right) \xi \,ds 
    \lesssim  \left( R \left( \frac{2 \pi }{\lth(\Gamma)} \right)^2 \int_{[0,2 \pi ]} \left( \llm(\phi) -  R \right)^2 \,d\phi \right)^{1/2}
\end{equation}
holds, which can be combined with $R \sim \lth(\Gamma)$ and \eqref{eq:nearlysphericalestimatethree} to obtain
\begin{align}
    \sup_{\xi \in \dot{H}^{1}(\Gamma), \left\| \xi \right\|_{\dot{H}^{1}(\Gamma)} \leq 1} \int_{\Gamma} \left( \frac{\partial \phi}{\partial s} - \overline{\kappa} \right) \xi \,ds 
    \overset{R \sim \lth(\Gamma)}{\lesssim} \left( \frac{1}{R}  \int_{[0,2 \pi ]} \left( \llm(\phi) -  R \right)^2 \,d\phi \right)^{1/2}
    \overset{\eqref{eq:nearlysphericalestimatethree}}{\lesssim} E^{1/2}.
\end{align}
Inserting the above estimate into \eqref{eq:dualhminusone} establishes the desired estimate \eqref{eq:hminusonekappa}. The equivalence \eqref{eq:honehalfkappa} follows immediately from Theorem \ref{thm:traceestimate}, while the interpolation 
\begin{align}
    \left\| \kappa - \overline{\kappa} \right\|^2_{L^2(\Gamma)} \overset{\eqref{eq:interpolation}}&{\lesssim} \left( \left\| \kappa - \overline{\kappa} \right\|^2_{\dot{H}^{-1}(\Gamma)}\right)^{1/3} \left( \left\| \kappa - \overline{\kappa} \right\|^2_{\dot{H}^{1/2}(\Gamma)}\right)^{2/3}
    \overset{\eqref{eq:hminusonekappa}, \eqref{eq:honehalfkappa}}{\lesssim}\left( E D^2 \right)^{1/3}
\end{align}
establishes \eqref{eq:ltwokappa}.

Finally, it remains to establish \eqref{eq:ltwoomegatwo}:
A direct computation reveals
\begin{align}
    \left\| \frac{\partial}{\partial \phi}  \frac{\rho_{\phi}(\phi)}{\rho(\phi)} \right\|^2_{L^2([0,2 \pi])} \overset{\eqref{eq:lengthrhoequivalence}, \eqref{eq:localizedanglederivative}}&{\sim} \int_{[0,2 \pi ]} \left( \frac{\partial \omega}{\partial \phi} \right)^2 \,d\phi \overset{\llm \sim R}{\sim}  R \left\| \frac{\partial \omega}{\partial s}  \right\|^2_{L^2(\Gamma)}.
    \label{eq:ltwoomegaderiv}
\end{align}
Using the arc-length version of \eqref{eq:kappapolar}, i.e.,
\begin{equation}
    \kappa(s) = -\frac{\partial \omega}{\partial s} + \frac{\partial \phi}{\partial s}
    \label{eq:kappaarc} 
\end{equation}
for $s \in [0, \lth(\Gamma)]$, we can rewrite and split the right-hand side of \eqref{eq:ltwoomegaderiv} as
\begin{align}
    \left\| \frac{\partial \omega}{\partial s}  \right\|^2_{L^2(\Gamma)} \overset{\eqref{eq:kappaarc}}&{=} \left\| - \kappa + \frac{\partial \phi}{\partial s} - \overline{\kappa} + \overline{\kappa}  \right\|^2_{L^2(\Gamma)}
    \overset{\Delta, \,\eqref{eq:ltwokappa}}{\lesssim} \left( ED^2 \right)^{1/3} + \left\| \frac{\partial \phi}{\partial s} - \overline{\kappa} \right\|^2_{L^2(\Gamma)}. \label{eq:tangleestimate}
\end{align}

To estimate the second term on the right-hand side of \eqref{eq:tangleestimate}, the domain is first changed from $\Gamma$ to $[0, 2 \pi]$ via
\begin{equation}
    \left\| \frac{\partial \phi}{\partial s} - \overline{\kappa} \right\|^2_{L^2(\Gamma)} \sim R \left\| \frac{1}{\llm(\phi)} - \overline{\kappa} \right\|^2_{L^2([0, 2 \pi])} \overset{\text{G.B.}}{=} R \left\| \frac{1}{\llm(\phi)} - \frac{2 \pi}{\lth(\Gamma)} \right\|^2_{L^2([0, 2 \pi])}.
\end{equation}
This can further be estimated as
\begin{align}
    R \left\| \frac{1}{\llm(\phi)} - \frac{2 \pi}{\lth(\Gamma)} \right\|^2_{L^2([0, 2 \pi])} &\lesssim R \left( \int_{[0, 2 \pi]} \left( \frac{1}{\llm(\phi)} - \frac{1}{R} \right)^2 + \left( \frac{1}{R} - \frac{2 \pi }{\lth(\Gamma)} \right)^2 \,d\phi \right) \\
    \overset{R \sim \llm(\phi) \sim \lth(\Gamma)}&{\sim} R^{-3} \left( \int_{[0, 2 \pi]} \left( \llm(\phi) - R \right)^2 + 2 \pi E^2 \right) 
    \overset{\eqref{eq:nearlysphericalestimatethree}}{\lesssim} R^{-3} \left( R E + E^2 \right).
    \label{eq:rree}
\end{align}

By means of \eqref{eq:energydissipationrelation}, we can estimate powers of $E$ by corresponding powers of $R^3D$ to obtain
\begin{equation}
    RE \lesssim R^3 \left( E D^2 \right)^{1/3} \quad \text{ and } \quad 
     E^2 \lesssim R^3 \left( E^2 D \right)^{1/3} \left( E D^2 \right)^{1/3}.
\end{equation}
Substituting this into \eqref{eq:rree} and making use of the $\varepsilon$-smallness of $E^2D$ establishes \eqref{eq:ltwoomegatwo} and concludes the proof.
\end{myproof}

The next corollary allows us to replace the $L^2$-norm of $\frac{\partial}{\partial \phi} \frac{\rho_{\phi}}{\rho}$ by that of $\rho_{\phi \phi}$:
\begin{corollary}[Global bound on $\rho_{\phi}$]
    Let $\rho \in \mathcal{R}_R$ parametrize a nearly circular curve. There exists a universal $\varepsilon \in \R_{>0}$ such that if
    \begin{equation}
     E^2D \leq \varepsilon,
     \label{eq:EDsmallfour}
    \end{equation}
    there holds
    \begin{equation}
     \frac{1}{R^3}\left\| \rho_{\phi \phi} \right\|^2_{L^2([0, 2 \pi])} \sim \frac{1}{R} \left\| \partial_{\phi} \frac{\rho_{\phi}}{\rho} \right\|^2_{L^2([0, 2 \pi])} \lesssim \left( E D^2 \right)^{1/3}.
     \label{eq:htworho}
    \end{equation}
    \label{cor:rhohtwo}
 \end{corollary}
 \begin{myproof}
    We only need to show the first equivalence in \eqref{eq:htworho}, as the inequality follow from \eqref{eq:ltwoomegatwo} and \eqref{eq:EDsmallfour}:

    One direction of the equivalence is established by
    \begin{equation}
        \begin{aligned}
            \frac{1}{R} \left\| \partial_{\phi} \frac{\rho_{\phi}}{\rho} \right\|^2_{L^2([0, 2 \pi])} 
            \overset{R \sim \rho}&{\lesssim} \frac{1}{R^{5}} \left( R^2 \int_{[0, 2 \pi]} \rho_{\phi \phi}^2 \,d\phi + \left\| \rho_{\phi} \right\|^2_{L^{\infty}([0, 2 \pi])} \int_{[0, 2 \pi]} \rho_{\phi}^2 \,d\phi \right) \\
            \overset{\eqref{eq:boundedrhoslope}}&{\lesssim} \frac{1}{R^3}\left( \int_{[0, 2 \pi]} \rho_{\phi \phi}^2 \,d\phi + \left( E^2D \right)^{1/3} \int_{[0, 2 \pi]} \rho_{\phi}^2 \,d\phi\right) \\
            \overset{\text{Poincaré, } \eqref{eq:EDsmallfour}}&{\lesssim} \frac{1}{R^3} \int_{[0, 2 \pi]} \rho_{\phi \phi}^2 \,d\phi,
        \end{aligned}
        \end{equation}
    For the opposite direction, we obtain
    \begin{equation}
        \begin{aligned}
            \frac{1}{R} \left\| \partial_{\phi} \frac{\rho_{\phi}}{\rho} \right\|^2_{L^2([0, 2 \pi])} 
            &\gtrsim  \frac{1}{R^{5}} \int_{[0, 2 \pi]} \rho_{\phi \phi}^2 \rho^2 -  \frac{1}{3}\partial_{\phi}\left( \rho_{\phi}^3 \right) \rho \,d\phi
            \overset{\text{I.P., } R \sim \rho}{\gtrsim} \frac{1}{R^3} \int_{[0, 2 \pi]} \rho_{\phi \phi}^2 \,d\phi.
        \end{aligned}
    \end{equation}
    Applying\eqref{eq:ltwoomegatwo} concludes the argument.
\end{myproof}

We can now lift the estimates \eqref{eq:hminusoneomega} and \eqref{eq:ltwoomegatwo} from $\rho$ to $\rho^2 - R^2$ and prove Corollary \ref{cor:sobolevconservedrho}:
\begin{myproof}[Proof of Corollary \ref{cor:sobolevconservedrho}]
    We compute
    \begin{align}
        \frac{1}{R^3} \left\| \rho^2 - R^2 \right\|^2_{\dot{H}^{1}([0, 2 \pi])}
        \overset{\eqref{eq:admissmassconservation}}&{\sim} \frac{1}{R} \int_{[0, 2 \pi]} \rho_{\phi}^2 \,d\phi
        \overset{\eqref{eq:nearlysphericalestimateone}}{\sim} E
    \end{align}
    to obtain the first inequality. To derive \eqref{eq:Htwoforconservedrho} we calculate
    \begin{align}
        \frac{1}{R^5} \left\| \rho^2 - R^2 \right\|^2_{\dot{H}^{2}([0, 2 \pi])}
        \overset{\eqref{eq:boundedrhoslope}}&{\lesssim} \frac{1}{R^{3}} \left( \int_{[0, 2 \pi]} \rho_{\phi \phi}^2 \,d\phi + \left( E^2 D \right)^{1/3} \int_{[0, 2 \pi]} \rho_{\phi}^2 \,d\phi \right) \\
        \overset{\text{Poincaré}}&{\lesssim} \frac{1}{R^{3}} \left( 1 + \left( E^2 D \right)^{1/3} \right) \int_{[0, 2 \pi]} \rho_{\phi \phi}^2 \,d\phi
        \overset{\eqref{eq:EDsmallconservedrho}, \eqref{eq:htworho}}{\lesssim} \left( E D^2 \right)^{1/3}.
    \end{align}
\end{myproof}
The following corollary will be used below in Corollary \ref{cor:hminusonehalf} to relate the $\dot{H}^{-1/2}$-norm of $\rho^2 - R^2$ to the squared distance $H$:
\begin{corollary}
    For a nearly circular curve parametrized by a radial function $\rho \in \mathcal{R}_R$, there holds
    \begin{align}
     \int_{[0, 2 \pi]}  \left| \ln\left( \frac{\rho}{R} \right)  \right| \left( \rho^2 - R^2 \right)^2  \,d\phi
     \lesssim  \frac{1}{R^2} \left\| \rho^2 - R^2 \right\|^{2}_{\dot{H}^{-1/2}([0, 2 \pi])}  \left\| \rho^2 - R^2 \right\|^{1/2}_{\dot{H}^{1}([0, 2 \pi])} \left\| \rho^2- R^2 \right\|^{1/2}_{\dot{H}^{2}([0, 2 \pi])}.
     \label{eq:conservedinequalitylthree}
    \end{align}
    There further exists a universal $\varepsilon \in \R_{>0}$ such that if
    \begin{equation}
     E^2 D \leq  \varepsilon,
     \label{eq:conservedlthreeEDsmall}
    \end{equation}
    then
    \begin{equation}
        \int_{[0, 2 \pi]}  \left| \ln\left( \frac{\rho}{R} \right)  \right| \left( \rho^2 - R^2 \right)^2  \,d\phi \lesssim \left( E^2 D \right)^{1/6} \left\| \rho^2 - R^2 \right\|^{2}_{\dot{H}^{-1/2}([0, 2 \pi])}.
     \label{eq:conservedinequalitylthreetwo}
    \end{equation}
    \label{cor:conservedlthree}
  \end{corollary}
  \begin{myproof}
    The mean value theorem implies
    \begin{equation}
        \left| x \right| \leq 2 \left| e^{2x} -1 \right|   \quad \text{ for all } x \in [\ln(1/2), \ln(3/2)],
   \end{equation}
    which, for $x = \ln(\rho/R)$, gives
    \begin{equation}
        \left| \ln\left(\frac{\rho(\phi)}{R}\right) \right| \leq \frac{2}{R^2}  \left| \rho^2(\phi) - R^2 \right|  \quad \text{ for all } \phi \in [0, 2 \pi],
    \end{equation}
    since condition \eqref{eq:admissannulusbound} guarantees that $\ln(\rho/R)$ lies in the admissible interval.
    Thus,
    \begin{equation}
        \int_{[0, 2 \pi]}  \left| \ln\left( \frac{\rho(\phi)}{R} \right) \right| \left| \rho^2 - R^2 \right|^2 \,d\phi \lesssim \frac{1}{R^2} \left\| \rho^2 - R^2 \right\|_{L^{\infty}[0, 2 \pi]} \left\| \rho^2 - R^2  \right\|^{2}_{L^2([0, 2 \pi])}
        \label{eq:lnrhoinequality}
    \end{equation}
    holds, whose right-hand side can further be estimated via Gagliardo-Nirenberg and interpolation as
    \begin{align}
        &\left\| \rho^2 - R^2 \right\|_{L^{\infty}([0, 2 \pi])} \left\| \rho^2 - R^2  \right\|^{2}_{L^2([0, 2 \pi])}
        \overset{\text{G.-N.}}{\lesssim} \left\| \rho^2 - R^2 \right\|^{1/2}_{\dot{H}^{1}([0, 2 \pi])} \left\| \rho^2 - R^2 \right\|^{5/2}_{L^2([0, 2 \pi])} \\
        \overset{\eqref{eq:interpolation}}&{\lesssim} \left\| \rho^2 - R^2 \right\|^{1/2}_{\dot{H}^{1}([0, 2 \pi])} \left( \left\| \rho^2 - R^2 \right\|^{4/5}_{\dot{H}^{-1/2}([0, 2 \pi])} \left\| \rho^2 - R^2 \right\|^{1/5}_{\dot{H}^{2}([0, 2 \pi])} \right)^{5/2}.
    \end{align}
    Inserting this into \eqref{eq:lnrhoinequality} above establishes \eqref{eq:conservedinequalitylthree}. With the smallness assumption on $E^2D$, \eqref{eq:Honeforconservedrho} and \eqref{eq:Htwoforconservedrho} can be inserted above to obtain \eqref{eq:conservedinequalitylthreetwo}.
  \end{myproof}

  Proving control of the $\dot{H}^{-1/2}$-norm of $\rho^2-R^2$ through the squared distance $H$ concludes the section:
\begin{myproof}[Proof of Corollary \ref{cor:hminusonehalf}]
    The proof is based on the dual formulation of negative Sobolev norms. We start with the representation
    \begin{equation}
        \left\| \rho^2 - R^2 \right\|_{\dot{H}^{-1/2}([0, 2 \pi])} = \sup_{ \left\| \zeta \right\|_{\dot{H}^{1/2}([0, 2 \pi])} \neq 0} \frac{\int_{[0, 2 \pi]} \left( \rho^2(\phi) -R^2 \right) \zeta(\phi) \,d\phi}{ \left\| \zeta \right\|_{\dot{H}^{1/2}([0, 2 \pi])}}.
        \label{eq:dualityhonehalfone}
    \end{equation}
    By means of the trace estimate for the circle, cf. \eqref{eq:traceequalcircle}, we can replace the $\dot{H}^{1/2}([0, 2 \pi])$-norm of $\zeta$ by the $\dot{H}^{1}$-norm of a function $\xi \in \dot{H}^{1}(\mathbb{T}^2_{2L}\setminus \partial B_R(c(t))) $ which solves
    \begin{equation}
        \begin{aligned}
            \Delta \xi &= 0 &&\quad \text{ in } \mathbb{T}^2_{2L}\setminus \partial B_R(c(t)) \\
                   \xi &= \zeta \circ \gamma_{\text{pol}} &&\quad \text{ on } \partial B_R(c(t)), 
        \end{aligned}
    \label{eq:dirichletproblemdisk}
    \end{equation}
    i.e., $\xi$ is a harmonic extension of $\zeta$ to the complement of a circle of radius $R$ centered at $c(t)$. To keep notation shorter, we introduce
    \begin{equation}
        \mathcal{D} = \mathbb{T}^2_{2L} \setminus \partial B_R(c(t)).
    \end{equation}
    Notice that the dilation invariance of both the $\dot{H}^{1/2}$-norm and the Dirichlet energy, together with the trace estimate, implies 
    \begin{equation}
        \left\| \zeta \right\|_{\dot{H}^{1/2}([0, 2 \pi])} \sim \left\| \nabla \xi  \right\|_{L^{2}(\mathcal{D})},
    \end{equation}
    so that we infer
    \begin{equation}
        \left\| \rho^2 - R^2 \right\|_{\dot{H}^{-1/2}([0, 2 \pi])} \sim \sup_{ \left\| \xi \right\|_{\dot{H}^{1}(\mathcal{D})} \neq 0} \frac{\int_{[0, 2 \pi]} \left( \rho^2(\phi) -R^2 \right) \mathring{\xi}(R,\phi) \,d\phi}{ \left\| \xi \right\|_{\dot{H}^{1}(\mathcal{D})}}
        \label{eq:dualityhonehalftwo}
    \end{equation}
    from \eqref{eq:dualityhonehalfone}, where we have used the convenient shorthand
    \begin{equation}
        \mathring{\xi}(r,\phi) := \xi \circ P (r,\phi)
        \label{eq:shorthandpolar}
    \end{equation}
    with the polar parametrization $P$ from \eqref{eq:polarparam}. Another subsequently used shorthand is
    \begin{equation}
        \{\rho > R\} := \{\phi \in [0, 2 \pi]: \rho(\phi) > R\}.
    \end{equation}
    The numerator of the right-hand-side of \eqref{eq:dualityhonehalftwo} can be rewritten as
    \begin{equation}
    \begin{aligned}
        \int_{[0, 2 \pi]} \left( \rho^2(\phi) -R^2 \right) \mathring{\xi}(R,\phi) \,d\phi =  &2\int_{\{\rho>R\} } \int_{R}^{\rho} \mathring{\xi}(R,\phi) r \,dr d\phi
        - 2\int_{\{\rho < R\} } \int_\rho^{R} \mathring{\xi}(R,\phi) r \,dr d\phi.
        \label{eq:hminusonehalfone}
    \end{aligned}
    \end{equation}
    We replace (at a cost to be estimated below) the trace data $\mathring{\xi}(R,\phi)$ in the area integral by the full harmonic extension $\overline{\xi}(r,\phi)$ and estimate
    \begin{align}
        &\int_{\{\rho>R\} } \int_{R}^{\rho} \mathring{\xi}(r,\phi) r \,dr d\phi - \int_{\{\rho < R\} } \int_\rho^{R} \mathring{\xi}(r,\phi) r \,dr d\phi
        = \int_{\mathbb{T}^2_{2L}} \left( \chi_{\Omega_{\text{in}}} - \chi_{\Omega_{B_R(c(t))}}  \right) \xi \,dV \\
        \overset{\eqref{eq:harmonicdistance}}&{=} -\int_{\mathbb{T}^2_{2L}} \Delta \vartheta \xi \,dV
        \overset{\text{I.P., C.S.}}{\leq}  H^{1/2} \left( \int_{\mathbb{T}^2_{2L}} \left| \nabla \xi \right|^2 \,dV \right)^{1/2}. \label{eq:hminusonehalftwo}
    \end{align}
    Next, we estimate the error made by replacing $\mathring{\xi}(R,\phi)$ with $\mathring{\xi}(r,\phi)$, beginning with the integral over $\{\rho > R\}$:
    \begin{align}
        \int_{\{\rho > R\}} \int_{R}^{\rho} \left( \mathring{\xi}(r,\phi) - \mathring{\xi}(R,\phi)  \right)r \,dr d\phi
        &= \int_{\{\rho > R\}} \int \int \chi_{\{R < r < \rho(\phi)\} } \chi_{\{R < r' < r\} } \partial_{r} \mathring{\xi}(r',\phi) \,dr' r \,dr d\phi \\
        \overset{\text{Fubini}}&{=} \int_{\{\rho > R\}} \int \int \chi_{\{R < r' < \rho(\phi)\} } \chi_{\{r' < r < \rho(\phi)\} } \partial_{r} \mathring{\xi}(r',\phi) r \,dr dr' d\phi \\
        &= \frac{1}{2} \int_{\{\rho > R\}} \int_{R}^{\rho} \sqrt{r'} \partial_{r}  \mathring{\xi}(r',\phi) \frac{1}{\sqrt{r'}} \left( \rho^2 - \left( r' \right)^2 \right) \,dr' d\phi \\
        \overset{\text{C.S.}}&{\leq} \left( \int_{\{\rho > R\}} \int_{R}^{\rho}  \left|   \partial_{r} \mathring{\xi}(r',\phi) \right|^2 r' \,dr' d\phi \int_{\{\rho > R\}} \int_{R}^{\rho} \left( \rho^2 - \left( r' \right)^2 \right)^2 \frac{1}{r'} \,dr' d\phi \right)^{1/2} \\
        \overset{R \leq r'\leq \rho}&{\leq} \left( \int_{\{\rho > R\}} \int_{R}^{\rho}  \left|   \partial_{r} \mathring{\xi}(r',\phi) \right|^2 r' \,dr' d\phi \int_{\{\rho > R\}} \int_{R}^{\rho} \left( \rho^2 - R^2 \right)^2 \frac{1}{r'} \,dr' d\phi \right)^{1/2} \\
        &\lesssim \left( \int_{\mathbb{T}^2_{2L}} \left| \nabla \xi \right|^2 \,dV  \int_{\{\rho > R\}} \left( \rho^2 - R^2 \right)^2 \ln(\frac{\rho}{R}) \,d\phi \right)^{1/2}.
    \end{align}
    A similar computation for $\{\rho < R\}$ yields -- notice that $\ln(\rho/R) < 0$ holds in this case --
    \begin{equation}
        \int_{\{\rho < R\}} \int_{\rho}^{R} \left(\mathring{\xi}(r,\phi) - \mathring{\xi}(R,\phi) \right) r \,dr d\phi \lesssim \left( \int_{\mathbb{T}^2_{2L}} \left| \nabla \xi \right|^2 \,dV  \int_{\{\rho < R\}} -\left( \rho^2 - R^2 \right)^2 \ln(\frac{\rho}{R}) \,d\phi \right)^{1/2},
    \end{equation} 
    so that we arrive at
    \begin{align}
        &\int_{\{\rho > R\}} \int_{R}^{\rho} \left(\mathring{\xi}(r,\phi) - \mathring{\xi}(R,\phi) \right)r \,dr d\phi - \int_{\{\rho < R\}} \int_{\rho}^{R} \mathring{\xi}(r,\phi) - \mathring{\xi}(R,\phi) r \,dr d\phi \\
        &\lesssim \left( \int_{\mathbb{T}^2_{2L}} \left| \nabla \xi \right|^2 \,dV  \int_{[0, 2 \pi]} \left|  \ln(\frac{\rho}{R}) \right| \left( \rho^2 - R^2 \right)^2   \,d\phi \right)^{1/2} \\
        \overset{\eqref{eq:conservedinequalitylthreetwo}}&{\lesssim} \left(  \int_{\mathbb{T}^2_{2L}} \left| \nabla \xi \right|^2 \,dV\right)^{1/2} \left( E^2D \right)^{1/12} \left\| \rho^2 - R^2 \right\|_{\dot{H}^{-1/2}([0, 2 \pi])}. \label{eq:hminusonehalfthree}
    \end{align}
    
    Expanding  \eqref{eq:hminusonehalftwo} as
    \begin{align}
        &\int_{[0, 2 \pi]} \left( \rho^2(\phi) -R^2 \right) \mathring{\xi}(R,\phi) \,d\phi \\
            &=  -2\left(\int_{\{\rho>R\} } \int_{R}^{\rho} \left( \mathring{\xi}(r,\phi) - \mathring{\xi}(R,\phi) \right) r \,dr d\phi \right) + 2\left(\int_{\{\rho < R\} } \int_\rho^{R} \left( \mathring{\xi}(r,\phi) - \mathring{\xi}(R,\phi) \right) r \,dr d\phi\right) \label{eq:approximationerror}\\
            &- 2\left(\int_{\{\rho > R\} } \int_\rho^{R} \mathring{\xi}(r,\phi) r \,dr d\phi - \int_{\{\rho < R\} } \int_\rho^{R} \mathring{\xi}(r,\phi) r \,dr d\phi\right),
            \label{eq:hminusonehalffour}
    \end{align}
    and inserting into \eqref{eq:dualityhonehalftwo}, we can absorb the error introduced by the approximation, cf. the first line of \eqref{eq:approximationerror}, into the left-hand side of \eqref{eq:dualityhonehalftwo} due to $\varepsilon$-smallness of $E^2D$. The estimate \eqref{eq:hminusonehalftwo} can be applied to the remaining term to obtain the result
    \begin{equation}
        \left\| \rho^2 - R^2 \right\|_{\dot{H}^{-1/2}([0, 2 \pi])} \lesssim H^{1/2}.
    \end{equation}
    We emphasize that the exact same argument carries through in the case of the plane.
\end{myproof}

    \section{Differential relations} \label{sec:differential}
Before we state and prove the \emph{differential relations} for $H$, $E$, and $D$ in the main proposition of this section, we recall that the perimeter is non-increasing and the area is preserved under the Mullins-Sekerka flow:
\begin{lemma}
    Let $T \in \R_{>0}$ and let $\left( \Gamma(t) \right)_{t \in [0,T)}$ be a family of sufficiently smooth, simple, closed, separating curves in $\mathbb{T}^2_{2L}$ or $\R^2$ evolving according to the MS dynamics, cf. \eqref{eq:msnormal} and \eqref{eq:harmonicproblem}. 

    Then
    \begin{equation}
         \left. \frac{d \lth(\Gamma(t))}{dt} \right|_t = -D(t) \leq 0 \quad \forall t \in [0,T)
         \label{eq:perimeterderivative}
    \end{equation}
    and
    \begin{equation}
         \left. \frac{d  \left|\Omega_{\text{in}}(t) \right|}{dt} \right|_t = 0 \quad \forall t \in [0,T)
         \label{eq:areaderivative}
    \end{equation}
    hold, i.e., the evolution is perimeter-decreasing and area-preserving. Moreover, the energy gap satisfies
    \begin{equation}
     \frac{dE}{dt} = -D. \label{eq:energyderivative}
    \end{equation}
    \label{lemma:geometricvariations}
\end{lemma}
\begin{myproof}
    The first two relations can easily be checked with the help of the first variation formulas of the perimeter and the enclosed area of evolving curves, cf. \cite[Section 2]{[Pruess2016]}. The last relation is a direct consequence of inserting the definition of the energy gap into \eqref{eq:perimeterderivative}.
\end{myproof}
The following estimate, an adaptation of \cite[Lemma 7.2]{[Chen1993]} from the plane to the flat torus, provides an important improvement upon the standard Sobolev embedding of $\dot{H}^{1}$ into $L^2$ for the normal velocity $V$ in the case of small $\left\| \kappa - \overline{\kappa} \right\|_{L^{1}(\Gamma)}$:
\begin{lemma}[Improved Sobolev embedding for normal velocity]
    Let $\rho \in \mathcal{R}_R$ parametrize a nearly circular curve. Assume further that
    \begin{equation}
        \left\| \kappa - \overline{\kappa} \right\|_{L^{1}(\Gamma)} \leq \frac{1}{5}
    \end{equation}
    is satisfied.

    Then there exists a universal constant $C > 0$ such that
    \begin{equation}
        \left\| V \right\|^2_{L^2(\Gamma)} \leq \frac{1}{4 \overline{\kappa}^2\left( 1 - C \left( \left\| \kappa - \overline{\kappa} \right\|_{L^{1}(\Gamma)} + \left( R/L \right)^2 \right) \right)} \left\| V_s \right\|^2_{L^2(\Gamma)}
        \label{eq:improvednormalvelocityestimate}
    \end{equation}
    holds, where $L$ is the length scale of the torus.
    \label{lemma:improvedsobolevnormalvelo}
\end{lemma}
\begin{remark}
    Notice that the original estimate 
    \begin{equation}
        \left\| V \right\|^2_{L^2(\Gamma)} \leq \frac{1}{4 \overline{\kappa}^2\left( 1 - C \left\| \kappa - \overline{\kappa} \right\|_{L^{1}(\Gamma)}  \right)} \left\| V_s \right\|^2_{L^2(\Gamma)}.
        \label{eq:improvednormalvelocityestimatechen}
    \end{equation}
    from \cite[Lemma 7.2]{[Chen1993]} is obtained in the limit $L \to \infty$. $\hfill \bigtriangleup$
\end{remark}
\begin{myproof}
    We will point out the main modifications required to transfer the result from the plane to the torus, capitalizing on Chen's analysis whenever possible.
    
    The core idea of the underlying proof is to establish control of lower order modes in terms of higher order modes for the normal velocity  $V$ along the curve $\Gamma$. Without loss of generality, we can rescale the curve such that $\lth(\Gamma) = 2\pi$ and thus by Gauss-Bonnet $\overline{\kappa} = 1$ holds. Adopting Chen's notation, we further assume that $\Gamma$ may be parametrized via
    \begin{equation}
        X(s) := \begin{pmatrix} R(\theta(s)) \cos(\theta(s)) \\ R(\theta(s)) \sin(\theta(s)) \end{pmatrix},
        \label{eq:chenparam}
    \end{equation}
    where $s$ is the arc-length parameter, $\phi(s)$ the inclination angle, $\theta(s)$ the polar angle at arc-length $s$, and $R(\theta)$ is the radial function. Notice that the roles of $\phi$ and $\theta$ are switched with respect to our convention. 

    To obtain an estimate of the form \eqref{eq:improvednormalvelocityestimate}, it suffices to estimate the moduli of the Fourier modes $\hat{V}(1)$ and $\hat{V}(-1)$ in terms of $\left\| V \right\|_{L^2(\Gamma)}$. This follows from applying a Poincaré-Wirtinger type inequality as in
    \begin{align}
        4 \left\| V \right\|^2_{L^2(\Gamma)}
        \overset{\int_{\Gamma}Vds = 0}&{=} 4 \sum_{k \in \Z \setminus \{0\} } \left| \hat{V}(k) \right|^2 
        \overset{ 4 \leq \left| k \right|, \text{ if } 2\leq \left|k\right|}{\leq} \sum_{k \in \Z \setminus \{-1,0,1\} }  \left| k \right|^2 \left| \hat{V}(k) \right|^2 - 3 \left(  \left| \hat{V}(-1) \right|^2 + \left| \hat{V}(1) \right|^2 \right) \\
        \overset{\text{Parseval}}&{=} \left\| V_s\right\|^2_{L^2(\Gamma)} - 3 \left(  \left| \hat{V}(-1) \right|^2 + \left| \hat{V}(1) \right|^2 \right),
        \label{eq:normalveloparseval}
    \end{align}
    where $V_s$ denotes the arc-length derivative of $V$.

    The smallness of these first two non-vanishing modes is established by relating them to the first two modes of the curvature function by means of single-layer potential theory. This will require the fundamental solution of the Laplacian on the torus, cf. Lemma \ref{lemma:periodicpotential}. We will subsequently identify $\R^2$ with $\C$ via the mappings
    \begin{equation}
        z = x + iy
    \end{equation}
    and
    \begin{equation}
        \vek{x} = \begin{pmatrix} x \\ y \end{pmatrix} = \begin{pmatrix} \Re(z) \\ \Im(z) \end{pmatrix}
    \end{equation}
    for $z \in \C$ and $\vek{x} = (x,y) \in \R^2$. The quantity $\tilde{L}$ denotes the length scale of the torus after rescaling so that $\lth(\Gamma) = 2 \pi$.
    
    As Chen did in the plane, we solve the problem of harmonically extending the curvature $\kappa$ to the bulk via single-layer potential theory: The function
    \begin{equation}
        u_{1}(\vek{x}) = \frac{1}{2\pi} \int_{[0, 2 \pi]} \Lambda(\vek{x} - X(s)) V(X(s)) \,ds
        \label{eq:uone}
    \end{equation}
    is equal to the harmonic extension of curvature $u$ up to a constant $m \in \R$ according to Corollary \ref{lemma:uniquepotential}. Note that, under the assumption $\left\| V_{s} \right\|_{L^2(\Gamma)}^2 < \infty$ (otherwise, there is nothing to prove), Morrey's embedding guarantees the required Hölder regularity of $V$ to apply the corollary. Thus, the estimate 
    \begin{equation}
        \left| \widehat{u_{1}}( \pm 1) \right| \lesssim \left\| \kappa - \overline{\kappa} \right\|_{L^{1}(\Gamma)} \left\| V \right\|_{L^{1}(\Gamma)}.
        \label{eq:smallcurvaturemodeshort}
   \end{equation}
    follows from the same analysis that Chen applied in the case of the plane.

    To relate the first modes of the normal velocity to those of curvature, we \emph{linearize} the harmonic problem around the circle of radius $1$, i.e., we consider the single-layer potential
    \begin{equation}
        u_{2}(\vek{x}) = \frac{1}{4 \pi} \int_{\Gamma} \Lambda  \left( \vek{x} - \vek{y}  \right) V(\vek{y}) dS(\vek{y}) = \frac{1}{4 \pi} \int_{[0, 2 \pi]} \Lambda  \left( \vek{x} - \tilde{X}(s)  \right) V(X(s)) \,ds \quad \forall \vek{x} \in \R^2,
        \label{eq:linearizedsingle}
    \end{equation}
    where 
    \begin{equation}
        \tilde{X}: [0, 2 \pi] \to \R^2, s \mapsto \begin{pmatrix} \cos(s) \\ \sin(s) \end{pmatrix}
    \end{equation}
    is a parametrization of the unit circle $\partial B_{1}(0)$. A comparison with \eqref{eq:uone} shows that the jump data $V$ are now pulled back to the unit circle $\partial B_1(0)$. By single-layer potential theory, it follows that $u_2$ is, up to a constant, the unique harmonic function on $\R^2 \setminus \partial B_1(0)$ with jump $V$ along the interface in the normal derivatives, and continuous along $\partial B_{1}(0)$.In contrast to the case of the plane, the linearization here does not lead to the identity
    \begin{equation}
        \hat{V}(\pm 1) = 2 \widehat{u_{2}}(\pm 1),
    \end{equation}
    but rather to a result of the form
    \begin{equation}
         \left| \hat{V}(\pm 1) \right| \leq  2 \left| \widehat{u_{2}}(\pm 1)  \right| + \frac{C'}{\tilde{L}^2} \left\| V \right\|_{L^{1}(\Gamma)},
         \label{eq:utwogoal}
    \end{equation}
    for some universal constant $C' > 0$.
    
    Combining \eqref{eq:potentialseriestwo} from Corollary \ref{cor:sigmaexpansion} with the fact that $ \left| \vek{x}\right| < 2\tilde{L}$ for any $\vek{x} \in \mathbb{T}_{2\tilde{L}}$, the linearized single-layer potential can be expressed as
    \begin{equation}
        u_{2}(\vek{x}) = u_{2a}(\vek{x}) + u_{2b}(\vek{x}) + u_{2c}(\vek{x}),
        \label{eq:utwosplit}
    \end{equation}
    with the definitions
    \begin{align}
        u_{2a}(\vek{x}) &:= \frac{1}{2\pi} \int_{[0, 2 \pi]} \log(|\vek{x} - \tilde{X}(s)|) V(X(s)) \,ds, \label{eq:utwoa}\\
        u_{2b}(\vek{x})&:= \frac{1}{2\pi} \int_{[0, 2 \pi]} \left( \sum_{m,n \in \Z \setminus \{0\} } \sum_{k=3}^{\infty} - \frac{1}{k!} \left( \left( \frac{\vek{x} - \tilde{X}(s)}{z_{m,n}} \right)^{k}  + \left( \overline{\frac{\vek{x} - \tilde{X}(s)}{z_{m,n}}} \right)^{k} \right) \right) V(X(s)) \,ds, \label{eq:utwob}\\
        u_{2c}(\vek{x})&:= - \frac{1}{8 \tilde{L}^2} \int_{[0, 2 \pi]}  \left| \vek{x} - \tilde{X}(s) \right|^2 V(X(s)) \,ds.  \label{eq:utwoc}
    \end{align}
    The function $u_{2a}$ coincides with the restriction of the linearized single-layer potential in the plane to the fundamental cell of the torus. Thus, Chen's analysis for the planar case applies and the Fourier transform of $u_{2a}$ along the unit circle satisfies
    \begin{equation}
        \hat{V}(\pm 1) = 2 \widehat{u_{2a}}(\pm 1).
        \label{eq:utwoaresult}
    \end{equation}
    The first modes $\widehat{u_{2c}}(\pm 1)$ of $u_{2c}$ can be estimated as
    \begin{equation}
    \begin{aligned}
         \left|  \widehat{u_{2c}}(\pm 1) \right| &= \frac{1}{8\tilde{L}^2} \left| \int_{[0, 2 \pi]} \int_{[0, 2 \pi]}  \left| \tilde{X}(t) - \tilde{X}(s) \right|^2 V(X(s)) e^{i\pm t} \,ds dt \right| \\
         &\leq \frac{1}{4\tilde{L}^2} \sup_{s,t \in [0, 2 \pi]} \left| \tilde{X}(t) - \tilde{X}(s)  \right| \left\| V \right\|_{L^{1}(\Gamma)} = \frac{\pi}{\tilde{L}^2} \left\| V \right\|_{L^{1}(\Gamma)}.
    \end{aligned}
    \label{eq:utwocest}
    \end{equation}
    We begin the estimate for the first modes of $u_{2b}$ analogously:
    \begin{equation}
    \begin{aligned}
        &\left|  \widehat{u_{2b}}(\pm 1) \right| \leq \sup_{s,t \in [0, 2 \pi]}  \left| \sum_{m,n \in \Z \setminus \{0\} } \sum_{k=3}^{\infty} \frac{1}{k!} \left( \left( \frac{\tilde{X}(t) - \tilde{X}(s)}{z_{m,n}} \right)^{k}  + \left( \overline{\frac{\tilde{X}(t) - \tilde{X}(s)}{z_{m,n}}} \right)^{k} \right) \right| \left\| V \right\|_{L^{1}(\Gamma)}.
    \end{aligned}
    \label{eq:utwobest}
    \end{equation}
    To proceed, we fix arbitrary $s,t \in [0, 2 \pi]$ and estimate the triple series in the absolute value above. We further introduce the notation
    \begin{equation}
        x(m,n,s,t) := \frac{\tilde{X}(t) - \tilde{X}(s)}{z_{m,n}}.
    \end{equation}
    Combining invariance of the modulus under complex conjugation with the triangle inequality we obtain
    \begin{align}
        \left| \sum_{m,n \in \Z \setminus \{0\} } \sum_{k=3}^{\infty} \frac{1}{k!} \left( \left( x(m,n,s,t) \right)^{k}  + \left( \overline{x(m,n,s,t)} \right)^{k} \right) \right|
        \overset{ \Delta, \,\left| a \right| =  \left| \overline{a} \right|}{\leq} 2  \sum_{m,n \in \Z \setminus \{0\} } \sum_{k=3}^{\infty} \frac{1}{k!} \left| x(m,n,s,t)  \right|^{k}. \label{eq:triplesum}
    \end{align}
    Notice that by our choice of the length scale $\tilde{L}$ and the definition of the lattice points $z_{m,n}$,
    \begin{equation}
         \left| x(m,n,s,t) \right| \leq \frac{1}{\tilde{L}} < 1
    \end{equation}
    holds for any $m,n \in \Z \setminus \{0\}$ and $s,t \in [0, 2 \pi]$. Thus, the inner series can be estimated via
    \begin{align}
        \sum_{k=3}^{\infty} \frac{1}{k!} \left| x(m,n,s,t) \right|^{k}
        \leq \frac{\left| x(m,n,s,t) \right|^{3}}{3!} \left(  \sum_{k=0}^{\infty} \frac{\left| x(m,n,s,t) \right|^{k} }{k!} \right)
        \leq \frac{e^{1/\tilde{L}}}{3!} \left| x(m,n,s,t) \right|^{3}
        \overset{\tilde{L} > 1}{\leq} \frac{1}{2} \left| x(m,n,s,t) \right|^{3}.
    \end{align}
    In light of \eqref{eq:triplesum}, this implies
    \begin{equation}
        \left| \sum_{m,n \in \Z \setminus \{0\} } \sum_{k=3}^{\infty} \frac{1}{k!} \left( \left( x(m,n,s,t) \right)^{k}  + \left( \overline{x(m,n,s,t)} \right)^{k} \right) \right| \leq \sum_{m,n \in \Z \setminus \{0\} } \left| x(m,n,s,t) \right|^{3}.
        \label{eq:triplesumtwo}
    \end{equation}
    As the exponent of the summands exceeds $2$, the lattice sum on the right-hand side of \eqref{eq:triplesumtwo} can be bounded from above via
    \begin{align}
        \sum_{m,n \in \Z \setminus \{0\} } \left| x(m,n,s,t) \right|^{3}
        \overset{ \left| \tilde{X} \right| \leq 1}&{\leq} \frac{1}{\tilde{L}^3} \sum_{m,n \in \Z \setminus \{0\} } \frac{1}{ \left(  \left| m \right| +  \left| n \right| \right)^3}
        \overset{\text{AM-GM}}{\leq} \frac{1}{\tilde{L}^3} \sum_{m,n \in \Z \setminus \{0\} } \frac{1}{ \left|m\right|^{3/2}\left|n \right|^{3/2}} \\
        &= \frac{4}{\tilde{L}^3} \left( \sum_{m=1}^{\infty} \frac{1}{ \left| m \right|^{3/2}} \right)^2 \leq \frac{\tilde{C}}{\tilde{L}^3}, \label{eq:triplesumthree}
    \end{align}
    where $\tilde{C}$ is some universal constant and AM-GM denotes the arithmetic-geometric mean inequality.

    The estimates \eqref{eq:triplesumtwo} and \eqref{eq:triplesumthree} can be plugged into \eqref{eq:utwobest} to obtain
    \begin{equation}
        \left|  \widehat{u_{2b}}(\pm 1) \right| \leq \frac{\tilde{C}}{\tilde{L}^3} \left\| V \right\|_{L^{1}(\Gamma)}.
        \label{eq:utwobfinal}
    \end{equation}

    Next, the proclaimed relation \eqref{eq:utwogoal} is established via
    \begin{align}
         \left| \hat{V}(\pm 1) \right| \overset{\eqref{eq:utwoaresult}}&{=} 2 \left|  \widehat{u_{2a}}(\pm 1) \right| \overset{\eqref{eq:utwosplit}}{=} 2 \left|  \widehat{u_{2}}(\pm 1) - \widehat{u_{2b}}(\pm 1) - \widehat{u_{2c}}(\pm 1) \right| 
         \overset{\eqref{eq:utwobfinal}, \eqref{eq:utwocest}}{\leq} 2 \left|  \widehat{u_{2}}(\pm 1) \right| + \frac{C'}{\tilde{L}^2} \left\| V \right\|_{L^{1}(\Gamma)},\quad
         \label{eq:normalvelomode}
    \end{align}
    where $C' > 0$ is again some universal constant.
    
    We proceed to estimate the linearization error. Using the definitions \eqref{eq:uone} and \eqref{eq:linearizedsingle} together with Corollary \ref{cor:sigmaexpansion}, we can express the supremum of the error between $u_{1}$ and $u_{2}$ as
    \begin{align}
        &\sup_{t \in [0, 2 \pi]}  \left| u_{1}(X(t)) - u_{2}(\tilde{X}(t))  \right| \\
        &\leq \frac{1}{2 \pi} \sup_{t,s \in [0, 2 \pi]}  \left| \Lambda(X(t) - X(s)) - \Lambda\left( \tilde{X}(t) - \tilde{X}(s) \right) \right| \left\| V \right\|_{L^{1}(\Gamma)} \\
        &\leq \frac{1}{2 \pi} \left( \sup_{t,s \in [0, 2 \pi]}  \left|   \log \left( \frac{ \left| X(t) - X(s)  \right|}{ \left| \tilde{X}(t) - \tilde{X}(s) \right|} \right)\right| \right.  \label{eq:errorone} \\
        &+ \sup_{t,s \in [0, 2 \pi]}  \left| \sum_{m,n \in \Z \setminus \{0\} } \sum_{k=3}^{\infty} - \frac{1}{k!} \left( \left( \frac{{X}(t) - {X}(s)}{z_{m,n}} \right)^{k}  + \left( \overline{\frac{X(t) - X(s)}{z_{m,n}}} \right)^{k} \right)  \right| \\
        &+ \sup_{t,s \in [0, 2 \pi]}  \left| \sum_{m,n \in \Z \setminus \{0\} } \sum_{k=3}^{\infty} - \frac{1}{k!} \left( \left( \frac{\tilde{X}(t) - \tilde{X}(s)}{z_{m,n}} \right)^{k}  + \left( \overline{\frac{\tilde{X}(t) - \tilde{X}(s)}{z_{m,n}}} \right)^{k} \right)  \right| \\
        &+ \left. \sup_{t,s \in [0, 2 \pi]}  \left| \left| X(t) - X(s) \right|^2 - \left| \tilde{X}(t) - \tilde{X}(s) \right|^2 \right| \right) \left\| V \right\|_{L^{1}(\Gamma)}.
    \end{align}

    The first term on the right-hand side of the inequality, cf. \eqref{eq:errorone}, can be estimated just as in the planar case, i.e.,
    \begin{equation}
        \frac{1}{2 \pi}  \sup_{t,s \in [0, 2 \pi]}  \left|   \log \left( \frac{ \left| X(t) - X(s)  \right|}{ \left| \tilde{X}(t) - \tilde{X}(s) \right|} \right)\right|  \lesssim \left\| \kappa - \overline{\kappa} \right\|_{L^{1}(\Gamma)}.
    \end{equation}
    All the remaining terms on the right-hand side of the inequality can be estimated in exactly the same way as for estimates \eqref{eq:utwobfinal} and \eqref{eq:utwocest}. Hence, the linearization error satisfies
    \begin{equation}
        \sup_{t \in [0, 2 \pi]}  \left| u_{1}(X(t)) - u_{2}(\tilde{X}(t))  \right| \lesssim \left( \left\| \kappa - \overline{\kappa} \right\|_{L^1(\Gamma)}  + \frac{1}{\tilde{L}^2} \right) \left\| V \right\|_{L^{1(\Gamma)}}.
        \label{eq:linerrorest}
    \end{equation}
    As a result, the first modes of normal velocity $\hat{V}(\pm 1)$ obey
    \begin{align}
         \left| \hat{V}(\pm 1) \right| \overset{\eqref{eq:normalvelomode}}&{\leq} 2 \left|  \widehat{u_{2}}(\pm 1) \right| + \frac{C'}{\tilde{L}^2} \left\| V \right\|_{L^{1}(\Gamma)} \\
         &= \frac{\sqrt{2}}{\sqrt{\pi}} \left( \left| \int_{[0, 2 \pi]} \left( u_{1}(X(s)) - u_{2}(\tilde{X}(s))  \right) \right|  + \left| \int_{[0, 2 \pi]}  u_{2}(\tilde{X}(s)) e^{\mp i s} \,ds \right| \right) + \frac{C_{1}}{\tilde{L}^2} \left\| V \right\|_{L^2(\Gamma)} \\
         \overset{\eqref{eq:linerrorest}, \eqref{eq:smallcurvaturemodeshort}}&{\lesssim} \left( \left\| \kappa - \overline{\kappa} \right\|_{L^{1}(\Gamma)}  + \frac{1}{\tilde{L}^2} \right) \left\| V \right\|_{L^2(\Gamma)}. \label{eq:finalnormalvelo}
    \end{align}
    Inserting \eqref{eq:finalnormalvelo} into \eqref{eq:normalveloparseval} and rearranging establishes the claim in the case of the torus.
\end{myproof}
We now come to the main proposition of the section:
\begin{prop}[Differential relations among $H$, $E$ and $D$]
    There exists a universal $\varepsilon  \in \R_{>0}$ such that if a given nearly circular curve parametrized by $\rho \in \mathcal{R}_R$ satisfies
    \begin{equation}
        E^2 D \leq \varepsilon \quad \text{ and } \quad \frac{R}{L} \leq \varepsilon,
        \label{eq:edsmalldd}
    \end{equation}
    then
    \begin{align}
     \frac{dH}{dt} &\lesssim \left( H D \right)^{1/2}, \label{eq:distancederivative} \\
     \frac{d D}{dt} + \left\| V_s \right\|^2_{L^2(\Gamma)} &\lesssim  ED^3 + D^{5/2}. \label{eq:dissipationderivative}
    \end{align}
    \label{prop:diffrelations}
\end{prop}
\begin{myproof}
    
    We restrict ourselves to the case of the torus, as the proof for the planar case follows analogously.

     Next, we turn to \eqref{eq:dissipationderivative} and begin by explicitly differentiating
     \begin{align}
        \frac{dD}{dt} &= \frac{d}{dt} \int_{\mathbb{T}_{2L}^2} \left|\nabla u(t,x) \right|^2 \,dV \\
        \overset{\text{Reynolds}}&{=} \int_{\Omega_{\text{in}}(t)} \frac{\partial}{\partial t}  \left|\nabla u_{\text{in}}(t,x) \right|^2 \,dV + \int_{\Omega_{\text{out}}(t)} \frac{\partial}{\partial t} \left|\nabla u_{\text{out}}(t,x) \right|^2 \,dV \\
        &+ \int_{\Gamma(t)} \left( \left|\nabla u_{\text{in}}(t,x) \right|^2 - \left|\nabla u_{\text{out}}(t,x) \right|^2 \right) V(t,s) \,ds,
    \end{align}
      where we have applied the Reynolds transport theorem, cf. \cite[Section 2.11]{[MFO2022]}. In light of the observation
      \begin{align}
        \int_{\Omega_{i}(t)} \frac{\partial}{\partial t} \left|\nabla u_{i}(t,x) \right|^2 \,dV &= 2 \int_{\Omega_{i}(t)} \left( \nabla u_i \cdot  \nabla \left( \frac{\partial u_i}{\partial t} \right) \right) \,dV
          \overset{\text{I.P.}}{=} 2  \int_{\Omega_{i}(t)}\text{div}(\frac{\partial u_i}{\partial t} \nabla u_i) - \frac{\partial u_i}{\partial t} \Delta u_i \,dV,
      \end{align}
      for $i \in \{\text{in}, \text{out}\}$, harmonicity of $u$ and the divergence theorem imply
      \begin{align}
          \frac{dD}{dt} &= \int_{\Gamma(t)} 2 \frac{du}{dt} \left( \nabla u_{\text{in}} -\nabla  u_{\text{out}}  \right) \cdot n + \left( |\nabla u_{\text{in}}|^2 - |\nabla u_{\text{out}}|^2 \right) V \,ds \\
          \overset{ \left. u \right|_{\Gamma} = \kappa}&{=}  \int_{\Gamma(t)} 2 \frac{d \kappa}{d t} \left( \nabla u_{\text{in}} -\nabla  u_{\text{out}}  \right) \cdot n + \left( |\nabla u_{\text{in}}|^2 - |\nabla u_{\text{out}}|^2 \right) V\,ds. \label{eq:labdissideriv}
      \end{align}
      Continuity over the interface gives rise to
      \begin{align}
          |\nabla u_{\text{in}}|^2 - |\nabla u_{\text{out}}|^2 &=  \left( \nabla u_{\text{in}} \cdot  n \right)^2 - \left( \nabla u_{\text{out}} \cdot  n \right)^2\\
          &= \left(\left( \nabla u_{\text{in}}\cdot  n \right) + \left( \nabla u_{\text{out}}\cdot  n \right) \right) \, \left(\left( \nabla u_{\text{in}}\cdot  n \right) - \left( \nabla u_{\text{out}}\cdot  n \right)\right).
      \end{align}
      Substituting this into \eqref{eq:labdissideriv}, we obtain the identity
      \begin{equation}
          \frac{dD}{dt} = -2\int_{\Gamma(t)} \frac{d \kappa}{d t} V + \left( \left( \nabla u_{\text{in}}\cdot  n \right) + \left( \nabla u_{\text{out}}\cdot  n \right) \right) V^2 \,ds.
      \end{equation}
      Making use of the evolution equation of curvature for a family of curves evolving with normal velocity $V$, cf. \cite[p.443]{[Kimia1990]},
      \begin{equation}
          \frac{d \kappa}{d t} = -\partial_{ss} V - \kappa^2 V,
      \end{equation}
       we find
       \begin{align}
        \frac{dD}{dt}
        \overset{\text{I.P.}}&{=} 2\int_{\Gamma(t)}  - \left( \partial_s V \right)^2 + \kappa^2 V^2 - \left( \frac{\partial u_{\text{in}}}{\partial n} + \frac{\partial u_{\text{out}}}{\partial n} \right) V^2 \,ds,
    \end{align}
       which can be rewritten as
       \begin{equation}
        \frac{dD}{dt} + 2 \left\| V_s \right\|^2_{L^2(\Gamma(t))} = 
        2\int_{\Gamma(t)}  \kappa^2 V^2 - \left( \frac{\partial u_{\text{in}}}{\partial n} + \frac{\partial u_{\text{out}}}{\partial n} \right) V^2 \,ds.
        \label{eq:dissiidentity}
       \end{equation}
    To establish the estimate \eqref{eq:dissipationderivative}, it suffices to control the terms on the right-hand side as in
    \begin{align}
        \int_{\Gamma(t)}  \kappa^2 V^2 \, ds &\leq C_{1}\left( D^{5/2} + ED^3 \right) + \left( \frac{1}{3} + \frac{1}{C_{1}} \right) \left\| V_{s} \right\|^2_{L^2(\Gamma(t))} \label{eq:dissistratone}\\
        \int_{\Gamma(t)} - \left( \frac{\partial u_{\text{in}}}{\partial n} + \frac{\partial u_{\text{out}}}{\partial n} \right) V^2 \,ds &\leq C_{2} D^{5/2} + \frac{1}{C_{2}} \left\| V_{s} \right\|^2_{L^2(\Gamma)}, \label{eq:dissipationaimtwo}
    \end{align}
    where $C_{1}$ and $C_{2}$ are (arbitrarily large) universal constants, and absorb the $\left\| V_s \right\|^2_{L^2(\Gamma)}$ terms in the left-hand side.
    
    We begin with the easier of the two estimates, which is \eqref{eq:dissipationaimtwo}. From interpolation, we obtain
    \begin{equation}
        \int_{\Gamma(t)} - \left( \frac{\partial u_{\text{in}}}{\partial n} + \frac{\partial u_{\text{out}}}{\partial n} \right) V^2 \,ds \overset{\eqref{eq:interpolation}}{\lesssim}
        \left( \left\| \frac{\partial u_{\text{in}}}{\partial n}  \right\|_{\dot{H}^{-1/2}(\Gamma(t))} +  \left\| \frac{\partial u_{\text{out}}}{\partial n}  \right\|_{\dot{H}^{-1/2}(\Gamma(t))} \right) \left\| V^2 \right\|_{\dot{H}^{1/2}(\Gamma(t))}.
        \label{eq:dissipationone}
    \end{equation}
    Both, the normal derivatives of $u$ and $V$, are first order derivatives of $u$. Considering that the dissipation is equivalent to $\left\| u \right\|^2_{\dot{H}^{1/2}(\Gamma(t))}$, one can expect that the $\dot{H}^{-1/2}$-norms of the $ \frac{\partial u_{i}}{\partial n} $ and $V$ are of order $D^{1/2}$. This is indeed the case, as we will show in the following, beginning with the estimate
     \begin{equation}
          \left\| \frac{\partial u_{i}}{\partial n}  \right\|_{\dot{H}^{-1/2}(\Gamma(t))} \lesssim D^{1/2} \quad      \text{ for } i \in \{\text{in}, \text{out}\}. \label{eq:dissizero}
     \end{equation}
     From the dual representation of negative Sobolev norms and the trace estimate (Theorem \ref{thm:traceestimate}) we obtain
     \begin{align}
        \left\| \frac{\partial u_{i}}{\partial n}  \right\|_{\dot{H}^{-1/2}(\Gamma(t))} &= \sup_{ \left\|\zeta \right\|_{\dot{H}^{1/2}(\Gamma)(t)} \neq 0} \frac{\displaystyle \int_{\Gamma(t)} \frac{\partial u_{i}}{\partial n} \zeta \,ds}{\left\| \zeta \right\|_{\dot{H}^{1/2}(\Gamma(t))}}
        \overset{\eqref{eq:dissipationtrace}}{\lesssim} \sup_{ \left\|\xi \right\|_{\dot{H}^{1}( \Omega_i(t) )} \neq 0} \frac{\displaystyle \int_{\Gamma(t)} \frac{\partial u_{i}}{\partial n}  \left. \xi \right|_{\Gamma(t)} \,ds}{\left\| \xi \right\|_{\dot{H}^{1}(\Omega_i(t))}}, \label{eq:dissione}
     \end{align}
     where $\xi$ is the harmonic extension of $\zeta$ to $\Omega_{i}(t)$ and $ \left. \xi \right|_{\Gamma(t)}$ denotes its trace.
     Applying the divergence theorem to the numerator on the right-hand-side for any such $\xi$ yields
     \begin{align}
        \int_{\Gamma(t)} \frac{\partial u_{i}}{\partial n}  \left. \xi \right|_{\Gamma(t)} \,ds 
        \overset{\text{div thm, C.S.}}&{\leq} \left( \int_{\Omega_{i}(t)} \left\| \nabla \xi \right\|^2 \,dV \int_{\Omega_i(t)} \left\| \nabla u_i \right\|^2 \,dV \right)^{1/2} 
        \leq \left\| \xi \right\|_{\dot{H}^{1}\left( \Omega_{i}(t) \right)} D^{1/2}. \label{eq:dissitwo}
     \end{align}
     Inserting \eqref{eq:dissitwo} into \eqref{eq:dissione} establishes \eqref{eq:dissizero}. This can be combined with the triangle inequality to estimate
     \begin{align}
        \left\| V \right\|_{\dot{H}^{-1/2}(\Gamma(t))} &= \left\| \frac{\partial u_{\text{out}}}{\partial n} - \frac{\partial u_{\text{in}}}{\partial n}  \right\|_{\dot{H}^{-1/2}(\Gamma(t))} 
        \overset{\Delta}{\leq} \left\| \frac{\partial u_{\text{in}}}{\partial n}  \right\|_{\dot{H}^{-1/2}(\Gamma(t))} + \left\| \frac{\partial u_{\text{out}}}{\partial n} \right\|_{\dot{H}^{-1/2}(\Gamma(t))}
        \overset{\eqref{eq:dissizero}}{\lesssim} D^{1/2}. \, \, \, \,\label{eq:dissifour}
     \end{align}
    To estimate the remaining term on the right-hand side of \eqref{eq:dissipationone}, we make use of interpolation to obtain
    \begin{align}
        \left\| V^2 \right\|_{\dot{H}^{1/2}(\Gamma(t))} \overset{\eqref{eq:interpolation}}&{\lesssim}  \left\| V \right\|_{L^{\infty}(\Gamma(t))} \left\| V \right\|^{1/2}_{\dot{H}^{1}(\Gamma(t))} \left\| V \right\|^{1/2}_{L_2(\Gamma(t))} 
        \overset{\text{G.N.}}{\lesssim}  \left\| V \right\|_{\dot{H}^{1}(\Gamma(t))} \left\| V \right\|_{L^2(\Gamma(t))} \\
        \overset{\eqref{eq:interpolation}}&{\lesssim} \left\| V \right\|_{\dot{H}^{1}(\Gamma(t))} \left( \left\| V \right\|^{2/3}_{\dot{H}^{-1/2}(\Gamma(t))} \left\| V \right\|^{1/3}_{\dot{H}^{1}(\Gamma(t))} \right)
        \overset{\eqref{eq:dissifour}}{\lesssim} \left\| V_s \right\|^{4/3}_{L^2(\Gamma(t))} D^{1/3}.\label{eq:dissiseven}
    \end{align}
    Applying the above estimate and \eqref{eq:dissizero} to \eqref{eq:dissipationone} yields
    \begin{equation}
          \left( \left\| \frac{\partial u_{\text{in}}}{\partial n}  \right\|_{\dot{H}^{-1/2}(\Gamma(t))} +  \left\| \frac{\partial u_{\text{out}}}{\partial n}  \right\|_{\dot{H}^{-1/2}(\Gamma(t))} \right) \left\| V^2 \right\|_{\dot{H}^{1/2}(\Gamma(t))} \lesssim D^{5/6} \left\| V_s \right\|^{4/3}_{L^2(\Gamma(t))},
          \label{eq:dissipationaimone}
    \end{equation}
    which leads to \eqref{eq:dissipationaimtwo} after applying Young's inequality with exponents $(3, 3/2)$.

    We now turn to the proof of \eqref{eq:dissistratone}: Due to the fact that the curvature, unlike its derivatives, does not vanish for circles, we effectively lose an order of decay for the left-hand side of \eqref{eq:dissistratone} compared with \eqref{eq:dissipationaimone}. We are essentially dealing with a quadratic form in the normal derivatives of $u$ instead of a cubic one as before, which calls for a modified approach.
    
    First, we decompose $\left\| \kappa V \right\|^2_{L^2(\Gamma(t))}$ as in
    \begin{equation}
        \left\| \kappa V \right\|^2_{L^2(\Gamma)} = \left( \left\| (\kappa - \overline{\kappa}) V \right\|^2_{L^2(\Gamma(t))} + 2 \overline{\kappa} \int_{\Gamma(t)} \left( \kappa - \overline{\kappa} \right) V^2 \,ds + \overline{\kappa}^2 \left\| V \right\|^2_{L^2(\Gamma(t))} \right).
        \label{eq:dissisix}
    \end{equation}
    Notice that the first and second term of the decomposition are quartic or cubic in derivatives of $u$ (by Poincaré), while the third term is only quadratic. To deal with the quadratic term,
    we proceed with a case distinction to resolve two different regimes for the dissipation, each of which will require a slightly different approach. Note that the strictly positive real number $\varepsilon_{D} \in \R_{>0}$ introduced below will be fixed later in the treatment of the second regime, as its size is irrelevant for the first regime.

    In the case of large dissipation, i.e.,
    \begin{equation}
        D > \varepsilon_{D} \overline{\kappa}^2,
    \end{equation}
    the inequality
    \begin{equation}
        \overline{\kappa}^2 < \frac{D}{\varepsilon_{D}}
    \end{equation}
    provides us with the estimate 
    \begin{align}
        \overline{\kappa}^2 \left\| V \right\|^2_{L^2(\Gamma(t))} < \frac{2D}{\varepsilon_{D}} \left\| V \right\|^2_{L^2(\Gamma(t))}.
        \label{eq:dissidissitwo}
    \end{align}
    This can be combined with interpolation
    \begin{align}
        \left\| V \right\|^2_{L^2(\Gamma(t))} &\lesssim \left( \left\| V \right\|^{2/3}_{\dot{H}^{-1/2}(\Gamma(t))} \left\| V_s \right\|^{1/3}_{L^2(\Gamma(t)) }  \right)^2
        \overset{\eqref{eq:dissifour}}{\lesssim} D^{2/3} \left\| V_{s} \right\|^{2/3}_{L^2(\Gamma(t))}
    \end{align}
    and Young's inequality with exponents $\left( 3/2, 3 \right)$ to arrive at
    \begin{equation}
        \overline{\kappa}^2 \left\| V \right\|^2_{L^2(\Gamma(t))} \leq C_{1} D^{5/2} +  \frac{1}{C_{1}} \left\| V_s \right\|^2_{L^2(\Gamma(t))},
        \label{eq:dissiptwo}
    \end{equation}
    where $C_{1} \in \R_{>0}$ can be chosen arbitrarily large.

    Let us now assume
    \begin{equation}
        D \leq \varepsilon_{D} \overline{\kappa}^2.
    \end{equation}
    The Gauss-Bonnet theorem implies
    \begin{equation}
        D \leq \varepsilon_{D} \overline{\kappa}^2 = \varepsilon_{D} \left( \frac{2 \pi}{\lth(\Gamma(t))} \right)^2.
    \end{equation}
    At the same time, the Cauchy-Schwarz inequality provides us with
    \begin{equation}
        \left\| \kappa - \overline{\kappa} \right\|_{L^{1}(\Gamma(t))} \leq \sqrt{\lth(\Gamma(t))} \left\| \kappa - \overline{\kappa} \right\|_{L^2(\Gamma(t))},
    \end{equation}
    which can be combined with a Poincaré-Wirtinger-inequality and the trace estimate, cf. Theorem \ref{thm:traceestimate}, to obtain
    \begin{equation}
        \left\| \kappa - \overline{\kappa} \right\|_{L^{1(\Gamma(t))}} \leq \lth(\Gamma(t)) \sqrt{D}.
    \end{equation}
    In light of the condition $D \leq \varepsilon_{D} \overline{\kappa}^2$, the above estimate yields 
    \begin{equation}
        \left\| \kappa - \overline{\kappa} \right\|_{L^{1}(\Gamma)} \leq 2 \pi \varepsilon_{D}.
    \end{equation}
    Choosing $\varepsilon_{D}$ such that
    \begin{equation}
        \varepsilon_{D} \leq \frac{1}{2 \pi} \min\left\{ \frac{1}{5}, \frac{1}{4C} \right\},
    \end{equation}
    where $C$ is the positive universal constant from Lemma \ref{lemma:improvedsobolevnormalvelo}, we thus have the necessary control over $\left\|\kappa - \overline{\kappa} \right\|_{L^{1}(\Gamma(t))}$ to make use of the aforementioned lemma. In view of the smallness of $R/L$ and the choice $\left\| \kappa - \overline{\kappa} \right\|_{L^1(\Gamma(t))} \leq \frac{1}{4C}$, the lemma can be applied to obtain
    \begin{equation}
        \overline{\kappa}^2 \left\| V \right\|^2_{L^2(\Gamma(t))} \leq \frac{1}{3} \left\| V_s \right\|^2_{L^2(\Gamma(t))}.
        \label{eq:dissiquadraticsmall}
    \end{equation}
    By \eqref{eq:dissiptwo} and \eqref{eq:dissiquadraticsmall}, the estimate
    \begin{equation}
        \overline{\kappa}^2 \left\| V \right\|^2_{L^2(\Gamma(t))} \leq C_{1} D^{5/2} + \left( \frac{1}{3} + \frac{1}{C_{1}} \right) \left\| V_s \right\|^2_{L^2(\Gamma(t))}
        \label{eq:quadraticvalid}
    \end{equation}
    is valid for both regimes of the dissipation, and it remains to estimate the first two terms on the right-hand side of \eqref{eq:dissisix}. The cubic term can can be estimated via
    \begin{align}
        \overline{\kappa} \int_{\Gamma(t)} \left( \kappa - \overline{\kappa} \right) V^2 \,ds 
        \overset{\text{C.S., G.B.}}&{\leq} \frac{2 \pi}{\sqrt{\mathscr{L}(\Gamma(t))}} \left\| \kappa - \overline{\kappa} \right\|_{L^2(\Gamma(t))} \left\| V \right\|^2_{L^{\infty}(\Gamma(t))} 
        \overset{\eqref{eq:ltwokappa}}{\lesssim} \frac{1}{\sqrt{\mathscr{L}(\Gamma(t))}} \left( ED^2 \right)^{1/6} \left\| V \right\|^2_{L^{\infty}(\Gamma(t))} \\
        \overset{\eqref{eq:energydissipationrelation}}&{\lesssim} D^{1/2} \left\| V \right\|^2_{L^{\infty}}(\Gamma(t)). 
        \label{eq:meancurvvelo}
    \end{align}
    In light of the chain of estimates
    \begin{align}
        \left\| V \right\|^2_{L^{\infty}(\Gamma(t))} \overset{\text{G.N.}}{\lesssim}  \left\| V \right\|_{L^2(\Gamma(t))} \left\| V \right\|_{\dot{H}^{1}(\Gamma(t))} 
        \overset{\eqref{eq:interpolation}}{\lesssim}  \left\| V \right\|^{2/3}_{\dot{H}^{-1/2}(\Gamma(t))} \left\| V_s \right\|^{4/3}_{L^2(\Gamma(t))} 
        \overset{\eqref{eq:dissifour}}{\lesssim} D^{1/3} \left\| V_s \right\|^{4/3}_{L^2(\Gamma(t))},
        \label{eq:veloinfinityestimate}
    \end{align}
    the relation \eqref{eq:meancurvvelo} implies
    \begin{equation}
        \overline{\kappa} \int_{\Gamma(t)} \left( \kappa - \overline{\kappa} \right) V^2 \,ds \lesssim D^{5/6} \left\| V_s \right\|^{4/3}_{L^2(\Gamma(t))}.
        \label{eq:togetherone}
    \end{equation}
    Splitting the quartic term $\left\|\left( \kappa - \overline{\kappa} \right) V \right\|^2_{L^2(\Gamma(t))}$ and estimating the factors via \eqref{eq:ltwokappa} and \eqref{eq:veloinfinityestimate} produces
    \begin{align}
        \left\| (\kappa - \overline{\kappa}) V \right\|^2_{L^2(\Gamma(t))} &\leq \left\| \kappa - \overline{\kappa} \right\|^2_{L^2(\Gamma(t))} \left\| V \right\|^2_{L^{\infty}(\Gamma(t))}
        \overset{\eqref{eq:ltwokappa}, \eqref{eq:veloinfinityestimate}}{\lesssim} E^{1/3} D \left\| V_s \right\|^{4/3}_{L^2(\Gamma(t))}.
        \label{eq:togethertwo}
    \end{align}
    Applying Young's inequality with the tuple $(3, 3/2)$ to \eqref{eq:togetherone} and \eqref{eq:togethertwo}, and inserting the result together with \eqref{eq:quadraticvalid} into \eqref{eq:dissisix} establishes \eqref{eq:dissistratone}.

    Now, we turn to the proof of \eqref{eq:distancederivative}. Differentiating the squared distance explicitly under the integral sign gives
    \begin{equation}
        \frac{d H}{d t} = \frac{d }{d t} \int_{\R^2} \nabla \vartheta \cdot \nabla \vartheta \,dV 
        = 2 \int_{\R^2} \nabla \vartheta_t \cdot \nabla \vartheta \,dV
        = -2 \int_{\R^2} \vartheta_{t} \Delta \vartheta \,dV, \label{eq:distancederivativeone}
    \end{equation}
    while an application of Reynolds transport theorem, cf. \cite[Section 2.11]{[MFO2022]}, shows that
    \begin{align}
        2 \frac{d H}{d t} &= - 2 \frac{d }{d t} \int_{\R^2} \vartheta \Delta \vartheta \,dV = - 2 \frac{d }{d t} \left( \int_{\Omega_{\text{in}}(t)} \vartheta  \,dV -\int_{B_{R}(c(t))} \vartheta \,dV \right) \\
        \overset{\text{Reynolds}}&{=}  - 2 \left( \int_{\R^2} \vartheta_t \Delta \vartheta \,dV  + \int_{\Gamma(t)} \vartheta V \,ds - \int_{\partial B_R(c(t))} \vartheta \dot{c} \cdot n_{\partial B_R(c(t))} \,ds \right) \label{eq:distancederivativetwo}
    \end{align}
    holds. Subtracting \eqref{eq:distancederivativeone} from \eqref{eq:distancederivativetwo} yields
    \begin{equation}
        \frac{d H}{d t} = - 2 \left( \int_{\Gamma(t)} \vartheta V \,ds - \int_{\partial B_R(c(t))} \vartheta \dot{c} \cdot n_{\partial B_R(c(t))} \,ds \right).
        \label{eq:distancederivativethree}
    \end{equation}
    We begin by estimating the first term on the left-hand-side of \eqref{eq:distancederivativethree}:

    From the definition of $\vartheta$, cf. \eqref{eq:harmonicdistance}, and the divergence theorem, we obtain
    \begin{align}
        - 2 \int_{\Gamma(t)} \vartheta V \,ds &= -2 \int_{\Gamma(t)} \vartheta \left( \frac{\partial u_{\text{out}}}{\partial n} - \frac{\partial u_{\text{in}}}{\partial n}   \right) \,ds \\
        \overset{\text{div thm}}&{=} 2 \int_{\mathbb{T}^2_{2L} \setminus \Gamma(t)} \text{div} \left( \nabla \vartheta u \right) - \Delta \vartheta u \,dV \\
        \overset{\text{div thm}}&{=} 2 \left( - \int_{\Gamma(t)} \left( \nabla \vartheta u_{\text{out}} \right) \cdot n \,ds +  \int_{\Gamma(t)} \left( \nabla \vartheta u_{\text{in}} \right) \cdot n \,ds  - \int_{\mathbb{T}^2_{2L} \setminus \Gamma(t)} \Delta \vartheta u \,dV\right) \\
        &= - 2 \int_{\mathbb{T}^2_{2L}} \Delta \vartheta u \,dV 
        \overset{\eqref{eq:harmonicdistance}}{=}  2 \int_{\mathbb{T}^2_{2L}} \left( \chi_{\Omega_{\text{in}}(t)} - \chi_{B_R(c(t))} \right) u \,dV.
        \label{eq:sqdistancezero}
    \end{align}
    In polar coordinates, cf. \eqref{eq:polarparam}, this can be expressed as
   \begin{align}
    \int_{\mathbb{T}^2_{2L}} \left( \chi_{\Omega_{\text{in}}(t)} - \chi_{B_R(c(t))} \right) u \,dV
    = \int_{\{ \rho >R \} } \int_{R}^{\rho(\phi)} \mathring{u}(r,\phi) r dr d\phi - \int_{\{ \rho <R \} } \int_{\rho(\phi)}^{R} \mathring{u}(r,\phi) r dr d\phi,
    \label{eq:sqdistanceone}
   \end{align}
   where we employ the shorthand
   \begin{equation}
       \mathring{u}(r,\phi) = u \circ P (r,\phi)
   \end{equation}
   defined in \eqref{eq:shorthandpolar}.

   Observing that 
   \begin{equation}
      \mathring{u}(\rho(\phi), \phi) = \kappa(\phi)
   \end{equation}
   holds for all $ \phi \in  [0, 2 \pi]$ and expanding the first integral gives
   \begin{equation}
   \begin{aligned}
    \int_{\{ \rho >R \} } \int_{R}^{\rho(\phi)} \mathring{u}(r, \phi) r \,dr d\phi
    = \int_{\{ \rho >R \} } \int_{R}^{\rho(\phi)} \left( \mathring{u}(r, \phi) - \mathring{u}(\rho(\phi), \phi) \right) r \,dr d\phi
    + \int_{\{ \rho >R \} } \int_{R}^{\rho(\phi)} \kappa(\phi) r \,dr d\phi.
   \end{aligned}
   \label{eq:urintegral}
   \end{equation}
   Equation \eqref{eq:sqdistancezero} followed by \eqref{eq:sqdistanceone} and \eqref{eq:urintegral} can be inserted into \eqref{eq:distancederivativethree} to obtain
   \begin{equation}
   \begin{aligned}
    \frac{1}{2} \frac{d H}{d t} &= \int_{\{ \rho >R \} } \int_{R}^{\rho(\phi)} \left( \mathring{u}(r, \phi) - \mathring{u}(\rho(\phi), \phi) \right) r \,dr d\phi
     - \int_{\{ \rho < R \} } \int^{R}_{\rho(\phi)} \left( \mathring{u}(r, \phi) - \mathring{u}(\rho(\phi), \phi) \right) r \,dr d\phi\\
    &+ \int_{\{ \rho >0 \} } \int_{R}^{\rho(\phi)} \kappa(\phi) rdr d\phi - \int_{\{ \rho <0 \} } \int_{\rho(\phi)}^{R} \kappa(\phi) r \,dr d\phi \\
    &+ \int_{\partial B_R(c(t))} \vartheta \dot{c} \cdot n_{\partial B_R(c(t))} \,ds.
   \end{aligned}
   \label{eq:sqdistancethree}
   \end{equation}
   The terms in the third line of \eqref{eq:sqdistancethree} can be estimated via
   \begin{align}
    \int_{\{ \rho >0 \} } \int_{R}^{\rho(\phi)} \kappa(\phi) rdr d\phi - \int_{\{ \rho <0 \} } \int_{\rho(\phi)}^{R} \kappa(\phi) r \,dr d\phi
    &= \frac{1}{2} \int_{[0, 2 \pi]} \kappa(\phi) \left(\rho(\phi)^2 - R^2\right) \,d\phi \\
    \overset{\eqref{eq:interpolation}, \eqref{eq:arcpolarcircle}}&{\sim}  \left\| \kappa \right\|_{\dot{H}^{1/2}(\Gamma(t))} \left\| \rho^2 - R^2 \right\|_{\dot{H}^{-1/2}([0, 2 \pi])} \\
    \overset{\eqref{eq:honehalfkappa}, \eqref{eq:conservedonehalf}}&{\lesssim} \left( HD \right)^{1/2}
    \label{eq:distanceone}
 \end{align}
   by including the integrals over $\{\rho = R\}$ and $\{\rho < R\}$ and integrating.

   Next, we consider the first line on the right-hand-side of \eqref{eq:sqdistancethree}. In a procedure similar to the proof of Corollary \ref{cor:hminusonehalf} we estimate
   \begin{align}
    &\left| \int_{\{ \rho >R \} } \int_{R}^{\rho(\phi)}  \left( \mathring{u}(r, \phi) - \mathring{u}(\rho(\phi), \phi) \right) r \,dr d\phi \right|\\
    \overset{\text{C.S.}}&{\leq} \frac{1}{2} \left( \int_{\{ \rho > R\} } \int_{[R, \rho(\phi)]}   \left| \partial_{r} \mathring{u}(r',\phi) \right|^2 r' \,dr' d\phi \int_{\{ \rho > R\} } \int_{[R, \rho(\phi)]}  \left( \left( r' \right)^2 - R^2 \right)^2 \frac{1}{r'} \,dr' d\phi \right)^{1/2} \\
    \overset{r' \leq \rho(\phi)}&{\leq} \frac{1}{2} \left( \int_{\mathbb{T}^2 }   \left| \nabla u \right|^2 \,dV \,  \int_{\{ \rho > R\} } \int_{[R, \rho(\phi)]}  \left( \rho(\phi)^2 - R^2 \right)^2 \frac{1}{r'} \,dr' d\phi \right)^{1/2} \\
    \overset{\eqref{eq:honehalfkappa}}&{\lesssim} \frac{D^{1/2}}{2} \left( \int_{\{ \rho > R\} } \ln\left( \frac{\rho(\phi)}{R} \right) \left( \rho(\phi)^2 - R^2 \right)^2  \,d\phi \right)^{1/2}. \label{eq:dHdtnonerror}
\end{align}

  An analogous computation for the integral over $\{\rho <R\}$ reveals 
  \begin{align}
    \left| \int_{\{ \rho <R \} } \int_{\rho(\phi)}^{R}  \left( \mathring{u}(r, \phi) - \mathring{u}(\rho(\phi), \phi) \right) r \,dr d\phi \right| 
    \lesssim \frac{D^{1/2}}{2} \left( \int_{\{ \rho < R\} } - \ln\left( \frac{\rho(\phi)}{R} \right) \left( \rho(\phi)^2 - R^2 \right)^2  \,d\phi \right)^{1/2},
  \end{align}
  so that we can estimate the integral over the combined domain via
  \begin{align}
    &\left| \int_{\{ \rho > R\} } \int_{R}^{\rho(\phi)}  \left( \mathring{u}(r, \phi) - \mathring{u}(\rho(\phi), \phi) \right) r \,dr d\phi  - \int_{\{ \rho(\phi) < R\} } \int_{\rho(\phi)}^{R}  \left( \mathring{u}(r, \phi) - \mathring{u}(\rho(\phi), \phi) \right) r \,dr d\phi \right| \\
    &\leq \frac{1}{2} D^{1/2} \left( \int_{[0, 2 \pi] }  \left| \ln\left( \frac{\rho((\phi))}{R} \right) \right| \left( \rho(\phi)^2 - R^2 \right)^2  \,d\phi \right)^{1/2} \\
    \overset{\eqref{eq:conservedinequalitylthreetwo} \eqref{eq:conservedonehalf}}&{\lesssim} D^{1/2} \left( \left( E^2 D \right)^{1/6} H \right)^{1/2} \lesssim \left( HD \right)^{1/2}. \label{eq:distancetwo}
  \end{align}
  To conclude the proof, it remains to estimate the last term on the right-hand-side of \eqref{eq:sqdistancethree}, which is a consequence of the motion of the barycenter during the evolution.  To this end, we write
    \begin{equation}
        \dot{\vek{c}} = \begin{pmatrix} \dot{c}_{(x)}\\ \dot{c}_{(y)}  \end{pmatrix} = \nabla \left( \dot{c}_{(x)} x + \dot{c}_{(y)} y \right),
    \end{equation}
    so that we may apply the divergence theorem
    \begin{align}
        \int_{\partial B_R(\vek{c}(t))} \vartheta \dot{\vek{c}} \cdot n_{\partial B_R(\vek{c}(t))} \,ds 
        \overset{\text{C-S}}&{\leq} \left( \int_{B_{R}(\vek{c}(t))}  \left| \nabla \vartheta \right|^2 \,dV \int_{B_{R}(\vek{c}(t))} \left| \dot{\vek{c}} \right|^2 \,dV \right)^{1/2}
        = H^{1/2} \left( \int_{B_{R}(\vek{c}(t))}  \left| \dot{\vek{c}} \right|^2 \,dV \right)^{1/2}\\
        \overset{\eqref{eq:translationdissipation}}&{\lesssim} H^{1/2} \left( \frac{D}{ \left|\Omega_\text{in} \right|} \int_{B_{R}(\vek{c}(t))} 1 \,dV  \right)^{1/2}
        \overset{ \left|\Omega_{\text{in}} \right| =  \left| B_R \right| }{=} \left( H D  \right)^{1/2},
        \label{eq:distancethree}
    \end{align}
    where we have made use of the estimate \eqref{eq:translationdissipation} from Lemma \ref{lemma:barycentermotion}, whose statement and proof can be found after this proof.

    The entire right-hand side of \eqref{eq:sqdistancethree} can now be estimated against (a multiple) of $\left( HD \right)^{1/2}$ by applying \eqref{eq:distanceone} to the first two lines, \eqref{eq:distancetwo} to the third line and \eqref{eq:distancethree} to the fourth line of the right-hand side of \eqref{eq:sqdistancethree}. This establishes the final claim of the proposition.
\end{myproof}
We conclude this section with an estimate for the velocity of the barycenter:
\begin{lemma}[Velocity of the barycenter]
    Let $\rho \in \mathcal{R}_R$ parametrize a nearly circular curve. There exists a universal $\varepsilon \in \R_{>0}$ such that for any $t_{1} > 0$ if
    \begin{equation}
        E^2D(t_{1}) \leq \varepsilon
    \end{equation}
    holds, the velocity of the barycenter under the MS dynamics satisfies the estimate
    \begin{equation}
        \left| \dot{\vek{c}}(t_{1}) \right|^2 \lesssim \frac{D}{ \left| \Omega_{\text{in}} \right|}.
    \label{eq:translationdissipation}
    \end{equation}
    \label{lemma:barycentermotion}
\end{lemma}
\begin{myproof}
    For any $t_{1}>0$, the bounds \eqref{eq:boundedrho} and \eqref{eq:boundedrhoslope} on $\rho$ and $\rho_{\phi}$ combined with the smallness of $E^2D$ and the continuity of $\rho$ in time imply the existence of a $\delta(t_{1})$ such that $\Gamma(t)$ can be parametrized via
\begin{equation}
    \tilde{\gamma}(t) = \vek{c}(t_{1}) + \tilde{\rho}(t,\phi) \begin{pmatrix} \cos(\phi) \\  \sin(\phi) \end{pmatrix}.
    \label{eq:frozenpole}
\end{equation}
Here $\tilde{\rho}(t)$ is an admissible radial function for all $t \in (t_{1} - \delta(t_{1}),t_{1}+\delta(t_{1}))$ satisfying 
\begin{equation}
    \tilde{\rho}(t_{1}) = \rho(t_{1}).
\end{equation}
Thus, \eqref{eq:frozenpole} parametrizes $\Gamma(t)$ during the interval $(t_{1} - \delta(t_{1}),t_{1} + \delta(t_{1}))$ while freezing the pole in time. As a consequence, the entire evolution of the curve during this period is captured by the evolution of $\tilde{\rho}$. In particular, the kinematic coupling from \eqref{eq:normalveloparam} simplifies to
\begin{equation}
    V(t,\phi) \tilde{\llm}(t,\phi) = \partial_t \tilde{\rho}(t,\phi) \tilde{\rho}(t,\phi), \label{eq:simplecoupling}
\end{equation}
where $\tilde{\llm}(t,\phi)$ is the corresponding length element
\begin{equation}
    \tilde{\llm}(t,\phi) = \sqrt{\tilde{\rho}(t,\phi)^2 + \partial_\phi \tilde{\rho}(t,\phi)^2}.
\end{equation}

Notice that the barycenter of $\tilde{\rho}(t, \cdot)$ differs from the center function $\vek{c}(t)$ exactly by $\vek{c}(t_{1})$, i,e,
\begin{equation}
    \vek{c}(t) = \vek{c}(t_{1}) + \frac{1}{\left| \Omega_{\text{in}} \right|} \begin{pmatrix} \displaystyle \int_{[0, 2 \pi]} \tilde{\rho}(f,\phi)^3 \cos(\phi) \, d\phi\\[2ex] \displaystyle \int_{[0, 2 \pi]} \tilde{\rho}(t,\phi)^3 \sin(\phi)\, d\phi \end{pmatrix},
    \label{eq:newbarycentercoords}
\end{equation}

Differentiating \eqref{eq:newbarycentercoords} in time we obtain for the velocity of the barycenter in the $x$-direction
\begin{align}
    \dot{c}_{(x)}(t) &= \frac{3}{\left| \Omega_{\text{in}} \right|} \int_{[0, 2 \pi]} \partial_t \tilde{\rho}(t,\phi) \tilde{\rho}^2(t,\phi) \cos(\phi) \,d\phi  \overset{\eqref{eq:simplecoupling}}{=} \frac{3}{\left| \Omega_{\text{in}} \right|} \int_{[0, 2 \pi]} V(t, \phi) \tilde{\llm}(t,\phi) \tilde{\rho}(t,\phi) \cos(\phi) \,d\phi \\
    &= \frac{3}{\left| \Omega_{\text{in}} \right|} \int_{\Gamma(t)} V(t, s) \left(x(t,s) -c_{(x)}(t_{1}) \right) \,ds = \frac{3}{\left| \Omega_{\text{in}} \right|} \int_{\Gamma(t)} \left( \frac{\partial u_{\text{out}}}{\partial n} - \frac{\partial u_{\text{in}}}{\partial n}  \right)(t,s) \left(x(t,s) -c_{(x)}(t_{1}) \right) \,ds.
\end{align}
Let $X(t,s)$ denote an extension of $x(t,s) - c_{(x)}(t_{1})$ to a compact radial neighborhood of $\Gamma(t)$ (the precise construction will be provided below). Then the divergence theorem can be used to establish
\begin{equation}
    \dot{c}_{(x)}(t) = \frac{-3}{\left| \Omega_{\text{in}} \right|} \int_{\R^2} \text{div} \left( \nabla \left( u(t,\sigma) \right) X(t,\sigma) \right)\, dV 
    \overset{C.S.}{\leq} D^{1/2} \frac{3}{\left| \Omega_{\text{in}} \right|} \left( \int_{\R^2} \left|\nabla X(t, \sigma) \right|^2 \,dV \right)^{1/2}.
\end{equation}
We will now show
\begin{equation}
    \int_{\R^2} \left|\nabla X(t, \sigma) \right|^2 \,dV \sim \left| \Omega_{\text{in}} \right|,\
    \label{eq:barycentervelointegral}
\end{equation}
for an adequately chosen extension $X(t, \cdot)$, from which follows that
\begin{equation}
     \left| \dot{c}_{(x)}(t) \right| \lesssim \left( \frac{D}{\left| \Omega_{\text{in}} \right|} \right)^{1/2}.
\end{equation}
Together with the analog estimate for $\left| \dot{c}_{(x)}(t) \right|$ we have estimate \eqref{eq:translationdissipation}.
\begin{figure}[ht]
    \centering
    \includegraphics[width=0.4\textwidth]{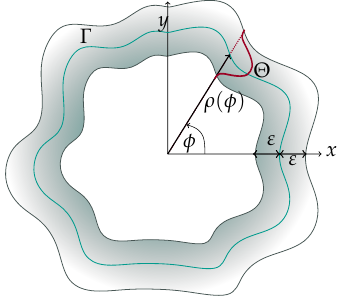}
    \caption{Radial neighborhood of $\Gamma(t)$: The support of the extension coincides with the shaded neighborhood. The bold red graph visualizes the bump function $\Theta$, cf. \eqref{eq:bumpdef} and \eqref{eq:boundednessbump}, along a radial slice of the neighborhood.}
    \label{fig:radialneighbor}
\end{figure}

It remains only to give the details of the extension $X$. We extend $x(t,s)-c_x(t_{1})$ in a compact annular neighborhood of $\Gamma(t)$. To cope with the lack of pointwise control over the curvature of $\Gamma(t)$, we make use of the available parametrization of the curve in polar coordinates as well as strong uniform bounds on the radial function and construct a \emph{radial neighborhood} around $\Gamma(t)$ of fixed width $2 \varepsilon_{\Theta}$; see Figure \ref{fig:radialneighbor}.

Let $\Theta$ be a smooth bump function in the radial direction with compact support
\begin{equation}
    \Theta: (-\varepsilon_{\Theta}, \varepsilon_{\Theta}) \to [0,1], 
    \label{eq:bumpdef}
\end{equation}
satisfying 
\begin{equation}
    \Theta(0) = 1
\end{equation}
and
\begin{equation}
    \left\| \partial_r \Theta \right\|_{L^{\infty}\left( \left( -\varepsilon_{\Theta}, \varepsilon_{\Theta} \right) \right)} \leq \frac{2}{\varepsilon_{\Theta}}.
    \label{eq:boundednessbump}
\end{equation}
Extend the bump function by $0$ outside of $(-\varepsilon_{\Theta}, \varepsilon_{\Theta})$. Using the bump function, the shifted $x$-coordinate can be extended by
\begin{equation}
    X(t, \cdot, \cdot): [0, 2 \pi] \times \left( -\varepsilon_{\Theta}, \varepsilon_{\Theta} \right) \to \R_{>0}, \, (\phi, r) \mapsto \tilde{\rho}(t,\phi) \cos(\phi) \Theta(\tilde{\rho}(\phi) - r).
\end{equation}
For notational simplicity, we omit the time variable in the subsequent lines.
The Dirichlet-energy of $X$ with respect to polar coordinates takes the form
\begin{equation}
    \int_{\mathbb{T}^2_{2L}} \left|\nabla X(\sigma) \right|^2 \,dV = \int_{[0, 2 \pi]} \int_{\tilde{\rho}(\phi) - \varepsilon_{\Theta}}^{\tilde{\rho}(\phi) + \varepsilon_{\Theta}} \left( \left( \frac{\partial X}{\partial r}  \right)^2 + \left( \frac{1}{r} \frac{\partial X}{\partial \phi}  \right)^2  \right) r \,dr d\phi,
\end{equation}
with partial derivatives
\begin{equation}
     \left. \frac{\partial X}{\partial r} \right|_{\phi, r} = \tilde{\rho}(\phi) \cos(\phi) \Theta_r(\tilde{\rho}(\phi)  - r)
\end{equation}
and
\begin{equation}
    \left. \frac{\partial X}{\partial \phi} \right|_{\phi, r} = \tilde{\rho}_{\phi}(\phi) \cos(\phi) \left( \Theta(\tilde{\rho}(\phi) - r)  - \Theta_{r}(\tilde{\rho}(\phi)  - r) \tilde{\rho}_{\phi}(\phi) \right) - \tilde{\rho}(\phi) \sin(\phi) \Theta(\tilde{\rho}(\phi) - r).
\end{equation}
Setting 
\begin{equation}
    \varepsilon_{\Theta} = \frac{R}{4},
\end{equation}
and making use of the boundedness of both $\Theta$ and its derivative, cf. \eqref{eq:bumpdef} and \eqref{eq:boundednessbump}, we can estimate
\begin{equation}
     \left| \left. \frac{\partial X}{\partial r} \right|_{\phi, r} \right|^2 \leq \tilde{\rho}^2 \frac{64}{R^2}
\end{equation}
and 
\begin{equation}
    \left| \left. \frac{\partial X}{\partial \phi} \right|_{\phi, r} \right|^2 \leq 2 \tilde{\rho}_{\phi}(\phi)^2 \left( 2 + \frac{128}{R^2} \tilde{\rho}_{\phi}(t,\phi)^2 \right) + 2 \tilde{\rho}(t,\phi)^2.
\end{equation}
In light of $E^2D \leq \varepsilon$, inserting the bounds \eqref{eq:admissannulusbound} and \eqref{eq:boundedrhoslope} above produces the desired estimate
\begin{equation}
    \int_{\R^2} \left|\nabla X(\sigma) \right|^2 \,dV \lesssim \int_{[0, 2 \pi]} \tilde{\rho}(\phi)^2 \,d\phi =  2 \left| \Omega_{\text{in}} \right|.
\end{equation}
\end{myproof}
    \section{Proof of the main theorem} \label{sec:main}

Before proving the main theorem of this article, Theorem  \ref{thm:main}, we first show that smallness of $E^2D(0)$ is preserved in time, as first observed in \cite[Lemma 5.1]{[Chugreeva2018]}:
\begin{corollary}
    There exists a universal $\varepsilon \in \R_{>0}$ such that any nearly circular curve satisfying the initial condition 
    \begin{equation}
        E^2 D(0) \leq \varepsilon,
    \end{equation}
    satisfies
    \begin{equation}
        \frac{d \left( E^2 D \right)}{dt} \leq 0
        \label{eq:EDDderivative}
    \end{equation}
    for all $t \in [0,\infty)$.
    \label{cor:stableeed}
\end{corollary}
\begin{myproof}

    A direct computation reveals
    \begin{align}
        \frac{d}{dt} (E^2D) &= 2 ED \frac{d E}{dt} + E^2 \frac{d D}{dt}
        \overset{\eqref{eq:dissipationderivative}}{\leq} - 2 ED^2 + E^2 C \left( ED^3 + D^{5/2} \right) \\
        &\leq ED^2 \left(-2 - C\left( E^2D + \left( E^2D \right)^{1/2} \right)\right)
        \overset{E^2D \leq \varepsilon}{\leq} 0, 
    \end{align}
    where $C$ is a universal constant.
\end{myproof}
\begin{remark}[Preservation of graph structure]
    Combined with Lemma \ref{lemma:globalbounds} which estimates the $L^{\infty}$-norms of $\rho - R$ and $\rho_{\phi}$, the above corollary implies that graph structure is preserved throughout the evolution if present initially, provided that $E^2D(0) \leq \varepsilon$.
    $\hfill \bigtriangleup$
    \label{rmk:graph}
\end{remark}

We can now prove the main theorem
\begin{myproof}[Proof of Theorem \ref{thm:main}] 
    The stability result \eqref{eq:stableeed} is an immediate consequence of Corollary \ref{cor:stableeed}.

    We proceed with the upper bound of the squared distance in terms of initial data, cf. \eqref{eq:squareddistancebound}:

    Let $t >0$ be arbitrary. Viewing $H$ as a function of $E$ as in the proof of \cite[Theorem 1.3]{[Chugreeva2018]} the chain rule implies
    \begin{align}
        \frac{dH}{dE} \frac{dE}{dt} \overset{\eqref{eq:dissipationderivative}}&{=} - \frac{dH}{dE} D \lesssim  H^{1/2} D^{1/2},
    \end{align}
    which can be combined with \eqref{eq:energyderivative} to obtain
    \begin{equation}
        - \frac{dH}{dE} \lesssim H^{1/2} D^{-1/2} \overset{\eqref{eq:energydissipationrelation}}{\lesssim} H^{1/2} R^{3/2} E^{-1/2}.
        \label{eq:dhde}
    \end{equation}
    Integration from $E(0)$ to $E(t)$ produces
    \begin{equation}
        H(t)^{1/2} - H(0)^{1/2} \lesssim R^{3/2} \left( E(0)^{1/2} - E^{1/2}(t) \right) \leq R^{3/2} E(0)^{1/2}
    \end{equation}
    so that indeed
    \begin{equation}
        H(t) \lesssim \max \left\{ H(0), R^3 E(0) \right\} =: \alpha_0 
    \end{equation}
    holds, which is \eqref{eq:squareddistancebound}.

    To derive the algebraic decay rate for the energy gap, cf. \eqref{eq:algebraicenergy}, we use \eqref{eq:squareddistancebound} to bring the interpolation inequality \eqref{eq:energygapinterpol} into the form
    \begin{equation}
        E(t) \lesssim \left( \alpha_0 D(t) \right)^{1/2}.
    \end{equation}
    Squaring the above inequality and inserting \eqref{eq:energyderivative} leads to
    \begin{equation}
        \frac{dE}{dt} \lesssim \alpha_{0}^{-1} E(t)^2,
    \end{equation}
    which yields \eqref{eq:algebraicenergy} upon integration. The exponential decay rate in turn follows immediately from the combination of \eqref{eq:energyderivative} and \eqref{eq:energydissipationrelation} and integration.

    Next, we consider the decay rates of the dissipation. Let $T_D > t$ be arbitrary. From \eqref{eq:energyderivative} it follows that
    \begin{equation}
        \int_{t}^{T_D} D(s) \,ds = E(t) - E(T) \leq E(t)
        \label{eq:energydiffdiss}
    \end{equation}
    must hold. Insertion of \eqref{eq:algebraicenergy} yields
    \begin{equation}
        \int_{t}^{T_D} D(s) \,ds \lesssim \frac{\alpha_{0}}{t},
        \label{eq:onealgdiss}
    \end{equation}  
    while at the same time \eqref{eq:dissipationderivative} and the smallness of $E^2D$ imply
    \begin{equation}
        \frac{dD}{dt} \lesssim D^{5/2}.
    \end{equation}
    This leads to
    \begin{equation}
        \frac{D(T_D)}{\left( 1 + (T_D-t)D(T_D)^{3/2} \right)^{2/3}} \lesssim D(t)
        \label{eq:twoalgdiss}
    \end{equation}
    through integration from $t$ to $T_D$. Combining \eqref{eq:onealgdiss} and \eqref{eq:twoalgdiss} and integrating gives
    \begin{align}
        \frac{\alpha_{0}}{t} &\gtrsim \frac{D(T_{D})\left( T_{D}-t \right)}{\left( 1+\left( T_{D}-t \right)D^{3/2}(T_{D}) \right)^{2/3}} \gtrsim \min \left\{ D(T_{D})\left( T_{D}-t \right), \left( T_{D} - t \right)^{1/3} \right\}.
    \end{align}
    A choice of $t = T_{D}/2$ reveals the dichotomy
    \begin{equation}
        D(T) \lesssim \frac{\alpha_{0}}{T_{D}^2} \quad \text{ or } \quad T_{D}^{4/3} \lesssim \alpha_{0},
    \end{equation}
    which establishes \eqref{eq:algebraicdissipation}. Inserting the exponential decay estimate \eqref{eq:exponentialenergy} instead of the algebraic one into \eqref{eq:energydiffdiss} and proceeding analogously proves \eqref{eq:exponentialdissipation}.

    The last claim to establish is the bound on the motion of the barycenter, i.e., \eqref{eq:barycenterbound}. This will be accomplished in two steps. To this end, let $T \in \R_{>0}$ be an arbitrary point in time, and assume first that $T \leq T_{1}$ holds, where $T_{1}$ is the turnover time from Theorem \ref{thm:main}.

    In light of estimate \eqref{eq:translationdissipation} on the velocity of the barycenter and the preservation of the area $ \left| \Omega_{\text{in}} \right|$, integration suffices to obtain
    \begin{equation}
    \begin{aligned}
        \left| c(T) \right| \overset{c(0) = 0}&{\leq}  \int_{0}^{T}  \left| \dot{c}(s)  \right| \,ds
         \overset{\eqref{eq:translationdissipation}}{\leq} \int_{0}^{T} \left( \frac{D(s)}{ \left|\Omega_{\text{in}} \right|} \right)^{1/2} \,ds
         \overset{\text{C-S}}{\leq} \left( \int_{0}^{T} D(s) \,ds \int_{0}^{T} \,ds \frac{1}{ \left| \Omega_{\text{in}} \right|} \right)^{1/2} \\
         \overset{\eqref{eq:energyderivative}}&{\lesssim} \left( \frac{\left( E(0) - E(T) \right) T}{ \left| \Omega_{\text{in}} \right|} \right)^{1/2}
         \overset{T \leq T_{1} \sim R^3}{\lesssim} \left( \frac{E(0) R^3}{ \left| \Omega_{\text{in}}  \right|} \right)^{1/2}
         \overset{ \left| \Omega_{\text{in}} \right| \sim R^2}{\sim} \left( E(0) R \right)^{1/2}.
    \end{aligned}
    \label{eq:baryfirst}
    \end{equation}
    In the case $T > T_{1}$, we make use of the improved (now exponential) decay rate of the energy gap to extend the integral past $T_{1}$. This is similar in spirit to the analysis carried out in the proof of \cite[Theorem 1.1]{[Otto2007]}. More concretely, defining a weight
    \begin{equation}
        w(s) = \exp\left(\frac{t}{2CR^3}\right),
        \label{eq:baryweighting}
    \end{equation}
    with $C$ being the same constant as in \eqref{eq:exponentialenergy}, we rewrite and estimate the integral
    \begin{equation}
        \left|\int_{T_{1}}^{T}  \left( \frac{D(s)}{ \left| \Omega_{\text{in}}  \right|} \right)^{1/2} \,ds  \right|
        \overset{\text{C-S}}{\leq} \left( \left| \int_{T_{1}}^{T} \frac{1}{w(s)}  \,ds \right|  \left|\int_{T_{1}}^{T}  \left( -w(s) \frac{dE}{dt}(s) \right) \,ds\right| \frac{1}{ \left| \Omega_{\text{in}} \right|} \right)^{1/2}. 
        \label{eq:baryint}
    \end{equation}
    The first integral on the right-hand-side of \eqref{eq:baryint} can be computed explicitly, and it satisfies
    \begin{equation}
        \left| \int_{T_{1}}^{T} \frac{1}{w(s)}  \,ds \right| = 2CR^3  \left| \exp\left(\frac{-T}{2CR^3}\right) - \exp\left(\frac{-T_{1}}{2CR^3}\right) \right| \leq 2CR^3,
        \label{eq:barytwo}
    \end{equation}
    while the exponential decay of $E$, cf. \eqref{eq:exponentialenergy}, helps to estimate the second integral via
    \begin{equation}
    \begin{aligned}
        \left|\int_{T_{1}}^{T}  \left( -w(s) \frac{dE}{dt}(s) \right) \,ds\right| &\leq \left|\int_{T_{1}}^{T}  \frac{dw}{dt}(s) E(s)  \,ds\right| +  \left| w(T_{1}) E(T_{1}) - w(T) E(T) \right| \\
        \overset{\eqref{eq:exponentialenergy}, \eqref{eq:baryweighting}}&{\lesssim} \left|\int_{T_{1}}^{T} \frac{1}{2CR^3}  E(0) \exp\left( -\frac{s}{2CR^3} \right)  \,ds\right| \\
          &+ \left| E(0) \left(\exp\left( -\frac{T_{1}}{2CR^3} \right)- \exp\left( -\frac{T}{2CR^3} \right) \right) \right|
        \lesssim E(0).
    \end{aligned}
    \label{eq:barythree}
    \end{equation}
    Inserting $ \left|\Omega_{\text{in}} \right| = \pi R^2$, \eqref{eq:barytwo} and \eqref{eq:barythree} into \eqref{eq:baryint} yields the bound
    \begin{equation}
        \left|\int_{T_{1}}^{T}  \left( \frac{D(s)}{ \left| \Omega_{\text{in}}  \right|} \right) \,ds  \right| \lesssim \left( E(0) R \right)^{1/2}.
    \end{equation}
     Complementing with estimate \eqref{eq:baryfirst} for times prior to $T_{1}$ establishes the claimed bound on the motion of the barycenter \eqref{eq:barycenterbound}.
\end{myproof}

    \section*{Acknowledgments}
The author was supported by the DFG-Graduiertenkolleg \emph{Energy, Entropy, and
Dissipative Dynamics (EDDy)}, project no. 320021702/GRK2326. He is indebted to his supervisors Maria G. Westdickenberg and Umberto Hryniewicz for their patience, guidance and support throughout the project. Moreover, he would like to thank the Max Planck Institute for Mathematics in the Sciences (MPI MIS) for its hospitality and especially Felix Otto for his time and valuable discussions. For further stimulating discussions, the author wishes to thank Wenhui Shi, Richard Schubert and Olaf Wittich.
    \newpage
    \bibliography{lib}
    \bibliographystyle{alpha}
    \appendix
    \section{Appendix} \label{sec:appendix}

\subsection{Homogeneous fractional Sobolev spaces}
To introduce the homogeneous fractional Sobolev spaces used throughout the document, we fix a convention for the Fourier transform of periodic functions. For $P \in \R_{>0}$ and any $2P$-periodic function $f \in L^2(\mathbb{T}^{1}_{2P})$, we introduce the Fourier transform of $f$ as
\begin{equation}
    \hat{f}:\Z \to \R, k \mapsto \hat{f}(k) = \frac{1}{\sqrt{2P} } \int_{\mathbb{T}^{1}_{2P}} f(x) e^{-\frac{2 \pi}{2P}i k x}dx.
\end{equation}
We define the Fourier transform operator
\begin{equation}
    \mathcal{F}:L^2(\mathbb{T}^{1}_{2P}) \to l^{2}(\Z), f \mapsto \left( k \mapsto  \hat{f}(k) \right), 
\end{equation}
with inverse
\begin{equation}
    \mathcal{F}^{\ast}:l^{2}(\Z) \to L^2(\mathbb{T}^{1}_{2P}), \hat{f} \mapsto \left( x \mapsto \frac{1}{\sqrt{2P} } \sum_{k \in \Z} \hat{f}(k) e^{\frac{2 \pi}{2P}i k x} \right),
\end{equation}
cf. \cite[Chapter 3.1]{[Taylor2011]}, arriving at the following definition:
\begin{definition}[Homogeneous fractional Sobolev space]
    Let $L \in \R_{> 0}$, $f \in L^{2}(\mathbb{T}^{1}_{2P})$ and $\sigma \in \R$. 
    The \textbf{homogeneous derivative operator} is defined via
    \begin{equation}
         |\partial_x|^{\sigma} f (x) = \mathcal{F}^{\ast} \left( | \frac{2 \pi }{2P} k |^{\sigma} \hat{f}(k) \right) (x)= \frac{1}{\sqrt{2P} } \sum_{k \in \Z} \left(  | \frac{2 \pi}{2P} k  |^{\sigma} \hat{f}(k) e^{\frac{2 \pi}{2P} ikx} \right),
         \label{eq:homogeneousderivative}
    \end{equation}
    for all  $x \in \mathbb{T}^{1}_{2P}$.

    The \textbf{homogeneous fractional Sobolev space} $\dot{H}^{\sigma}(\mathbb{T}_{2P}^{1})$ consists of the  $ f\in L^{2}(\mathbb{T}_{2P}^{1})$ such that 
    \begin{equation}
        \left\| f \right\|_{\dot{H}^{\sigma}(\mathbb{T}_{2P}^{1})} := \left\| | \partial_x|^{\sigma} f \right\|_{L^2(\mathbb{T}_{2P}^{1})} =  \left( \sum_{k \in \Z} | \frac{2 \pi }{2P} k|^{2\sigma} |\hat{f}(k)|^2 \right)^{1/2} < \infty
        \label{eq:homogeneoussobolevnorm}
    \end{equation}
    holds. If further $\Gamma \in \mathcal{M}$ is a closed curve of length $\lth(\Gamma)$, arc-length-parametrized by some
    \begin{equation}
        \gamma: \mathbb{T}^1_{\lth(\Gamma)} \to \mathcal{M},
    \end{equation}
    we say that a function
    \begin{equation}
        f:\Gamma \to \R
    \end{equation}
    is in $\dot{H}^{\sigma}(\Gamma)$ for some $\sigma \in  \R$, if the function
    \begin{equation}
        \tilde{f}:= f \circ \gamma: \mathbb{T}^1_{\lth(\Gamma)} \to \R
    \end{equation}
    is in $\dot{H}^{\sigma}(\mathbb{T}^{1}_{\lth(\Gamma)})$. In this case we define
    \begin{equation}
        \left\| f \right\|_{\dot{H}^{\sigma}(\Gamma)} := \left\| \tilde{f} \right\|_{\dot{H}^{\sigma}(\mathbb{T}^{1}_{\lth(\Gamma)})}. 
        \label{eq:arclengthsobolev}
    \end{equation}
\end{definition}
One of the major benefits of working with homogeneous fractional Sobolev spaces is the availability of interpolation inequalities of the following kind, where rational orders suffice for our needs:
\begin{corollary}[Interpolation inequalities]
    Let $P \in \R_{>0}$ and $\sigma \in \Q$. Let further $\alpha, \beta \in \Q$ such that $\alpha < \sigma < \beta$ holds. Then $f \in \dot{H}^{\alpha}(\mathbb{T}^{1}_{2P}) \cap \dot{H}^{\beta}(\mathbb{T}^{1}_{2P})$ implies $f \in \dot{H}^{\sigma}(\mathbb{T}^{1}_{2P})$ and 
    \begin{equation}
        \left\| f \right\|_{\dot{H}^{\sigma}(\mathbb{T}^{1}_{2P})} \leq \left\| f \right\|^{\frac{1}{p}}_{\dot{H}^{\alpha}(\mathbb{T}^{1}_{2P})} \left\| f \right\|^{\frac{1}{q}}_{\dot{H}^{\beta}(\mathbb{T}^{1}_{2P})}
        \label{eq:interpolation}
    \end{equation}
    for
    \begin{equation}
        p:= \frac{\beta - \alpha}{\beta - \sigma}, \quad \text{ and } \quad q := \frac{\beta - \alpha}{\sigma - \alpha}.
        \label{eq:hoelderexp}
    \end{equation}
    \label{cor:interpolation}
\end{corollary}
\begin{myproof}
    The estimate is an immediate consequence of the Hölder inequality paired with the observation that $p$ and $q$ from \eqref{eq:hoelderexp} satisfy
    \begin{equation}
        \frac{1}{p} + \frac{1}{q} = 1 \quad \text{ and } \quad \sigma = \frac{\alpha}{p} + \frac{\beta}{q}.
        \label{eq:hoelderone}
    \end{equation}
\end{myproof}
In addition to interpolation inequalities, homogeneous fractional Sobolev semi-norms provide Poincaré-type inequalities for fractional order derivatives:
\begin{lemma}[Poincaré-Wirtinger-type inequality]
    Let $\sigma, P \in \R_{>0}$ and let $f \in \dot{H}^{\sigma}(\mathbb{T}^{1}_{2P}) \cap L^{1}(\mathbb{T}^{1}_{2P})$. Then $f$ satisfies
    \begin{equation}
        \left\| f - \overline{f} \right\|^2_{L^2(\mathbb{T}^{1}_{2P})} \leq \left( \frac{P}{\pi} \right)^{2\sigma } \left\| f \right\|^2_{\dot{H}^{\sigma}(\mathbb{T}^{1}_{2P})},
        \label{eq:poincarewirt}
    \end{equation}
    where $\overline{f}$ denotes the integral mean of $f$.
\end{lemma}
\begin{myproof}
    The estimate follows from the Fourier representation of the homogeneous fractional Sobolev semi-norm, cf. \eqref{eq:homogeneoussobolevnorm}, and the observation that the zero-mode vanishes for a function with zero mean. 
\end{myproof}

As the definition of homogeneous Sobolev semi-norms of functions on curves is based on arc-length parametrization, cf. \eqref{eq:arclengthsobolev}, it is useful to understand how homogeneous Sobolev semi-norms of periodic functions behave under reparametrizations. This is captured in Lemma \ref{lm:conversiondotH} below:

\begin{lemma}[Conversion of homogeneous Sobolev semi-norms]
    Let $\rho \in \mathcal{R}_R$, let $\gamma$ denote the corresponding nearly circular polar parametrization and let further 
    \begin{equation}
        s(\phi):\mathbb{T}^{1}_{2 \pi} \to \mathbb{T}^{1}_{\lth(\Gamma)}
    \end{equation}
    denote the corresponding arc-length function with inverse $\phi(s)$.

    Then, for any $\sigma \in  (0,1)$, the respective $\dot{H}^{\sigma}$ semi-norms are related via
    \begin{equation}
        \left\| \left( f \circ \gamma \right) \right\|_{\dot{H}^{\sigma}([0,2 \pi])} \sim \lth(\Gamma)^{\sigma - 1/2} \left\|  f \right\|_{\dot{H}^{\sigma}(\Gamma)} \sim R^{\sigma - 1/2} \left\|  f \right\|_{\dot{H}^{\sigma}(\Gamma)}.
        \label{eq:arcpolarcircle}
    \end{equation}

    The analogous estimates for $\sigma \in (-1,0)$ read
    \begin{equation}
        \left\| \left( f \circ \gamma \right) \right\|_{\dot{H}^{\sigma}([0,2 \pi])} \sim \lth(\Gamma)^{\sigma + 1/2} \left\|  f \right\|_{\dot{H}^{\sigma}(\Gamma)} \sim R^{\sigma + 1/2} \left\|  f \right\|_{\dot{H}^{\sigma}(\Gamma)}.
        \label{eq:neghomogeneousnormconversion}
    \end{equation}
    \label{lm:conversiondotH}
\end{lemma}
    The behavior under reparametrization can be analyzed based on the equivalence between Fourier-based homogeneous Sobolev norms and the Solobodecki-semi-norms, cf. \cite[Prop. 1.3]{[Benyi2012]}. 
\begin{remark}
    For the sake of simplicity, we simply write \emph{norm} instead of \emph{semi-norm} in the main text.$\hfill \bigtriangleup$
\end{remark}

\subsection{Expressions for the Dirichlet energy} \label{app:bulkparams}

We state the explicit expressions for the Dirichlet energy of a function $f \in \dot{H}^{1}(\mathcal{M})$ under pullbacks by either bulk parametrization. For the polar parametrization $P$, it is given by
\begin{equation}
    \int_{\mathcal{M}}  \left| \nabla f  \right|^2 dV =\int_{\mathcal{P}_{\text{in}} \cup \mathcal{P}_{\text{out}}} \left( \left( \frac{\partial f}{\partial r} \right)^2 + \left( \frac{1}{r} \frac{\partial f}{\partial \phi}  \right)^2  \right) rdr d\phi. \label{eq:dirichletcircle}
\end{equation}
The expression for the flattening parametrization $F$ reads
\begin{equation}
    \int_{\mathcal{M}}  \left| \nabla f  \right|^2 dV = \int_{\mathcal{F}_{\text{in}} \cup \mathcal{F}_{\text{out}}} \left(\frac{c^2 + \left( ra \right)^2}{\left( a+b \right)ra}  \right) \left( \frac{\partial f}{\partial r} \right)^2 -2\frac{c}{ra}\frac{\partial f}{\partial r} \frac{\partial f}{\partial \phi}  + \left( \frac{a+b}{ra} \right) \left( \frac{\partial f}{\partial \phi}  \right)^2   dr d\phi\label{eq:dirichlethanzawa}
\end{equation}
where $a,b$ and $c$ are defined as
\begin{align}
    a(r,\phi) = \left(\frac{\rho(\phi)}{R} - 1\right) \beta(r) + 1, \quad
    b(r,\phi) = r \left( \frac{\rho(\phi)}{R} - 1 \right)\partial_r \beta(r), \quad
    c(r,\phi) = r \frac{\rho_{\phi}(\phi)}{R} \beta(r).
\end{align}

\subsection{Proof of the trace estimate:} \label{app:traceestimate}
\begin{myproof}[Proof of Theorem \ref{thm:traceestimate}]

    The proof of the trace estimate proceeds in three steps: We first establish the equivalence between the Dirichlet energy of a harmonic problem solved on $\Omega_{\text{in}}$ and $\Omega_{\text{out}}$ and the Dirichlet energy between a harmonic problem solved on a circle of radius $R$ and its complement for equal boundary data  (Lemma \ref{lemma:tracecircle}). It is then shown that the Dirichlet energy on the circular domain and its complement are equivalent to the $\dot{H}^{1/2}$-norm on the circular boundary (Lemma \ref{lemma:traceestimatetorus}). Lastly, we use the equivalence between the $\dot{H}^{1/2}$-norm on the circular boundary and the $\dot{H}^{1/2}$-norm on the interface from Lemma \ref{lm:conversiondotH}.

    More explicitly, we show
    \begin{align}
        \int_{\mathcal{M} \setminus \Gamma} \left\| \nabla f \right\|^2 dV
        \overset{\text{Lemma } \ref{lemma:tracecircle}}{\sim} \int_{\mathcal{M} \setminus \partial B_{R}(0)} \left\| \nabla \tilde{f} \right\|^2 dV 
        \overset{\text{Lemma } \ref{lemma:traceestimatetorus}}{\sim} \left\|  \left. \tilde{f} \right|_{\partial B_R(0)} \right\|_{\dot{H}^{1/2}([0, 2 \pi])} 
        \overset{\text{Lemma } \ref{lm:conversiondotH}}{\sim} \left\|  \left. \tilde{f} \right|_{\Gamma} \right\|_{\dot{H}^{1/2}(\Gamma)}.
    \end{align}

    Notice that, although Lemma \ref{lemma:traceestimatetorus} only provides the necessary inequality in one direction for $\Omega_{\text{out}}$, equivalence for $\Omega_{\text{in}}$ is sufficient to obtain an overall equivalence.
\end{myproof}

We begin with Lemma \ref{lemma:tracecircle}:
\begin{lemma}
    Let $\rho \in \mathcal{R}_{R}$ be the radial function of the polar parametrization $\gamma$ which parametrizes a nearly circular curve $\Gamma$ and let $f \in H^{1}(\Omega_{\text{in}} \cup \Omega_{\text{out}})$ satisfy
    \begin{equation}
        \begin{aligned}
            \Delta f(t,\cdot) &= 0 \quad &&\text{ on } \Omega_{\text{in}} \cup \Omega_{\text{out}}\\
            f(t,\cdot) &= g(t,\cdot) \quad &&\text{ on } \Gamma . \label{eq:generalharmonicproblem}
        \end{aligned}
    \end{equation}
    for $g \in \dot{H}^{1/2}(\Gamma)$. Let further $\tilde{f} \in H^{1}(B_{R}(0) \cup \mathbb{T}^2_{2L} \setminus \overline{B_R(0)}) $ solve the harmonic problem
    \begin{equation}
        \begin{aligned}
            \Delta \tilde{f}(t,\cdot) &= 0 \quad &&\text{ on } B_{R}(0) \cup \mathbb{T}^2_{2L} \setminus \overline{B_R(0)}\\
            \tilde{f}(t,\cdot) &= \tilde{g}(t,\cdot) \quad &&\text{ on } \partial B_R(0), \label{eq:circleharmonicproblem}
        \end{aligned}
    \end{equation}
    where $\tilde{g} = g \circ P(R, \cdot) \circ \gamma^{-1}$, cf. \eqref{eq:polarparamcurve} and \eqref{eq:polarparam} is the Dirichlet data of the former problem pulled back to the circle of radius $R$.

    Then, there holds
    \begin{equation}
          \int_{\mathcal{M}}  \left| \nabla f \right|^2 dV \sim \int_{\mathcal{M}}  \left| \nabla \tilde{f} \right|^2 dV.
    \end{equation}
    \label{lemma:tracecircle}
\end{lemma}
\begin{myproof}
    In a first step, we pull back both problems to the domain 
    \begin{equation}
        \mathcal{F} := \mathcal{F}_{\text{in}} \cup \mathcal{F}_{\text{out}} = \bigsqcup_{\phi \in [0, 2 \pi]} (0, R_{\text{max}, 2L}(\phi)) \times \{\phi\}. 
    \end{equation}
    For problem \eqref{eq:generalharmonicproblem}, this is achieved via the parametrizations \eqref{eq:interiorflat} and \eqref{eq:twopolarparambulkout}, while problem \eqref{eq:circleharmonicproblem} is pulled back to this domain by the polar parametrization \eqref{eq:polarparam}. After identifying $f$ and $\tilde{f}$ with their pullbacks, both solutions have the same domain and the exact same Dirichlet boundary data given by $g \circ \gamma^{-1}$. Thus, we can make use of Dirichlet's principle to estimate the Dirichlet energies (which differ due to the differences in pullback metrics) against each other. These estimates follow immediately from the explicit formulas \eqref{eq:dirichletcircle} and \eqref{eq:dirichlethanzawa}.
\end{myproof}

To establish the trace estimate, it remains to prove the following lemma. For notational simplicity, we assume without loss of generality by scaling $R=1$.
\begin{lemma}[Trace estimate on the disk and its complement]
    Let $L \in \R_{\geq \sqrt{2}}$ and  $g \in \dot{H}^{1/2}(S^{1})$ and suppose that $u \in \dot{H}^{1}(\Omega)$ for
    \begin{equation}
        \Omega \in \{ B_{1}(0), \R^2 \setminus B_{1}(0), \mathbb{T}^2_{2L} \setminus B_{1}(0) \} 
    \end{equation}
    solves the following boundary value problem:
    \begin{equation}
    \begin{aligned}
        \Delta u &= 0 \quad \text{ in } \Omega \\
               u &= g \quad \text{ on } \partial \Omega = S^{1}. 
    \end{aligned}
    \label{eq:dirichletproblemtorus}
    \end{equation}
    Then the estimate
    \begin{equation}
        \int_{\Omega}  \left| \nabla u \right|^2 dV \lesssim \left\| g \right\|^2_{\dot{H}^{1/2}(S^{1})}
        \label{eq:traceinequ}
    \end{equation}
    holds. For $\Omega = B_{1}(0)$ and $\Omega = \R^2 \setminus B_{1}(0)$ equality holds:
    \begin{equation}
        \int_{\Omega}  \left| \nabla u \right|^2 dV = \ \left\| g \right\|^2_{\dot{H}^{1/2}(S^{1})}. \label{eq:traceequalcircle}
    \end{equation}
    \label{lemma:traceestimatetorus}
\end{lemma}
\begin{myproof}
        For $\Omega \in \{ B_{1}(0), \R^2 \setminus B_{1}(0)\}$ identity \eqref{eq:traceequalcircle} can be obtained by explicitly solving the Laplace problem, computing the corresponding Dirichlet energy and relating that to the $\dot{H}^{1/2}$-norm of the Dirichlet data.

        For the case $\Omega = \mathbb{T}^2_{2L} \setminus B_{1}(0)$, one explicitly solves the Laplace problem on the annulus
        \begin{equation}
            A_{L} := \left\{(x,y) \in \R^2 \mid 1 \leq \sqrt{x^2 + y^2} \leq L \right\},
        \end{equation}
        where the Dirichlet data is given by $g$ on the inner boundary and by its integral mean $\overline{g}$ on the outer boundary. One can check that this solution satisfies \eqref{eq:traceinequ} (for the domain $A_{L}$). By constant extension to $M_{2L} \setminus \left(A_L \cup B_{1}(0)\right)$, one creates a competitor for the solution of the original Laplace problem. By the Dirichlet principle, the Dirichlet energy of the solution of the original problem must be bounded from above by that of the annular competitor, which proves the claim.
\end{myproof}

\subsection{Periodic potential theory in two dimensions} \label{app:potential}

We collect a few technical results on periodic potential theory in two dimensions in this section with a focus on the (quasi-)elliptic \emph{Weierstrass sigma function}.

\subsubsection{The Weierstrass functions}
Let $\omega_{1}, \omega_{3} \in \C$ be linearly independent over $\R$ and consider the lattice 
\begin{equation}
    \mathcal{G} := 2\Z \omega_{1} + i 2\Z \omega_{3} \subset \C.
\end{equation}
Then any meromorphic function $f:\C \to \C$ that is doubly periodic with respect to $\mathcal{G}$ is called an \textbf{elliptic function}. If $f$ is not periodic for any submultiple of $\omega_{1}$ or $\omega_{3}$, $\omega_{1}$ and $\omega_{3}$ are called the \textbf{fundamental pair of periods} for $f$. It is known that the elliptic functions form a subfield of the meromorphic functions, cf. \cite[Proposition 1.2.1]{[Krieg2007]}, which is generated by rational functions in the \textbf{Weierstrass} $\wp$-\textbf{function}
\begin{equation}
    \wp:\C \to \C, z \mapsto \frac{1}{z^2} + \sum_{m,n \in \Z \setminus \{0\} } \frac{1}{\left( z - z_{m,n} \right)^2 } -  \frac{1}{\left( z_{m,n}\right)^2},
\end{equation}
and its derivative $\wp'$, cf. \cite[Theorem 1.2.4]{[Krieg2007]}.  Here, the $z_{m,n}$ are the lattice points
\begin{equation}
    z_{m,n} = 2m \omega_{1} + 2n\omega_{3} \in \mathcal{G}
\end{equation}
for $m, n \in \Z$.
The Weierstrass $\wp$-function is the negative derivative of the (quasi-)elliptic \textbf{Weierstrass zeta-function} $\zeta$, i.e,
\begin{equation}
    \wp(z) = - \frac{d}{dz} \zeta(z) \quad \forall z \in \C \setminus \mathcal{G},
\end{equation}
which admits the series representation
\begin{equation}
    \zeta(z) = \frac{1}{z} + \sum_{m,n \in \Z \setminus \{0\} } \frac{1}{z - z_{m,n} }+ \frac{1}{z_{m,n}} + \frac{z}{\left( z_{m,n} \right)^2},
    \label{eq:zetafunction}
\end{equation}
which converges absolutely and uniformly for every $z \in \C \setminus \mathcal{G}$, cf. \cite[Chapter 1.6.1]{[Krieg2007]}. This meromorphic function with a single order-one pole at the origin is in turn the logarithmic derivative of the \textbf{Weierstrass sigma-function}, i.e.,
\begin{equation}
    \zeta(z) = \frac{d}{dz} \left( \log(\sigma (z)) \right) \quad \forall z \in \C \setminus \mathcal{G}.
    \label{eq:zetasigma}
\end{equation}
The Weierstrass sigma-function satisfies the product representation
\begin{equation}
    \sigma(z) = z \prod_{m,n \in \Z \setminus \{0\}} \left( 1 - \frac{z}{z_{m,n}} \right) \exp\left(\frac{z}{z_{m,n}} + \frac{1}{2} \left( \frac{z}{z_{m,n}} \right)^2\right),
    \label{eq:sigmaproduct}
\end{equation}
which implies the series representation
\begin{equation}
    \log(\sigma (z)) = \log(z) + \sum_{m,n \in \Z \setminus \{0\}} \log\left( 1 - \frac{z}{z_{m,n}} \right) + \frac{z}{z_{m,n}} + \frac{1}{2} \left( \frac{z}{z_{m,n}} \right)^2
    \label{eq:logsigmaseries}
\end{equation}
for its (principal branch) logarithm, cf. \cite[Chapter 8.4]{[Lawden1989]}.

The quasi-periodicity of $\zeta$  is captured by the quantities $\eta_{1}, \eta_{3} \in \C$ defined as
\begin{equation}
    \zeta(z + 2\omega_{j}) -  \zeta(z) =: 2\eta_{j},
\end{equation}
so that
\begin{equation}
    \log(\sigma(z + 2\omega_{j})) - \log(\sigma(z))= i \pi + 2\eta_{j}\left( z + \omega_j \right)
    \label{eq:quasilogsigma}
\end{equation}
holds with $j \in \{1,3\}$, cf. \cite{[Hasimoto2008]}.
The \emph{Legendre relation}
\begin{equation}
    \eta_1 \omega_{3} - \eta_3 \omega_{1} = \frac{i \pi}{2},
    \label{eq:legendre}
\end{equation}
cf. \cite[Chapter 1.6.1]{[Krieg2007]} for a proof, relates these quantities to the generators of the underlying lattice.

\begin{figure}[h]
    \centering
    \includegraphics[width=0.7\textwidth]{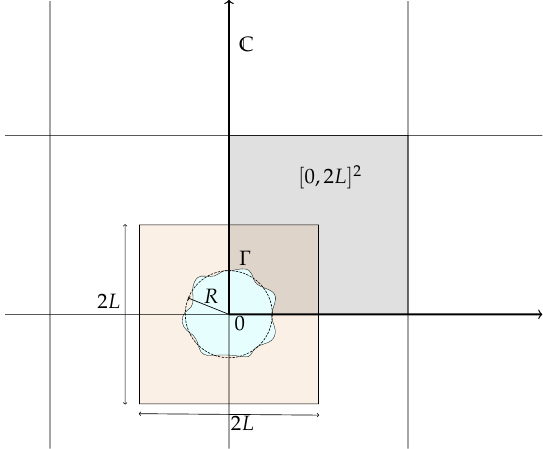}
    \caption{The orange rectangle $[-L, L]^2$ represents the torus $\mathbb{T}^2_{2L}$, while the grey rectangle represents the fundamental cell $[0, 2L]^{2}$ for the lattice of the elliptic functions.}
    \label{fig:toruscell}
\end{figure}

In this article, we restrict ourselves to the square lattice generated by $\omega_{1} = L$ and $\omega_{3} = iL$, with $L$ being the length scale of the torus. Moreover, we shift the lattice in such a way that its origin coincides with the origin of the torus $\mathbb{T}^2_{2L}$, cf. Figure \ref{fig:toruscell}.

Based on the logarithm of the sigma function, we can build a periodic analog of the fundamental solution of the Laplace equation on $\mathbb{T}^2_{2L}$ as the following lemma shows:

\begin{lemma}[Periodic analog of fundamental solution]
    The real-valued function
    \begin{equation}
        \Lambda(z): \C \to \R, \, z \mapsto \frac{1}{2} \left( \log(\sigma(z)) + \overline{\log(\sigma(z))} \right) - \frac{\pi}{4L^2}  \left| z \right|^2
        \label{eq:periodicpotential}
    \end{equation}
    is a periodic analog of the fundamental solution on the torus of edge length $2L$, i.e., it satisfies (in a distributional sense)
    \begin{equation}
        \Delta \Lambda(z) = 2 (\partial_{z} \partial_{\overline{z}} + \partial_{\overline{z}} \partial_z) \Lambda(z) = 4\pi \left( \sum_{m,n \in \Z} \delta(z - z_{m,n}) - \frac{1}{4L^2} \right),
        \label{eq:periodicpotentialdelta}
    \end{equation}
    where $\delta$ is the Dirac delta distribution.

    Moreover, $\Lambda(z)$ is doubly periodic with period $2L$ and singularities only at the lattice points. It thus descends to a (single-valued) function on the torus $\mathbb{T}^2_{2L}$ with a single singularity at the origin.
    \label{lemma:periodicpotential}
\end{lemma}
\begin{myproof}
    The following proof is based on the article \cite{[Hasimoto2008]}.

    For the Wirtinger operators
    \begin{equation}
        \partial_{z} = \frac{1}{2} \left( \partial_{x} - i \partial_{y} \right) \quad \text{and} \quad \partial_{\overline{z}} = \frac{1}{2} \left( \partial_{x} + i \partial_{y} \right),
    \end{equation}
    one can compute the relation
    \begin{equation}
        \Delta f = 4 \partial_{z} \partial_{\overline{z}} f = 4 \partial_{\overline{z}} \partial_{z} f.
    \end{equation}

    It is immediate that
    \begin{equation}
        \Delta \left| z \right|^2  = \left( 4 \partial_{z} \partial_{\overline{z}} \right) \left( z \overline{z} \right) = 4 
        \label{eq:modulusz}
    \end{equation}
    holds, so that only the singular part of the fundamental solution is left to be determined. To this end, we make use of the \emph{Cauchy-Pompeiu}-formula
    \begin{equation}
        f(u) = \frac{1}{2 \pi i} \left( \int_{\partial U} \frac{f(z)}{z - u} dz + \int_{U} \frac{\partial_{\overline{z}} f(z)}{z - u} dz \wedge d\overline{z} \right),
    \end{equation}
    with $u \in U$ for some open domain $U \subset \C$ with $C^{1}$-boundary and $f:\overline{U}\to \C$ continuous with bounded partial derivatives in $U$, cf. \cite[Theorem 4.3]{[Curry2024]}.

    Together with the computation
    \begin{equation}
        dz \wedge d\overline{z} = d\left( x + iy \right) \wedge d\left( x-iy \right) = -2i dx \wedge dy = -2i dV,
    \end{equation}
    where $dV$ is the volume form on the flat torus, this formula immediately establishes that
    \begin{equation}
        \partial_{\overline{z}} \frac{1}{z - u} = \pi \delta(z-u)
    \end{equation}
    holds in a distributional sense. Holomorphicity implies
    \begin{equation}
        \partial_{\overline{z}} \frac{1}{u} = 0 = \partial_{\overline{z}} \frac{z}{\left( u \right)^2}, 
    \end{equation}
    so that the above distributional equality can be transformed into
    \begin{equation}
        \partial_{\overline{z}} \left( \frac{1}{z - z_{m,n}} + \frac{1}{z_{m,n}} + \frac{z}{\left( z_{m,n} \right)^2} \right)= \pi \delta(z-z_{m,n}).
    \end{equation}
    A comparison with \eqref{eq:zetafunction} paired with the absolute convergence of the series indicates
    \begin{equation}
        \partial_{\overline{z}} \zeta(z) = \pi \sum_{m,n \in \Z} \delta(z - z_{m,n}),
    \end{equation}
    which, in light of \eqref{eq:zetasigma}, leads to
    \begin{equation}
       \frac{1}{4} \Delta \log(\sigma(z)) = \partial_{z} \partial_{\overline{z}} \log(\sigma(z)) = \pi \sum_{m,n \in \Z} \delta(z - z_{m,n}).
        \label{eq:zetadeltaone}
    \end{equation}
    The analogous relation
    \begin{equation}
        \frac{1}{4} \Delta \overline{\log(\sigma(z))} = \partial_{\overline{z}} \partial_{z} \overline{\log(\sigma(z))} = \pi \sum_{m,n \in \Z} \delta(z - z_{m,n}).
        \label{eq:zetadeltatwo}
    \end{equation}
    follows from symmetry in $\overline{z}$ and $z$. For the function
    \begin{equation}
        \Lambda: \C \to \R, \, z\mapsto  \frac{1}{2} \left( \log(\sigma(z)) + \overline{\log(\sigma(z))} \right) - \frac{\pi}{4L^2}  \left| z \right|^2,
        \label{eq:periodicpotentialtwo}
    \end{equation}
    equations \eqref{eq:modulusz}, \eqref{eq:zetadeltaone}, and \eqref{eq:zetadeltatwo} imply the defining equation \eqref{eq:periodicpotentialdelta}. 

    Next, we observe that the product representation \eqref{eq:sigmaproduct} of $\sigma$ implies that the singularities of $\Lambda(z)$ are restricted to the lattice $\mathcal{G}$.

    To conclude the proof, we need to show that $\Lambda$ is indeed doubly periodic in $2\omega_{1}=2L$ and $2\omega_{3} = i 2L$. The square structure of the lattice together with the Legendre relation \eqref{eq:legendre} implies
    \begin{equation}
        \eta_j \omega_j = \frac{\pi}{4}
    \end{equation}
    for $j \in \{1,3\}$, so that
    \begin{equation}
        \eta_{1} = \frac{\pi}{4L} \quad \text{and} \quad \eta_{3} = \frac{i \pi}{4L}
    \end{equation}
    hold for the lattice at hand, cf. \cite[(30),(31)]{[Hasimoto2008]}.

    These values for $\eta_{1}$ and $\eta_{3}$ together with \eqref{eq:quasilogsigma} can be used to compute the identities
    \begin{equation}
        \Lambda(z + 2\omega_j) - \Lambda(z) = 0
    \end{equation}
    for $z \in \C \setminus \mathcal{G}$ and $j \in \{1,3\}$, so that $\Lambda$ is indeed doubly periodic with periods $2L$ and $2iL$.

\end{myproof}
\begin{remark}
    By adding its complex conjugate to $\log(\sigma(z))$, we eliminate the jumps over the branch cuts of the logarithm. As a result, $\Lambda$ is holomorphic on $\C \setminus \mathcal{G}$. $\hfill \bigtriangleup$
\end{remark}

The following representation formula is an immediate consequence of the series representation \eqref{eq:logsigmaseries}, the invariance of the lattice under conjugation and the Laurent series expansion of the logarithm
\begin{equation}
    \log(1 - z) = - \sum_{k=1}^{\infty} \frac{z^{k}}{k!} 
\end{equation}
for $ \left| z \right| < 1$.

\begin{corollary}
    The periodic fundamental solution $\Lambda$ satisfies the following, absolutely convergent series expansion
    \begin{align}
        \Lambda(z) &= \log( \left| z \right|) \\
        &+ \sum_{m,n \in \Z \setminus \{0\}} \log\left( 1 - \frac{z}{z_{m,n}} \right) + \frac{z}{z_{m,n}} + \frac{1}{2} \left( \frac{z}{z_{m,n}} \right)^2 \\
        &+ \sum_{m,n \in \Z \setminus \{0\}} \log\left( 1 - \overline{\frac{z}{z_{m,n}}} \right) + \overline{\frac{z}{z_{m,n}}} + \frac{1}{2} \left( \overline{\frac{z}{z_{m,n}}} \right)^2 \\
        &- \frac{\pi}{4L^2}  \left| z \right|^2.
        \label{eq:potentialseriesone}
    \end{align}

    For $ \left| z \right| < 2L$, it takes the form
    \begin{equation}
        \Lambda(z) = \log( \left| z \right|) + \sum_{m,n \in \Z \setminus \{0\}} \sum_{k=3}^{\infty} - \frac{1}{k!} \left( \left( \frac{z}{z_{m,n}} \right)^{k} + \left( \overline{\frac{z}{z_{m,n}}} \right)^{k} \right) - \frac{\pi}{4L^2}  \left| z \right|^2.
        \label{eq:potentialseriestwo}
    \end{equation}
    \label{cor:sigmaexpansion}
\end{corollary}

The periodic fundamental solution $\Lambda$ can now be used to construct the \textbf{single-layer potential}: Let $\phi \in C^{0,\alpha}(\Gamma)$ be given. Then the single-layer potential $S[\phi]$ is defined as
\begin{equation}
    S[\phi](\vek{x}) = \int_{\Gamma} \Lambda(\vek{x} - \vek{y}) \phi(\vek{y}) d\sigma(\vek{y}) \quad \forall \vek{x} \in \mathbb{T}^2_{2L}.
    \label{eq:singlelayerpotential}
\end{equation}
The periodic single-layer potential satisfies all of the properties of the classical single-layer potential, and the reader is referred to \cite[Chapter 12]{[DallaRiva2021]} for a comprehensive overview.

The following proposition, which is adopted from \cite[Theorems 12.23, 12.25 ]{[DallaRiva2021]}, shows how single-layer potentials may be used to solve Laplace problems with Neumann boundary conditions:

\begin{prop}[Neumann-Laplace problems and single-layer potentials]
    Let $g \in C^{0,\alpha}(\Gamma)$ satisfy
    \begin{equation}
        \int_{\Gamma} g(s) ds = 0.
    \end{equation}
    Then the set of functions $u_{\text{out}} \in C^{1,\alpha}(\overline{\Omega_{\text{out}}})$ that satisfies the exterior Neumann boundary value problem
    \begin{alignat}{2}
        \Delta u_{\text{out}} &= 0 \quad &&\text{ in } \Omega_{\text{out}} \\
        \partial_{\nu} u_{\text{out}} &= g \quad &&\text{ on } \Gamma
    \end{alignat}
    is equal to
    \begin{equation}
        \left\{ S[\phi_{\text{out}}] + c: c \in \R \right\}, 
    \end{equation}
    where $\phi_{\text{out}} \in C^{0, \alpha}(\Gamma)$ is the unique function satisfying
    \begin{equation}
    \begin{aligned}
        \int_{\Gamma} \phi(s) ds &= 0, \\
        \left( \frac{1}{2} \text{Id} + W_{2L, \Gamma}^{t}\right) \phi_{\text{out}} &= g. \label{eq:jumpone}
    \end{aligned}
    \end{equation}
    Here $W_{2L, \Gamma}^{t}$ is the integral operator
    \begin{equation}
        W_{2L, \Gamma}^{t}: C^{0,\alpha}(\Gamma) \to  C^{0,\alpha}(\Gamma), \, \phi \mapsto \left( \vek{x} \mapsto \int_{\Gamma} \phi(\vek{y}) n(\vek{x}) \cdot \nabla \Lambda(\vek{x} - \vek{y}) ds(\vek{y}) \right).
    \end{equation}
    Analogously, the set of functions $u_{\text{in}} \in C^{1,\alpha}(\overline{\Omega_{\text{in}}})$ that satisfies the interior Neumann boundary value problem
    \begin{alignat}{2}
        \Delta u_{\text{in}} &= 0 \quad &&\text{ in } \Omega_{\text{in}} \\
        \partial_{\nu} u_{\text{in}} &= g \quad &&\text{ on } \Gamma
    \end{alignat}
    is equal to
    \begin{equation}
        \left\{ S[\phi_{\text{in}}] + c: c \in \R \right\}, 
    \end{equation}
    where $\phi_{\text{in}} \in C^{0, \alpha}(\Gamma)$ is the unique function satisfying
    \begin{align}
        \int_{\Gamma} \phi_{\text{in}}(s) ds &= 0, \\
        \left( -\frac{1}{2} \text{Id} + W_{2L, \Gamma}^{t}\right) \phi_{\text{in}} &= g, \label{eq:jumptwo}
    \end{align}
    with $W_{2L, \Gamma}^{t}$ defined as above. 
    \label{prop:existencesinglelayer}
\end{prop}
\begin{myproof}[Proof Sketch]
    The proposition is based on similar statements for classical single-layer potential theory, cf. \cite[Chapter 4]{[DallaRiva2021]} or \cite[Chapter 3]{[Folland1995]}. By decomposing the periodic single-layer potential into the classical part and into a remainder term, which is analytic outside the lattice, one can make use of the classical theory. 

    Notice that the periodic fundamental solution presented in \cite{[DallaRiva2021]} differs from the one presented in Corollary \ref{lemma:periodicpotential}. In \cite{[Hasimoto2008]}, however, it is shown that the difference is merely a constant and thus negligible. In \cite[Chapter 7.11]{[Taylor2]}, the author presents the result for a broader class of fundamental solutions on compact manifolds without any explicit formulae.
\end{myproof}

This last corollary allows us to relate the single-layer potential induced by the normal velocity of the Mullins-Sekerka evolution to the harmonic extensions of curvature:

\begin{corollary}[Periodic single-layer potential]
    Let $u_{i} \in C^{1,\alpha}(\Omega_{i})$ solve the Laplace problems
    \begin{alignat}{2}
        \Delta u_{i} &= 0 \quad &&\text{ in } \Omega_{i} \\
        \partial_{n} u_{i} &= f \quad &&\text{ on } \Gamma
    \end{alignat}
    for $i \in \{\text{in}, \text{out}\}$ with Dirichlet data $f \in C^{1,\alpha}(\Gamma)$.
    Then the single-layer potential $v := S[\phi]$ with density
    \begin{equation}
        \phi(\vek{x}) = \frac{\partial u_{\text{out}}}{\partial n} - \frac{\partial u_{\text{in}}}{\partial n}
        \label{eq:jumpdensity}
    \end{equation}
    is equal up to a constant to the composite function
    \begin{equation}
    \begin{aligned}
        u(\vek{x}) := \begin{cases}
            u_{\text{in}}(\vek{x}) &\quad \text{for } \vek{x} \in \Omega_{\text{in}}, \\
            u_{\text{out}}(\vek{x}) &\quad \text{for } \vek{x} \in \Omega_{\text{out}}, \\
            u(\vek{x})  &\quad \text{else},
        \end{cases}
    \end{aligned}
    \end{equation}
    where $u$ is defined as the limit
    \begin{equation}
        u(\vek{x}) := \lim_{\Omega_{\text{out}} \ni \vek{x}_{\text{out}}\to \vek{x} } u_{\text{out}}(\vek{x}_{out}) = \lim_{\Omega_{\text{in}} \ni \vek{x}_{\text{in}}\to \vek{x}} u_{\text{in}}(\vek{x}_{\text{in}}).
    \end{equation}
    \label{lemma:uniquepotential}
\end{corollary}
\begin{myproof}
    By Proposition \ref{prop:existencesinglelayer}, there exists a Hölder continuous density $\psi$ such that
    the single-layer potential
    \begin{equation}
        w := S[\psi]
    \end{equation}
    solves the Neumann Laplace problems
    \begin{alignat}{2}
        \Delta w_{i} &= 0 \quad &&\text{ in } \Omega_{i} \\
        \frac{\partial w_{i}}{\partial n} &=  \frac{\partial w_{i}}{\partial n} \quad &&\text{ on } \Gamma
    \end{alignat}
    for $i \in  \{ \text{in}, \text{out}\}$. In fact, this density can be obtained by applying the Dirichlet-to-Neumann operator to the density of the double-layer potential satisfying the corresponding Dirichlet problems, cf. \cite[(11.35)]{[Taylor2]}.

    The jump formulas \eqref{eq:jumpone} and \eqref{eq:jumptwo} together with the Neumann boundary conditions and definition \ref{eq:jumpdensity} imply 
    \begin{equation}
        \phi(\vek{x}) \overset{\eqref{eq:jumpdensity}}{=} \frac{\partial u_{\text{out}}}{\partial n}(\vek{x}) - \frac{\partial u_{\text{in}}}{\partial n}(\vek{x}) \overset{\eqref{eq:jumpone} - \eqref{eq:jumptwo}} = \psi(\vek{x})
    \end{equation}
    for all $\vek{x} \in \Gamma$.

    In combination with the uniqueness result from Proposition \ref{prop:existencesinglelayer}, this implies that $v$ and $w$  and thus $v$ and $u$ must be equal up to a constant.
\end{myproof}

\end{document}